\documentclass[a4paper,10pt]{article}
\usepackage[francais,english]{babel}
\usepackage{amsmath, amssymb, amsthm, amscd}
\usepackage{amsrefs}
\usepackage[dvips]{graphicx}
\usepackage[active]{srcltx}
\usepackage{pstricks}
\usepackage{pst-all}

  \usepackage{geometry}
\newcommand{\CQFD}{\hfill{$\Box$}}
\newcommand{\ds}{\displaystyle}
\newcommand{\bs}{\boldsymbol}

 \newtheorem{prop}{Proposition}
 \newtheorem{lemme}{Lemma}
 \newtheorem{theo}{Theorem}
 \theoremstyle{definition}
 \newtheorem{defi}{Definition}
 \newtheorem*{dem}{Proof}

 \newtheorem{cor}{Corollary}[section]
 
 \newtheorem{rqs}{Remark}
\newtheorem*{comment}{\textit{\textbf{Comment}}}
 \theoremstyle{remark}
 \newtheorem{ex}{Example}

\newcommand{\Nr}{\mathbb{N}^r \setminus \left\{\textbf{0}\right\}}
\newcommand{\Nn}{\mathbb{N}^n \setminus \left\{\textbf{0}\right\}}
\newcommand{\summ}{\displaystyle\sum}

\begin{document}


\title{Extension of Estermann's theorem to  Euler products associated to a multivariate polynomial}

\author{Ludovic Delabarre\footnote{Universit\'{e} de Saint-Etienne, Facult\'{e} des Sciences, Laboratoire de Math{\'e}matiques LaMuse, 23 rue du Docteur Paul Michelon 42023 Sant-Etienne cedex 2, FRANCE. Email: ludovic.delabarre(a)univ-st-etienne.fr 
\textbf{Keywords:} multivariables Euler products, meromorphic continuation, natural boundary, cyclotomic polynomial, rational point on a toric variety
\textbf{Classif. math.}:11M32 11M41 32D15 11N99 14G05}}


\maketitle

\abstract{
 Given a multivariate polynomial $h\left(X_1,\dots,X_n\right)$ with integral coefficients verifying an hypothesis of analytic regularity (and satisfying $h(\textbf{0})=1$), we determine the maximal domain of meromorphy of the  Euler product  $\prod_{p \ \textrm{prime}}h\left(p^{-s_1},\dots,p^{-s_n}\right)$ and the natural boundary is  precisely described when it exists.  In this way we extend  a well known result  for one variable polynomials due to Estermann from 1928.
As an application, we  calculate  the natural boundary of the multivariate  Euler products associated to a family of toric varieties.
}




\tableofcontents


\textbf{ACKNOWLEDGEMENTS:}

In this article, the main results of the author's thesis are presented. The author expresses his thanks to his advisor Driss Essouabri for his many helpful remarks and his careful reading of this paper.
 He also wants to thank Régis de la Bretèche and Ben Lichtin for his precious advises and his many relevant remarks concerning the writing of this paper.
  In addition he expresses his gratitude towards the University of Caen, in particular the laboratory of mathematics LMNO where he began his thesis.

\section{Introduction.}

A classic result of Estermann (\cite {est}) from 1928 characterized precisely when an Euler product 
$Z(s) = \prod_p \ h(p^{-s}),$ 
determined by a polynomial $h(X) \in \mathbb{Z}[X]$ with $h(0) = 1,$ admits a meromorphic extension to $\mathbb{C}.$ In addition, Estermann showed that if this property is not satisfied, then the Euler product has a natural boundary as a meromorphic function which he identified exactly.  

This article extends Estermann's theorem to 
all Euler products $Z(\mathbf{s}) = h(p^{-s_1},\dots, p^{-s_n})$ determined by any polynomial $h(X_1,\dots, X_n) \in \mathbb{Z} [X_1,\dots, X_n]$ verifying an hypothesis of analytic regularity which is mostly satisfied (see Definition \ref{anouvcond2}) when $n \ge 2$  and $h(\textbf{0}) = 1.$ Thus, we characterize precisely the natural boundary of $Z(\mathbf{s})$ when $h$ is not cyclotomic\footnote{If $h$ is cyclotomic one can check that $Z(\mathbf{s}) = h(p^{-s_1},\dots, p^{-s_n})$ continues meromorphically to $\mathbb{C}^n$ (see Remark \ref{cyclo1}).} (see Definition \ref{cyclo}), that is,  the boundary of a maximal domain  on which it  can be meromorphically continued.

%
%


\subsection{Notations.}\label{notation}

For two   positive integers $r$ and $n$ we define:
$$
h\left(X_1,\dots,X_n\right) :=  1+a_{1}X_{1}^{\alpha_{11}}X_{2}^{\alpha_{21}}\cdots X_{n}^{\alpha_{n1}}+..
.+a_{r}X_{1}^{\alpha_{1r}}X_{2}^{\alpha_{2r}}\cdots X_{n}^{\alpha_{nr}};
$$

$$
Z\left(\mathbf{s}\right) := \prod_{p}h\left(p^{-s_{1}},p^{-s_{2}},\dots,p^{-s_{n}}\right) \ \textrm{for} \ \mathbf{s} = \left(s_1,\dots,s_n\right) \in  \mathbb{C}^n,
$$
where  $a_{j}$ for  $ j=1,\dots,r$ are integers and $\alpha_{\ell j}$ for $ j=1,2,\dots,r$  and  $ \ell=1,\dots,n$ are non negative integers.

We also  fix the following notations throughout the article.

We put
$
\boldsymbol{\alpha}:=\left(\alpha_{\ell j}\right)_{(\ell, j)\in \{1,\dots,n\}\times \{1,\dots,r\}}\in \mathbb{M}_{n,r}(\mathbb{N})
$
the matrix encoding the exponents of $h$ whose rows for $\ell \in \{1,\dots,n\}$ are written $\boldsymbol{\alpha}_{\ell \cdot}:=\left(\alpha_{\ell 1},\dots, \alpha_{\ell r}\right)$  and columns for $j\in \{1,\dots,r\}$ are written $\boldsymbol{\alpha}_{\cdot j}:={}^t \!\left(\alpha_{1j},\dots,\alpha_{nj}\right)$.

For $j \in \left\{1,\dots,r\right\}$ we set $
\ds\textbf{X}^{\boldsymbol{\alpha}_{\cdot j}} := X_{1}^{\alpha_{1j}}X_{2}^{\alpha_{2j}}\cdots X_{n}^{\alpha_{nj}}
$ so that $h(\textbf{X})=1+\sum_{j=1}^{r}a_j \textbf{X}^{\boldsymbol{\alpha}_{\cdot j}}$.

We will write for all $\textbf{y} \in \mathbb{R}^{r}$ $
\ds\Vert \textbf{y} \Vert := \sum_{j=1}^{r}y_{j}$ and $\langle \mathbf y \rangle = \mathbb{R}\mathbf{y}.
$

 $\textrm{For} \  \mathbf{s}=\left(s_1,\cdots,s_n\right)\in \mathbb{C}^{n},$ and $\ell \in \left\{1,\cdots,n\right\}$, we put:
  $$
  \sigma_{\ell} := \Re\left(s_{\ell}\right); \ \tau_{\ell} := \Im\left(s_{\ell}\right); \ \boldsymbol{\sigma} := \Re\left( \mathbf{s}\right) := \left(\sigma_1,\cdots,\sigma_n\right); \ \boldsymbol{\tau} := \Im\left( \mathbf{s}\right) := \left(\tau_1,\dots, \tau_n\right).
 $$


We write for $\boldsymbol{\nu}=(\nu_1,\dots,\nu_m)$ and $\boldsymbol{w}={}^t \!(w_1,\dots,w_m)$ the classical matrix product between $\boldsymbol{\nu}$ and $\boldsymbol{w}$:
\begin{center}
$\ds \boldsymbol{\nu}\cdot\boldsymbol{w}  := \summ_{i=1}^{m}\nu_{i} w_{i}.$
\end{center}

For $\boldsymbol{\beta}=(\beta_1,\dots,\beta_r)$ and $\mathbf{s}=(s_1,\dots,s_n)$, we will use throughout the paper the following equality
$$
\sum_{j=1}^{r}\beta_j \left( \mathbf{s} \cdot \bs{\alpha}_{j \cdot}\right)=\sum_{\ell=1}^{n}s_{\ell} \left( \bs{\alpha}_{\ell \cdot} \cdot {}^t \!\bs{\beta}\right):= \mathbf{s}\cdot\boldsymbol{\alpha} \cdot {}^t \! \boldsymbol{\beta};
$$ which  results from the classical identity $\left(\mathbf{s} \cdot \boldsymbol{\alpha}\right) \cdot {}^t \! \boldsymbol{\beta}=\mathbf{s} \cdot \left(\boldsymbol{\alpha} \cdot {}^t \! \boldsymbol{\beta}\right)$. 


\subsection{Statement of main results.}

 We first recall the classic result  of Estermann 
\cite{est} for one variable polynomials. 
\vskip .1 in

\noindent {\bf Theorem. ({\rm Estermann})} {\it Let $h\left(X\right) = 1 + \sum_{m = 1}^r b_{m}X^{m} = \prod_{m = 1}^{r}{\left(1-\alpha_{m}X\right)} \in
\mathbb{Z}[X].$  Let $f\left(s\right)= \prod_{p}h\left(p^{-s}\right),$
which converges for $\Re\left(s\right)>1$. Then:}

\begin{enumerate}
\item[(i)] {\it $f\left(s\right)$ can be meromorphically  extended to $\Re\left(s\right)>0$.}
\item[(ii)] {\it If $|\alpha_{m} | = 1$ $\forall m = 1, .., r$, then $f\left(s\right)$ can be extended to $\mathbb{C}$. Otherwise, $\Re\left(s\right) = 0$ is a natural boundary
for $f$ (i.e. for each  point $s = it$ on this vertical line, 
$f$ cannot be extended as a meromorphic function on any neighborhood of $s=it$).}
\end{enumerate}



%

\begin{rqs} G.  Dahlquist \cite{dal}  generalized this result  to analytic functions  $h$   with isolated singularities within the unit circle.  Later, deep work by Kurokawa \cite{kurokawa2}, \cite{kurokawa3} and Moroz \cite{moroz} extended Estermann's result by allowing
polynomials $h(X)$ whose coefficients were integral linear combinations of complex numbers associated to characters of finite dimensional representations of a topological group. 
\end{rqs}
 Estermann's main result leads naturally to the following basic definition.
\begin{defi}\label{cyclo}
 A rational fraction $Q\left(X_1,\dots,X_n\right)\in \mathbb{Z}(X_1,\dots,X_n)$ is said to be cyclotomic if there exists a finite subset  $I$ of $\Nn$ such that 
\begin{equation}\label{ratcyclo}
\ds Q\left(X_1,\dots,X_n\right) = \prod_{\boldsymbol{\lambda}=\left(\lambda_1,\dots,\lambda_n\right) \in I}\left(1 - X_{1}^{\lambda_1}\cdots X_{n}^{\lambda_n}\right)^{\gamma\left(\bs{\lambda}\right)},
\end{equation}
where  $\gamma\left(\bs{\lambda}\right) \in \mathbb{Z}$ for each $\bs{\lambda} \in I.$ 

In particular when $Q\left(X_1,\dots,X_n\right)\in \mathbb{Z}[X_1,\dots,X_n]$ is a polynomial satisfying (\ref{ratcyclo}), we will say that $Q$ is a cyclotomic polynomial.
\end{defi}

\begin{rqs}\label{cyclo1}
If $Q\left(X_1,\dots,X_n\right)=\prod_{\boldsymbol{\lambda}=\left(\lambda_1,\dots,\lambda_n\right) \in I}\left(1 - X_{1}^{\lambda_1}\cdots X_{n}^{\lambda_n}\right)^{\gamma\left(\bs{\lambda}\right)}$ is a cyclotomic rational fraction, then we have for  $\sigma_{\ell}>1$  $(\ell \in \left\{1,\dots,n\right\})$:

$$ \ds\prod_{p}Q\left(p^{-s_1},\dots, p^{-s_n}\right) \\ 
  =  \ds\prod_{\bs{\lambda} \in I}\zeta\left(\mathbf{s} \cdot {}^t \!\bs{\lambda} \right)^{-\gamma\left(\bs{\lambda}\right)},$$

where $\zeta\left(z\right)$  denotes the Riemann zeta function.  As a result, it is clear that this Euler product
 meromorphically extends to $\mathbb{C}^n$ as a finite product of classical Riemann zeta functions.

  Moreover, if two polynomials $h_1\left(X_1,\dots,X_n\right)$ and $h_2\left(X_1,\dots,X_n\right)$ are such that:
$$
h_1\left(X_1,\dots,X_n\right) = h_2\left(X_1,\dots,X_n\right)  Q(X_1,\dots, X_n)
$$
with $Q$ a cyclotomic rational fraction, 
then the maximal domains of meromorphic continuation of the  Euler products $\prod_{p}h_1\left(p^{-s_{1}},p^{-s_{2}},\dots,p^{-s_{n}}\right)$ and $\prod_{p}h_2\left(p^{-s_{1}},p^{-s_{2}},\dots,p^{-s_{n}}\right)$ coincide.
\end{rqs}

 So from now on, it suffices to  assume the following:
\begin{center}
\vskip .1 in

 {\it \textbf{ The polynomial $h$ is not cyclotomic and has no  cyclotomic factors.}} 
\end{center}
\begin{defi}

Suppose that the polynomial $h$ satisfies the preceding property.

For all $\delta \geq 0$ we put
$
\ds \mathbf{W}(\delta) = \left\{\mathbf{s} \in \mathbb{C}^{n}:   \bs{\sigma} \cdot \bs{\alpha}_{\cdot j} >\delta, \forall  j \in \left\{1,\cdots,r\right\}  \right\}.
$

\end{defi}

%

 For a  polynomial $h$ in $n \geq 1$ variable(s), we first observe that  $Z\left(\mathbf{s}\right)$ defines a holomorphic function of $\mathbf{s}$ in the domain  $\bs{\sigma} \cdot \bs{\alpha}_{\cdot j}>1, (j=1,\dots,r)$. In \cite{bel},  G. Bhowmik, D. Essouabri and B. Lichtin showed that there is a meromorphic continuation of $Z\left(\mathbf{s}\right)$ to
   $\mathbf{W}(0)$. They did so by proving the following result.

\vskip .1 in
\noindent {\bf Theorem ({\rm Bhowmik-Essouabri-Lichtin})} {\it For each
  $\delta>0$, there exists a bounded Euler product $G_{\delta}\left(\mathbf{s}\right)$, absolutely
  convergent on $\mathbf{W}(\delta)$ such that:
\begin{equation} \label{etoile}
  Z\left(\mathbf{s}\right) = \prod_{\stackrel{\bs{\beta}= \left(\beta_{1},\cdots,\beta_{r}\right) \in \mathbb{N}^{r}}{1 \leq   \Vert\bs{\beta}\Vert \leq [\delta^{-1}]}}
{\zeta \left(\mathbf{s}\cdot\boldsymbol{\alpha} \cdot {}^t \! \boldsymbol{\beta}\right)^{\gamma(\boldsymbol{\beta})}G_{\delta}\left(\mathbf{s}\right)}; \ \textrm{where} \ \left\{\gamma(\boldsymbol{\beta}): \boldsymbol{\beta}\in \mathbb{N}^{r}\right\} \subset \mathbb{Z}. 
 \end{equation}}


 In fact, their result is somewhat stronger. They also showed that the function  $Z\left(\mathbf{s}\right)$ does not admit a meromorphic continuation to $\mathbf{W}(\delta)$ for any $\delta<0.$ This followed from the fact that $\textbf{0}$ is  an accumulation point of zeros or  poles of the one variable function $t \longmapsto Z\left(t\cdot\bs{\theta}\right)$ for almost all direction $\bs{\theta}\in \mathbb{R}^n$.


\vspace{0.2cm}

Before announcing the main result, we will first introduce a definition.

Since $\mathbf{W}(0) = \{\mathbf{s}\in \mathbb{C}^n : \Re(\mathbf{s} \cdot \bs{\alpha}_{\cdot j})\geq 0, \forall j=1,\dots,r\}$, then $\partial \mathbf{W}(0)$ is a polyhedron whose faces are of the form $
\ds\mathcal{F}(\boldsymbol{\alpha}_{\cdot e}) = \{\mathbf{s}\in \overline{\mathbf{W}(0)}:\Re(\mathbf{s} \cdot \bs{\alpha}_{\cdot e})=0\}
$
for a vector $\boldsymbol{\alpha}_{\cdot e}\in \{\bs{\alpha}_{\cdot 1},\dots,\bs{\alpha}_{\cdot r}\}$.

We will say by abuse of language that $\mathcal{F}(\boldsymbol{\alpha}_{\cdot e})$ is a face of polar vector $\boldsymbol{\alpha}_{\cdot e}$\footnote{In reality $\boldsymbol{\alpha}_{\cdot e}$ is a polar vector of $\mathcal{F}(\boldsymbol{\alpha}_{\cdot e})\cap \mathbb{R}^n=\{\textbf{x}=(x_1,\dots,x_n)\in\overline{\mathbf{W}(0)}\cap \mathbb{R}^n :  \textbf{x}\cdot\boldsymbol{\alpha}_{\cdot e}=0 \}$.}.

 Now let $\mathcal{F}(\boldsymbol{\alpha}_{\cdot e})$ be a face of $\partial \mathbf{W}(0)$ as above and consider in particular $\boldsymbol{\widehat{\alpha}}_{\cdot e}\in \mathbb{N}^n, \boldsymbol{\widehat{\alpha}}_{\cdot e} \in \mathbb{Q}\boldsymbol{\alpha}_{\cdot e}$ the vector collinear with  $\boldsymbol{\alpha}_{\cdot e}$ whose nonzero components are relatively prime.

We also put $
\Lambda_e:=\{j\in \{1,\dots,r\}: \boldsymbol{\alpha}_{\cdot j}\in \mathbb{Q}\boldsymbol{\alpha}_{\cdot e}\}.
$

It is clear that for all $j\in \Lambda_e$ there exists $q_j\in \mathbb{N}^{*}$ such that $\boldsymbol{\alpha}_{\cdot j}=q_j \boldsymbol{\widehat{\alpha}}_{\cdot e}$.

Then we define:
\begin{displaymath}
\begin{array}{lll}
\ds [h]_e(\textbf{X}) & := & \ds 1+\sum_{\boldsymbol{\alpha}_{\cdot j}\in \mathbb{Q}\boldsymbol{\alpha}_{\cdot e}}a_j \textbf{X}^{\boldsymbol{\alpha}_{\cdot j}} \\ 
\ds \widetilde{[h]_e}(T) & := & \ds 1+\sum_{\boldsymbol{\alpha}_{\cdot j}\in \mathbb{Q}\boldsymbol{\alpha}_{\cdot e}}a_j T^{q_j}\in \mathbb{Z}[T] \ \textrm{verifying} \ \widetilde{[h]_e}(\textbf{X}^{\boldsymbol{\widehat{\alpha}}_{\cdot e}}) = [h]_e(\textbf{X}).
\end{array}
\end{displaymath}

\begin{defi}\label{anouvcond2}

We will say that the face $\mathcal{F}(\boldsymbol{\alpha}_{\cdot e})$ is a non-degenerate face if the one variable polynomial $\widetilde{[h]_e}(T)$ has no multiple root.
 
\end{defi}

  Our main result is as follows.

\vspace{0.5cm}
\noindent {\bf Main Theorem ({\rm see Theorem \ref {aresultatprincipal} \S 3}).}

{\it Assume that the polynomial $h$ is not cyclotomic and admits at least one face $\mathcal{F}(\boldsymbol{\alpha}_{\cdot e})$ of $\partial \mathbf{W}(0)$ non-degenerate in the sense of Definition \ref{anouvcond2}.
 Let  $\mathcal{B}$ denote any open ball centered at any point $\mathbf{s}^0 \in \mathcal{F}(\boldsymbol{\alpha}_{\cdot e})$. Then the function  $Z(\mathbf{s})$ cannot be extended as a meromorphic function to any domain that contains the open ball $\mathcal{B}.$}
\vspace{0.5cm}

 The proof of this result  extends arguments used in  \cite{est}, \cite{dal}, \cite{bel} and \cite{sautoy3}  and adds  two new ideas. The first (see \S 2.2) is to write $h\left(X_1,\cdots,X_n\right)$ as an infinite product of cyclotomic polynomials. This allows us to manipulate the  Euler product  $\prod_{p}h\left(p^{-s_{1}},\cdots,p^{-s_{n}}\right)$  beyond $\mathbf{W}(1)$ with reasonable facility. In particular,  this allows us to give a different proof of the fact that  $Z\left(\mathbf{s}\right)$ meromorphically extends up to $\mathbf{W}(0)$. 

 The second new idea allows us to analyze with good precision how the zeroes of $Z\left(\mathbf{s}\right)$ can accumulate (i.e. when they cannot be cancelled out) inside any  open ball of  any  point of $\partial \mathbf{W}(0).$ In particular, we are able to do this provided that these zeroes  also belong to  a suitable  line (determined by a real direction vector in $\mathbb{R}^n$). Such zeroes must be  zeroes  of appropriate  factors of $Z(\textbf {s})$,  given the identity that is  derived in \S 2.2. One might therefore think that rather precise information about the location of the zeroes of $\zeta(s)$ would be needed to carry out such an analysis. However, this is not really the case. Indeed, the strategy of this work is to concentrate in a first time our attention to ``good'' points of $\partial \mathbf{W}(0)$ in the neighborhood of which we can find an accumulation of zeroes without assuming any assumption about the zeroes of $\zeta(s)$ in order to carry out this
analysis. But a recurrent difficulty in this article precisely consists in proving that these ``good'' points are generic points of $\partial \mathbf{W}(0)$ (see Definition \ref{defigeneric} and Remark \ref{defigeneric2}).


%
%
%
%
%
%
%

 \vspace{0.5cm}


 In \S 4, our goal is to refine the obstruction to continuing $Z(\mathbf{s})$ across $\partial \mathbf{W}(0)$ by restricting the types of sets along which such a continuation is impossible. 

Since $\partial \mathbf{W}(0)$ is a real hypersurface in $\mathbb{C}^n,$ it is natural to ask whether there exists an extension of $Z(\mathbf{s})$ to 
 a real hypersurface that contains $\partial \mathbf{W}(0).$ But for this to be reasonable, we must first
 clarify what type of function this extension should be. One such possible class consists of the 
 C-R functions.
\begin{defi}\label{CR}
 A function $f$ continuous on a real hypersurface $\frak{H}$ of class $\mathcal{C}^1$ in $\mathbb{C}^n$ is said to be C-R (Cauchy-Riemann) if for  any  differential form $\omega$ of bidegree $(n,n-2)$   that is  $\mathcal{C}^{\infty}$ in a neighborhood of $\frak{H}$  and satisfies $\textrm{supp} \  \omega \cap \frak{H}$ is compact, we have $
\ds\int_{\frak{H}}f\overline{\partial}\omega=0.
$ 
\end{defi}

\begin{ex}
 If $F$ is an holomorphic function in a neighborhood of a real hypersurface $\frak{H}$ of class $\mathcal{C}^1$ in $\mathbb{C}^n$, then $f = F\mid_{\frak{H}}$ is a real analytic  C-R function on $\frak{H}$.

\end{ex}

 By restricting our attention to real-analytic hypersurfaces, we are able to show the following.

\noindent {\bf Theorem 3.}

{\it Assume that the polynomial $h$ is not cyclotomic and that $\mathcal{F}(\boldsymbol{\alpha}_{\cdot e})$ is a non-degenerate face.  Let $\frak{H}$ denote any  real-analytic hypersurface that intersects $\mathcal{F}(\boldsymbol{\alpha}_{\cdot e})\subseteq\partial \mathbf{W}(0)$ and is not a subset of $\mathbf{W}(0).$ 

Then  no real-analytic  C-R extension of $Z(\mathbf{s})$ to  $\frak{H}$ can exist. }

\vspace{0.3cm}

In \cite{bel} p. 13, Bhowmik, Essouabri and Lichtin have defined a multivariable zeta function which is naturally associated to a projective toric variety via a projective embedding in $\mathbb{P}^m$ for some $m$. The analytic property of this function provide, via tauberian theorems, some results concerning the density of rational points on such variety. Some results in this direction have been also provided by de la Bretèche and Swinnerton-Dyer (\cite{delabreteche0}, \cite{delabreteche1}, \cite{delabreteche2}). They have obtained precise asymptotic estimations of counting functions relative to a height by studying the associated one variable height zeta function linked to this class of multivariable zeta functions.

Furthermore, in \cite{bel}, Bhowmik, Essouabri and Lichtin have described the form of these multiple zeta functions associated to a toric variety. They proved that such a zeta function is an Euler product $\prod_{p}h^{*}(p^{-s_1},\dots,p^{-s_m})$ where $h^{*}$ is an analytic function satisfying $h^{*}(\textbf{X})=\left(\prod_{\bs{\nu}\in K}\left(1-\textbf{X}^{\bs{\nu}}\right)^{-c(\bs{\nu})}\right)h(\textbf{X})$ with $h(\textbf{X})\in \mathbb{Z}[X_1,\dots,X_{m}]$, $K$ a finite subset of $\mathbb{N}^m$ and $\{c(\bs{\nu})\}_{\bs{\nu}\in K}$ a finite set of positive integers.

As an example of the main result of this paper, we improve a result of \cite{bel} by giving the exact domain of meromorphy of a family of multiple zeta functions associated to the hypersurfaces $x_1\cdots x_n=x_{n+1}^{n}$.

Indeed, we have in particular (see \cite{bel}, Theorem 7) the following expression for the multiple zeta function $Z_{n}(s_1,\dots,s_{n+1})$ asssociated to $x_1\cdots x_n=x_{n+1}^{n}$:
$$
Z_{n}(s_1,\dots,s_{n+1}) = \frac{\prod_{i=1}^{n}\zeta(n s_i +s_{n+1})}{\zeta(s_1 + \cdots + s_{n+1})}\prod_{p}V_n(p^{-s_1},\dots,p^{-s_{n+1}})
$$
with 
$$
V_n(X_1,\dots,X_{n+1}) =  \sum_{\textbf{r}\in \{0,\dots,n-1\}^n; n\mid \lVert \textbf{r}\lVert}X_1^{r_1} \cdots X^{r_n}X_{n+1}^{\lVert \textbf{r}\lVert/n}.
$$
And since all the faces of $\partial \mathbf{W}(0)$ computed from the polynomial $V_n(X_1,\dots,X_{n+1})$ are non-degenerate faces, we obtain the following result:

\begin{cor}
 {\it The exact maximal domain $\mathcal{D}$ of meromorphy (in the sense of Theorem 3) of the zeta function: 
\begin{center}$
\ds Z_n(s_1,\dots,s_{n+1}) = \sum_{\stackrel{x_i\in \mathbb{Z}; \textrm{gcd}(x_i; i=1,\dots,n)=1}{x_1\cdots x_n = x_{n+1}^{n}}}\frac{1}{x_1^{s_1}\cdots x_{n+1}^{s_{n+1}}}
$\end{center}

is given by $\ds\mathcal{D}=\left\{\mathbf{s}\in \mathbb{C}^{n+1}; \forall \textbf{r}\in \{0,\dots,n-1\}^{n}, n\mid \|\textbf{r}\|; s_1 r_1 + \cdots s_n r_n + s_{n+1} \frac{\|\textbf{r}\|}{n}>0  \right\}.$}
\end{cor}

\section{Rewriting $Z\left(\mathbf{s}\right)$ as a product of zeta functions and meromorphic continuation.}

\subsection{An inversion formula for a multivariate arithmetical function.}

The following result, which generalizes the inversion formula for a single variable arithmetical function, will be used  to prove a basic identity in \S 2.2.

\begin{defi}\label{aformuleinversion}
 Given a multivariate arithmetical function $
g: \Nn \longrightarrow \mathbb{C}
$
and a classical one variable arithmetical function
$
f: \mathbf {N}^* (= \mathbb{N} - \{0\})\longrightarrow \mathbb{C},
$
we  define $f \tilde{*} g: \Nn \longrightarrow \mathbb{C}$ as a multivariate arithmetical function  by setting:
$$
\forall \boldsymbol{\beta}\in \Nn, \ f \tilde{*} g\left(\bs{\beta}\right) = \sum_{\stackrel{\stackrel{ \textbf{b} \in \Nn }{m \in \mathbb{N}^*,}}{m \textbf{b} =\boldsymbol{\beta}}}f\left(m\right)g\left(\textbf{b}\right).
$$
\end{defi}

\begin{lemme}
 Given a multivariate arithmetical function $
g: \Nn \longrightarrow \mathbb{C},
$
and two arithmetical functions
$
f_1, f_2: \mathbb{N}^* \longrightarrow \mathbb{C},
$
we have the equality $
f_1\tilde{*}\left(f_2\tilde{*}g\right) = \left(f_1*f_2\right)\tilde{*}g,
$
where $*$  denotes the  standard convolution product between one variable arithmetical functions. 
\end{lemme}

\begin{dem}
 We have for all $ \boldsymbol{\beta}\in \Nn:$

\begin{center}
\begin{tabular}{lll}
$ f_1\tilde{*}\left(f_2\tilde{*}g\right)\left(\bs{\beta}\right) $ & $= \sum_{\stackrel{\stackrel{ \textbf{b} \in \Nn,}{m \in \mathbf{N^*},}}{m \textbf{b} =\boldsymbol{\beta}}}f_1\left(m\right)\left(f_2\tilde{*}g\right)\left(\textbf{b}\right)$ & $= \sum_{\stackrel{\stackrel{ \textbf{b} \in \Nn,}{m \in \mathbf{N^*},}}{m \textbf{b} =\boldsymbol{\beta}}}f_1\left(m\right) \sum_{\stackrel{\stackrel{ \textbf{e} \in \Nn,}{d \in \mathbf{N^*},}}{d \textbf{e} = \textbf{b}}}f_2\left(d\right)g\left(\textbf{e}\right)$ \\ 
\empty & $= \sum_{\stackrel{\stackrel{ \textbf{e} \in \Nn,}{\left(m, d\right) \in (\mathbf{N^*})^2,}}{md \textbf{e} =\boldsymbol{\beta}}}f_1\left(m\right)f_2\left(d\right)g\left(\textbf{e}\right)$  
 & $= \sum_{\stackrel{\stackrel{ \textbf{e} \in \Nn,}{k \in \mathbb{N}^{*},}}{k \textbf{e} =\boldsymbol{\beta}}}\left(\sum_{md=k}f_1\left(m\right)f_2\left(d\right)\right)g\left(\textbf{e}\right)$ \\ 
\empty & $= \sum_{\stackrel{\stackrel{ \textbf{e} \in \Nn,}{k \in \mathbb{N}^{*},}}{k \textbf{e} =\boldsymbol{\beta}}}\left(f_1*f_2\right)\left(k\right)g\left(\textbf{e}\right)$ & $= \left(f_1*f_2\right)\tilde{*}g\left(\bs{\beta}\right)$. \ This completes the proof. \CQFD
\end{tabular}
\end{center}

\end{dem}

Thus, if $f$ is  an invertible  single variable arithmetic function (with respect to convolution),  and if we know $f\tilde{*}g$, we are able to find $g$:

\begin{cor}
{\it  Let $f: \mathbb{N}^* \longrightarrow \mathbb{C}$  be an invertible arithmetic function with inverse $f^{-1}$.  Then, for any multivariate arithmetical function $g: \Nn \longrightarrow \mathbb{C}$ we have for all $ \boldsymbol{\beta}\in \Nn$, $
g\left(\bs{\beta}\right) = f^{-1}\tilde{*}\left(f\tilde{*}g\right)\left(\bs{\beta}\right). 
$}
\end{cor}

\subsection{Meromorphic continuation of $Z\left(\mathbf{s}\right)$.}\label{meromorphiccontinuation}

 In this subsection we give a different proof (from that in \cite{bel}) that $Z(\mathbf{s})$ has a meromorphic continuation in $\mathbf{W}(0).$ Our argument is based upon an expression  for any polynomial $h$ as in \S 1.1 as an
infinite product of cyclotomic polynomials.

 Consider the following quantity:
\begin{equation}\label{adist-h-zero}
C:=C(h) = \frac{1}{|a_1|+\cdots + |a_r|}.
\end{equation}

It is clear that if each $|Y_i|<C$ ($i=1,\dots,r$) then:
\begin{equation}\label{aCoptimal}
\left|\sum_{i=1}^{r}a_i Y_i\right|<1;
\end{equation}
and we verify that $C=C(h)$ is maximal among the  $C$ satisfying (\ref{aCoptimal}).

\begin{lemme}\label{aetoile2}
If each $|Y_i| < C=C(h)= \frac{1}{|a_1|+\cdots + |a_r|},$ then we have:
\begin{equation}
1+a_{1}Y_{1}+\cdots +a_{r}Y_{r} = \prod_{\boldsymbol{\beta}\in \Nr}\left(1-Y_{1}^{\beta_{1}}\cdots Y_{r}^{\beta_{r}}\right)^{\gamma(\boldsymbol{\beta})},
\end{equation}
where the right side converges absolutely and each $\gamma(\bs{\beta})\in \mathbb{Z}$ and satisfies the equality:

 \begin{equation}\label{expgama}\gamma(\boldsymbol{\beta}) = \displaystyle \sum_{\substack{  
\textbf{b} \in \mathbb{N}^r\setminus \left\{\textbf{0}\right\} \\ m \in \mathbb{N} \\ m \textbf{b} =\boldsymbol{\beta} }} \left(\left(-1\right)^{\Vert \textbf{b} \Vert } \frac{\mu\left(m\right)}{m} \frac{\left(\Vert \textbf{b} \Vert -1\right)!}{b_{1}!\cdots b_{r}!} a_{1}^{b_{1}}\cdots a_{r}^{b_{r}} \right)\in \mathbb{Z} \ \ (\mu\left(.\right) \textrm{denotes the  M\"{o}bius function).}\end{equation}

In addition we have $
\ds\left|\gamma(\boldsymbol{\beta})\right| \ll C^{-\lVert\boldsymbol{\beta}\lVert}
$
uniformly in $\boldsymbol{\beta}\in \Nr$.
\end{lemme}

\begin{rqs}
 The fact that $\gamma(\boldsymbol{\beta}) \in \mathbb{Z}$ can be proved by recurrence on $\Arrowvert\bs{\beta}\Arrowvert$ in the same way as in \cite{est} p. $448$ in the supplement of the appendix added on January 14, 1928.

\end{rqs}

\begin{dem}[\textbf{Lemma \ref{aetoile2}}]

Let us estimate $\gamma(\boldsymbol{\beta})$ by using the expression given in (\ref{expgama}).  We have:

\begin{displaymath}
\begin{array}{lll}
\left|\gamma(\boldsymbol{\beta})\right| & \leq & \ds \sum_{\substack{  
\textbf{b} \in \Nr \\ m \in \mathbb{N} \\ m \textbf{b} =\boldsymbol{\beta} }}\frac{1}{m}\frac{(\lVert \textbf{b}\lVert-1)!}{b_1! \cdots, b_r!}\prod_{j=1}^{r}|a_j|^{b_j} 
  \leq  \ds \frac{1}{\lVert\boldsymbol{\beta}\lVert}\sum_{\substack{  
\textbf{b} \in \Nr \\ m \in \mathbb{N} \\ m \textbf{b} =\boldsymbol{\beta} }}\frac{\lVert \textbf{b}\lVert!}{b_1! \cdots, b_r!}\prod_{j=1}^{r}|a_j|^{b_j} \\ 
 & \leq & \ds \frac{1}{\lVert\boldsymbol{\beta}\lVert}\sum_{m\mid \lVert\boldsymbol{\beta}\lVert}\sum_{\stackrel{\textbf{b}\in \Nr}{\lVert \textbf{b}\lVert} = \frac{\lVert\boldsymbol{\beta}\lVert}{m}}\frac{\lVert \textbf{b}\lVert!}{b_1! \cdots, b_r!}\prod_{j=1}^{r}|a_j|^{b_j}  
  \leq  \ds \frac{1}{\lVert\boldsymbol{\beta}\lVert}\sum_{m\mid \lVert\boldsymbol{\beta}\lVert} \left(|a_1| + \cdots + |a_r|\right)^{\frac{\lVert\boldsymbol{\beta}\lVert}{m}} \\ 
 & \leq & \ds \frac{\tau(\lVert\boldsymbol{\beta}\lVert)}{\lVert\boldsymbol{\beta}\lVert}\left(|a_1| + \cdots + |a_r|\right)^{\lVert\boldsymbol{\beta}\lVert} \ \textrm{where} \ \tau(\lVert\boldsymbol{\beta}\lVert) \ \textrm{is the number of divisors of} \ \lVert\boldsymbol{\beta}\lVert \\ 
 & \leq & \ds \frac{\tau(\lVert\boldsymbol{\beta}\lVert)}{\lVert\boldsymbol{\beta}\lVert}C^{-\lVert\boldsymbol{\beta}\lVert} 
  \ll  C^{-\lVert\boldsymbol{\beta}\lVert} \ \textrm{uniformly in} \ \boldsymbol{\beta}\in \Nr.
\end{array}
\end{displaymath}

Now put
$
\ds G(\textbf{Y}):= \prod_{\boldsymbol{\beta}\in \Nr}\left(1-Y_1^{\beta_1}\cdots Y_r^{\beta_r}\right)^{\gamma(\boldsymbol{\beta})}.
$

Verify that $G(\textbf{Y})$ is an holomorphic function in $\{\textbf{Y}\in \mathbb{C}^r : \max_i |Y_i|<C \}.$

So let $0<C_1<C$ and let us prove the convergence of $G$ for $\max_i |Y_i|<C_1$.

We have for $|Y_i|<C_1$:

\begin{displaymath}
\begin{array}{llllll}
\ds\sum_{\boldsymbol{\beta}\in \Nr}|\gamma(\boldsymbol{\beta})| |\textbf{Y}^{\bs{\beta}}| & \leq &\ds \sum_{\boldsymbol{\beta}\in \Nr}C^{-\lVert\boldsymbol{\beta}\lVert}C_1^{\lVert\boldsymbol{\beta}\lVert}  
  =  \ds\sum_{\boldsymbol{\beta}\in \Nr}\left(\frac{C_1}{C}\right)^{\lVert\boldsymbol{\beta}\lVert} \\ 
 & \leq & \ds \left(\sum_{k=0}^{+\infty}\left(\frac{C_1}{C}\right)^k\right)^{r} = \frac{1}{\left(1-\frac{C_1}{C}\right)^r} < +\infty.
\end{array}
\end{displaymath}

Hence $\textbf{Y}\longmapsto G(\textbf{Y})$ converges absolutely and defines an holomorphic function in $\mathcal{D}:=\{\textbf{Y}\in \mathbb{C}^r : \max_i |Y_i|<C\}$; moreover, for all $\textbf{Y}\in \mathcal{D}$ we have:

\begin{equation}\label{acvabs}
\begin{array}{llllll}
\log(G(\textbf{Y})) & = &\ds \sum_{\boldsymbol{\beta}\in \Nr}\gamma(\boldsymbol{\beta})\log(1-\textbf{Y}^{\beta})  
 & = & \ds-\sum_{\boldsymbol{\beta}\in \Nr}\gamma(\boldsymbol{\beta})\sum_{m=1}^{+\infty}\frac{1}{m}\textbf{Y}^{m\beta}. 
\end{array}
\end{equation}

And since (\ref{acvabs}) converges absolutely for $\textbf{Y}\in \mathcal{D}$, we have:

\begin{equation}\label{alogG}
\begin{array}{lll}
\log(G(\textbf{Y})) & = &\ds -\sum_{\textbf{b}\in \Nr}\left(\sum_{\boldsymbol{\beta}\in \Nr, m\in \mathbb{N}, m\bs{\beta}=\textbf{b}}\frac{\gamma(\boldsymbol{\beta})}{m}\right)\textbf{Y}^{\textbf{b}} \\ 
 & = & \ds -\sum_{\textbf{b}\in \Nr}\frac{1}{\lVert \textbf{b}\lVert}\left(\sum_{\boldsymbol{\beta}\in \Nr, m\in \mathbb{N}, m\bs{\beta}=\textbf{b}}\lVert\boldsymbol{\beta}\lVert\gamma(\boldsymbol{\beta})\right)\textbf{Y}^{\textbf{b}}.
\end{array}
\end{equation}

But we have:
 $$
 -\Vert\bs{\beta}\Vert\gamma(\boldsymbol{\beta}) = \summ_{m \in \mathbb{N}, \textbf{b} \in \Nr, m \textbf{b} =\boldsymbol{\beta}}\mu\left(m\right)\left(-1\right)^{\Vert \textbf{b} \Vert - 1}\frac{\Vert \textbf{b} \Vert !}{b_1 !\cdots b_r !}a_1^{b_1}\cdots a_r^{b_r};
 $$

 And by Definition \ref{aformuleinversion}, we have: 
 
  $$ -\Vert \bs{\beta}\Vert\gamma(\boldsymbol{\beta}) = \mu \tilde{*} g\left(\bs{\beta}\right),$$

 where $\mu$ denotes the  M\"{o}bius function, and the function  $g$ is defined by $g\left(\bs{\beta}\right) = \left(-1\right)^{\Vert \bs{\beta}\Vert - 1}\frac{\Vert\bs{\beta}\Vert !}{\beta_1 !\cdots \beta_r !}a_1^{\beta_1}\cdots a_r^{\beta_r}$ for $\boldsymbol{\beta}\in \Nr$.

Since the function $\textbf{1}$ being $1$ everywhere is the inverse of the function $\mu$ we obtain:
\begin{equation}\label{ainversemobius}
 \summ_{\exists m \in \mathbb{N}, \boldsymbol{\beta}\in \Nr, m\bs{\beta} = \textbf{b}}-\Vert \bs{\beta}\Vert\gamma(\boldsymbol{\beta}) = \left(-1\right)^{\Vert \textbf{b}\Vert - 1}\frac{\Vert \textbf{b} \Vert !}{b_1 !\cdots b_r !}a_1^{b_1}\cdots a_r^{b_r}.
 \end{equation}

So the identity (\ref{alogG}) provides:

\begin{displaymath}
\begin{array}{lll}
\log(G(\textbf{Y})) & = & \ds\sum_{\textbf{b}\in \Nr}\frac{1}{\lVert \textbf{b}\lVert}(-1)^{\lVert\textbf{b}\lVert}\frac{\lVert \textbf{b}\lVert!}{b_1!\cdots b_r!}\textbf{a}^{\textbf{b}}\textbf{Y}^{\textbf{b}}  
  =  \ds\sum_{k=1}^{+\infty}\frac{(-1)^k}{k}\left(\sum_{\lVert \textbf{b}\lVert=k} \frac{k!}{b_1!\cdots b_r!} \prod_{j=1}^{r}(a_j Y_j)^{b_j}\right) \\ 
 & = & \ds\sum_{k=1}^{+\infty} \frac{(-1)^k}{k}\left(a_1 Y_1 + \cdots + a_r Y_r\right)^k 
  =  \ds\log\left(1+a_1Y_1+\cdots + a_r Y_r\right);
\end{array}
\end{displaymath}

which completes the proof of this lemma.
 \CQFD

\end{dem}

 We fix a polynomial $h$ as in \S \ref{notation} and consider the quantity $C=C(h)$ introduced in (\ref{adist-h-zero}):

\begin{cor}\label{aexpansion}
{\it If each  $|\textbf{X}^{\boldsymbol{\alpha}_{\cdot j}}| < C$ for $j \in \{1,\dots,r\}$, then:
\begin{displaymath}
\begin{array}{lll}
1 + a_1 \textbf{X}^{\bs{\alpha}_{\cdot 1}}+ \dots + a_r \textbf{X}^{\bs{\alpha}_{\cdot r}} &=&\displaystyle \prod_{\boldsymbol{\beta}\in \Nr}\left(1-\prod_{\ell=1}^{n}X_{\ell}^{\boldsymbol{\alpha}_{\ell \cdot} \cdot {}^t \! \boldsymbol{\beta}}\right)^{\gamma(\boldsymbol{\beta})} 
 = \displaystyle \prod_{\boldsymbol{\beta}\in \Nr}\left(1-\textbf{X}^{\boldsymbol{\alpha} \cdot {}^t \! \boldsymbol{\beta}}\right)^{\gamma(\boldsymbol{\beta})}.
\end{array}
\end{displaymath}}
\end{cor}


\begin{theo} \label{aana}
The function $Z\left(\mathbf{s}\right)$ is meromorphic on $\mathbf{W}(0)$. 

Moreover, if we write for all $\delta>0$ $M_\delta=\left[C^{-\frac{1}{\delta}}\right]+1$ $(M_\delta \in \mathbb{N})$,  there exists $A_{M_{\delta}}$ meromorphic on $\mathbf{W}(\delta)$ with possible zeros or poles in the set:
$$
\Phi_{\delta} = \left\{\mathbf{s} \in \mathbf{W}(\delta) \mid \exists \boldsymbol{\beta}\in \mathbb{N}^r, \mathbf{s}\cdot\boldsymbol{\alpha} \cdot {}^t \! \boldsymbol{\beta} = \rho, \ \rho \ \textrm{zero \ or \ pole \ of \ } \zeta\left(s\right) \right\};
$$ 

and such that this relation holds on $\mathbf{W}(\delta)$:

$$
Z\left(\mathbf{s}\right) = \prod_{p \leq M_\delta}h\left(p^{-s_1},\dots,p^{-s_n}\right)A_{M_{\delta}}\left(\mathbf{s}\right).
$$
\end{theo}


\begin{dem}
 We show that $Z\left(\mathbf{s}\right)$ is meromorphic on $\mathbf{W}(\delta)$ for all $\delta>0$.

 We know that if each $|\textbf{X}^{\boldsymbol{\alpha}_{\cdot j}}| < C$ then
$
\ds h\left(X_1,\dots,X_n\right) = \prod_{\boldsymbol{\beta}\in \Nr}\left(1-\textbf{X}^{\boldsymbol{\alpha} \cdot {}^t \! \boldsymbol{\beta}}\right)^{\gamma(\boldsymbol{\beta})},
$
where the right side converges absolutely.

We know according to the previous lemma that $
\ds |\gamma(\boldsymbol{\beta})| = O\left(C^{-\Vert\bs{\beta}\Vert} \right).
$
Thus, for $\boldsymbol{\beta}\in \Nr$ and $\mathbf{s} \in \mathbf{W}(\delta)$ we have:

\begin{center}
\begin{tabular}{ll}
$\displaystyle\sum_{p>M_\delta}\Bigl|\gamma(\boldsymbol{\beta})p^{-\mathbf{s}\cdot\bs{\alpha}\cdot {}^t \! \bs{\beta}} \Bigl|$ & $\leq |\gamma(\boldsymbol{\beta})|\displaystyle\sum_{p>M_\delta}p^{-\bs{\sigma}\cdot\bs{\alpha}\cdot {}^t \! \bs{\beta}}$  
  $\leq |\gamma(\boldsymbol{\beta})|\displaystyle\int_{M_\delta}^{+\infty}x^{-\Vert\bs{\beta}\Vert\delta} dx$ \\ 
 & $= O\left(|\gamma(\boldsymbol{\beta})|M_\delta^{-\Vert\bs{\beta}\Vert\delta+1} \right)$  
  $= O\left(C^{-\Vert\bs{\beta}\Vert}M_\delta^{-\Vert\bs{\beta}\Vert\delta+1} \right).$

\end{tabular}
\end{center}

 Since  $|x|<1$ implies $
\ds\sum_{\boldsymbol{\beta}\in \Nr}x^{\Vert\bs{\beta}\Vert} = \left( \frac{1}{1-x} \right)^r - 1 < +\infty
$
 and since  $M_\delta> C^{-\frac{1}{\delta}}$, we see that $\epsilon$ can be supposed  small enough so  that $M_\delta>\left(C-2\epsilon\right)^{-\frac{1}{\delta}}.$ Then we have:
$$
\sum_{p>M_\delta}\sum_{\boldsymbol{\beta}\in \Nr}\left| \frac{\gamma(\boldsymbol{\beta})}{p^{\mathbf{s}\cdot\bs{\alpha}\cdot {}^t \! \bs{\beta}}} \right| < +\infty;
$$

and thus according to Fubini's theorem applied with the counting measure we obtain:
$$
\prod_{p>M_\delta}\prod_{\boldsymbol{\beta}\in \Nr}\left( 1-p^{-\mathbf{s}\cdot\bs{\alpha}\cdot {}^t \! \bs{\beta}} \right)^{\gamma(\boldsymbol{\beta})} = \prod_{\boldsymbol{\beta}\in \Nr}\prod_{p>M_\delta}\left( 1-p^{-\mathbf{s}\cdot\bs{\alpha}\cdot {}^t \! \bs{\beta}} \right)^{\gamma(\boldsymbol{\beta})}.
$$

 We then have, initially  for $\sigma_k > -\frac{\log\left(C\right)}{\log\left(2\right)}, \left(k=1,\dots,n\right)$ (i.e. $|p^{-s_k}| < C$ for each $k$), and subsequently by analytic continuation to  $\mathbf{W}(\delta)$, the following equality:
$$
\prod_{p>M_\delta}h\left(p^{-s_1},\dots,p^{-s_n}\right) = \prod_{\boldsymbol{\beta}\in \Nr} \zeta_{M_{\delta}} \left(\mathbf{s}\cdot\bs{\alpha}\cdot {}^t \! \bs{\beta} \right)^{-\gamma(\boldsymbol{\beta})},
$$
where 
\begin{equation} \label{azetamdelta} \zeta_{M_{\delta}} \left(\mathbf{s}\cdot\bs{\alpha}\cdot {}^t \! \bs{\beta}\right) = \zeta \left(\mathbf{s}\cdot\bs{\alpha}\cdot {}^t \! \bs{\beta}\right)\prod_{p \leq M_\delta}\left( 1-p^{-\mathbf{s}\cdot\bs{\alpha}\cdot {}^t \! \bs{\beta}} \right).
\end{equation}

We  then notice that for all $z$ with $\Re\left(z\right)>0$, $\zeta\left(z\right)$ and $\zeta_{M_\delta}\left(z\right)$  have exactly the same zeros with the same multiplicities  since $\Re\left(\mathbf{s}\cdot\bs{\alpha}\cdot {}^t \! \bs{\beta}\right)=\bs{\sigma}\cdot\bs{\alpha}\cdot  \bs{\beta} =  \sum_j \beta_j \left(\bs{\sigma}\cdot\bs{\alpha}_{\cdot j}\right) > \delta$ when $\mathbf{s} \in \mathbf{W}(\delta)$ and $\left(1-p^{-\Re\left(z\right)}\right)$ does not vanish when  $p \leq M_\delta$ and $\Re (z) > \delta.$

Set $\ds A_{M_\delta}\left(\mathbf{s}\right) = \prod_{\boldsymbol{\beta}\in \Nr} \zeta_{M_\delta} \left(\mathbf{s}\cdot\bs{\alpha}\cdot {}^t \! \bs{\beta} \right)^{-\gamma(\boldsymbol{\beta})}.$

The zeroes or  poles of $A_{M_\delta}$  must belong to  $\Phi_{\delta}$.
 Moreover, $A_{M_\delta}$ is meromorphic on $\mathbf{W}(\delta)$.

 Indeed, we write $
\ds A_{M_\delta}\left(\mathbf{s}\right) = A_{1,M_\delta}\left(\mathbf{s}\right)A_{2,M_\delta}\left(\mathbf{s}\right),
$
with $\ds A_{1,M_\delta}\left(\mathbf{s}\right) = \prod_{\Vert\bs{\beta}\Vert \leq [\delta^{-1}]} \zeta_{M_\delta} \left(\mathbf{s}\cdot\bs{\alpha}\cdot {}^t \! \bs{\beta} \right)^{-\gamma(\boldsymbol{\beta})}$
 and
 $\ds A_{2,M_\delta}\left(\mathbf{s}\right) = \prod_{\Vert\bs{\beta}\Vert>[\delta^{-1}]} \zeta_{M_\delta} \left(\mathbf{s}\cdot\bs{\alpha}\cdot {}^t \! \bs{\beta}\right)^{-\gamma(\boldsymbol{\beta})}.$

$A_{1,M_\delta}$ is clearly meromorphic on $\mathbb{C}^n$  since it equals a  finite product of meromorphic functions.

For $A_{2,M_\delta}$, we have for $\Vert\bs{\beta}\Vert \geq [\delta^{-1}]+1$, $\ds 
\Re\left(\mathbf{s}\cdot\bs{\alpha}\cdot {}^t \! \bs{\beta}\right)=\bs{\sigma}\cdot\bs{\alpha}\cdot {}^t \! \bs{\beta} \geq \Vert\bs{\beta}\Vert\delta>1.
$

Hence
$
\ds\left| \zeta_{M_\delta} \left(\mathbf{s}\cdot\bs{\alpha}\cdot {}^t \! \bs{\beta} \right)-1 \right| \leq \sum_{k=M_\delta+1}^{+\infty}k^{-\Vert\bs{\beta}\Vert\delta}<\int_{M_\delta}^{+\infty}x^{-\Vert\bs{\beta}\Vert\delta}dx = \frac{1}{\Vert\bs{\beta}\Vert\delta-1}\frac{1}{M_\delta^{\Vert\bs{\beta}\Vert\delta-1}}.
$

 Therefore,
\begin{displaymath}
\begin{array}{c}
\displaystyle \sum_{\Vert\bs{\beta}\Vert \geq [\delta^{-1}]+1}\left|\gamma(\boldsymbol{\beta})\right|\left| \zeta_{M_\delta} \left(\mathbf{s}\cdot\bs{\alpha}\cdot {}^t \! \bs{\beta} \right)-1 \right| 
  \displaystyle \leq \sum_{\Vert\bs{\beta}\Vert \geq [\delta^{-1}]+1}\left(\frac{C^{-1}}{\left(M_\delta\right)^{\delta}}\right)^{\Vert \beta\Vert}\frac{M_\delta}{\Vert\bs{\beta}\Vert\delta-1}< +\infty;
\end{array}
\end{displaymath}
which proves the meromorphy of $A_{2,M_\delta}$, and so of $A_{M_\delta}$ on $\mathbf{W}(\delta)$.

Finally,  since
$
\ds Z\left(\mathbf{s}\right) = A_{M_\delta}\left(\mathbf{s}\right)\prod_{p \leq M_\delta}h\left(p^{-s_1},\dots,p^{-s_n}\right);
$
 it follows that $Z\left(\mathbf{s}\right)$ is  meromorphic on $\mathbf{W}(\delta)$. 

This completes the proof. \CQFD
\end{dem}

\section{Natural boundary of $Z(\mathbf{s})$.}\label{naturalboundary}

In the preceding section we have proved that $Z(\mathbf{s})$ meromorphically extends to $\mathbf{W}(0)$.
The aim of this section is to prove that $\partial \mathbf{W}(0)$ is a natural boundary of $Z(\mathbf{s})$ when $h$ is not cyclotomic and admits at least one non-degenerate face $\mathcal{F}(\boldsymbol{\alpha}_{\cdot e})$ in the sense of Definition \ref{anouvcond2}.

Let us start by giving the definition of what we mean by a generic set:
\begin{defi}\label{defigeneric}
 A set $G\subseteq E$ is said to be generic in $E$ if the complement of $G$ in $E$ has empty interior.
\end{defi}

\begin{rqs}\label{defigeneric2}
 In the following when we will use the term ``generic point'' in a set $E$, it will be understood that we consider any point belonging to some generic set in $E$.
\end{rqs}

The  underlying idea behind the proof of the main theorem  is the observation that  if there exists a meromorphic extension  that is defined at some point of $\partial \mathbf{W}(0)$, then such a   function would also have to be defined in  some  open ball of a generic point of $\partial \mathbf{W}(0)$.  
  Consequently it suffices to restrict our attention from the beginning to any generic point $\mathbf{s}^0$ of $\partial \mathbf{W}(0)$, and then  prove the existence of an accumulation of zeroes or poles of the restriction of $Z(\mathbf{s})$ to  a suitable line through $\mathbf{s}^0$ with real direction vector $\bs{\theta}\in \mathbb{Q}^n$ (i.e. with parametrization $t \to \mathbf{s}^0 + t \bs{\theta}$). In fact, the generic set we end up 
using is constructed in  several  basic  steps (see the condition (\ref{agenericite}) and the proofs of Lemmas  \ref{arg-generique}, \ref{perturbation-sigma0}, and \ref{lemmeBe}).

 \textbf{Note:} Since the argument here is particularly difficult to follow because of technical difficulties, some \textit{\textbf{comments}} (clearly presented as such) will be add throughout this section to make the proof more readable. Obviously these comments are not part of the proof itself. 
\vspace{0.3cm}

We recall the assumption (see \S 1) that \textit{$h$ is not cyclotomic and does not contain any cyclotomic factor}. For $Z(\mathbf{s})$ determined by $h$, our main result is as follows.


\begin{theo}\label{aresultatprincipal}
 If $\mathcal{F}(\boldsymbol{\alpha}_{\cdot e})$ is a non-degenerate face of  $\partial \mathbf{W}(0)$ in the sense of Definition \ref{anouvcond2} then  $Z\left(\mathbf{s}\right)$ cannot be meromorphically continued into any open ball $\mathcal{B}$ centered at any  point  $\mathbf{s}^0 \in \mathcal{F}(\boldsymbol{\alpha}_{\cdot e})\subseteq\partial \mathbf{W}(0)$. 
\end{theo}

The proof of Theorem \ref{aresultatprincipal} will be done in two  basic parts.

In  Part 1, we will show  the existence of an  accumulation of zeroes coming from the factors  $t\longmapsto\prod_{p \leq M_\delta}h\left(p^{-s^0_1-t\theta_1},\cdots,p^{-s^0_n-t\theta_n}\right)$. 
 Here, the delicate point to verify will be   the existence of such zeroes with {\it positive} real part, meaning that all such zeroes lie  ``to the right of" (or ``above")  $\partial \mathbf{W}(0)$.

 In Part 2 we prove that  these zeroes are {\it not}  cancelled by possible poles coming from $\zeta$-factors of $A_{M_\delta}\left(\mathbf{s}^0+t \bs{\theta}\right)$. For this we first remark that the zeroes or poles of $A_{M_\delta}\left(\mathbf{s}^0+t \bs{\theta}\right)$ can be expressed as follows:
$$
t(\bs{\beta},\rho) = \frac{\rho - \mathbf{s}^0\cdot\bs{\alpha}\cdot {}^t \! \bs{\beta}}{\bs{\theta}\cdot\bs{\alpha}\cdot {}^t \! \bs{\beta}}, 
$$ 
where $\boldsymbol{\beta}\in \Nr$ and $\rho$ designates the pole $1$ or a non-trivial zero of the Riemann zeta function.
The key observation  in Part 2 is that we can always choose a generic $\mathbf{s}^0 \in \partial \mathbf{W}(0)$  (via a second generic condition we impose upon both  $\boldsymbol{\sigma}^0$ and $\boldsymbol{\tau}^0$,   see Lemmas \ref {perturbation-sigma0} and \ref {lemmeBe})   so that  for all $\rho$ and for all $\boldsymbol{\beta}\in \Nr$ $
\ds h\left(p^{-s^0_1 - t(\bs{\beta},\rho)\theta_1},\dots,p^{-s^0_n - t(\bs{\beta},\rho)\theta_n}\right) \neq 0.
$

\vspace{0.2cm}

 So let $\mathbf{s}^0\in \mathcal{F}(\boldsymbol{\alpha}_{\cdot e})$ verifying $\bs{\sigma}^0 \cdot \bs{\alpha}_{\cdot e}=0$ and consider an open ball $\mathcal{B}(\mathbf{s}^0) = \mathcal{B}$ of radius arbitrarly small around the point  $\mathbf{s}^0$.

Moving $\mathbf{s}^0 \in \partial \mathbf{W}(0)\cap \mathcal{B}$ if necessary, we can assume that $\boldsymbol{\sigma}^0$ verify the following genericity conditions given a column $\textbf{w}\in \mathbb{Q}^n$: there exists $e\in \{1,\dots,r\}$ such that
\begin{equation}\label{agenericite}
\boldsymbol{\sigma}^0\cdot\textbf{w} = 0 \Longleftrightarrow  \textbf{w} \in \mathbb{Q}\boldsymbol{\alpha}_{\cdot e}. \end{equation}

 In particular we have: 
if $\bs{\beta},\bs{\beta}'\in \Nr$ are such that $\bs{\sigma}^0\cdot\bs{\alpha}\cdot {}^t \! \bs{\beta}= \bs{\sigma}^0\cdot\bs{\alpha}\cdot {}^t \! \bs{\beta}'$, then $\bs{\alpha}\cdot {}^t \! \bs{\beta}' \in \bs{\alpha}\cdot {}^t \! \bs{\beta} + \mathbb{Q}\boldsymbol{\alpha}_{\cdot e}.$ 

Moreover, by rearranging the indexes if necessary, we will suppose from now on without loss of generality that
\begin{equation}\label{alphaenonnul}
 \alpha_{n e}\neq 0.
\end{equation}

\begin{defi}\label{ae-iemepartie}
 We define the $e$-th main part of $h$ to be the polynomial $
\ds [h]_{e}(\textbf{X}) =1+ \sum_{\boldsymbol{\alpha}_{\cdot j} \in \langle \boldsymbol{\alpha}_{\cdot e} \rangle} a_j \textbf{X}^{\boldsymbol{\alpha}_{\cdot j}}.
$
\end{defi}

\begin{defi}
We set  
$\displaystyle  \Lambda_{e}  = \big\{j \in \left\{1,\dots,r\right\} : \boldsymbol{\alpha}_{\cdot j} \in  \langle \boldsymbol{\alpha}_{\cdot e} \rangle \big\}$ and 
$\ds B_e = \big\{\boldsymbol{\beta}\in \textbf{N}^r : \beta_j = 0 \quad  \text{if} \quad j \notin \Lambda_e \big\}.$
\end{defi}

\noindent {\bf Proof of Part 1.}

\vskip .1 in

\begin{comment}
 The aim of this part is to find an expression of zeroes of the factors involving $h$ in $t\longmapsto Z(\mathbf{s}^0+t\bs{\theta})$. The fact we fix a point $\mathbf{s}^0\in \partial \mathbf{W}(0)$ in a direction $\bs{\theta}\in \mathbb{R}^n$ naturally reduces the study of the zeros of the multivariate polynomial $h$ to those of a generalized polynomial of two variables $W(X,Y)$ (depending on the parameters $\bs{\sigma}^0, \bs{\tau}^0, \bs{\theta}$ and $p$). The main difficulty of this part will be to prove the existence of Puiseux branches of such a generalized polynomial so to obtain an expression of the zeroes of $t\longmapsto h\left(p^{-s^0_1-t\theta_1},\dots, p^{-s^0_n-t\theta_n}\right)$ for all prime number $p$ large enough. This is what we do in Proposition \ref{athpuiseuxpolgen}.
\end{comment}

 Recall that thanks to \S \ref{meromorphiccontinuation},  we have, for any $\delta>0$, the following expression  for  $Z(\mathbf{s})$  whenever  $\mathbf{s} \in \mathbf{W}(\delta)$: $
\ds Z\left(\mathbf{s}\right) = \prod_{p \leq M_\delta}h\left(p^{-s_1},\cdots,p^{-s_n}\right)A_{M_\delta}\left(\mathbf{s}\right);
$
where $M_{\delta} =\left[C^{-\frac{1}{\delta}}\right]+1$ and $C = C(h) = \frac{1}{|a_1| + \cdots + |a_r|}
$.
 Then consider $Z|_{L}$, where $L $ denotes that part of the line $L(\mathbf{s}^0, \bs{\theta})$ parametrized by $\mathbf{s}^0 + t \bs{\theta}$ with $\Re (t) > 0.$ The aim is to prove the existence of an accumulation of zeroes of $Z|_L$ in any rectangle $\Xi_{u,\eta}$ depending on two parameters $(u,\eta>0)$:
\begin{center}
\begin{tabular}{ll}
$\Xi_{u,\eta}:$ & $0<\Re\left(t\right)< 1 $\\ 
 & $0<u<\Im\left(t\right)<u+\eta.$
\end{tabular}
\end{center}
To be able to use  the infinite product expression given by Theorem \ref{aana} for $Z|_L$, it suffices to impose the following  condition on  $\bs{\theta}\in\mathbb{Q}^n$: 
\begin{equation} \label{ahyp_theta1}   
          \bs{\theta} \cdot \bs{\alpha}_{\cdot j} \geq 1 \ \textrm{for all} \ j \in \{1,\dots,r\}.  
         \end{equation} 
Then it is  simple to verify that when $\bs{\theta}$ satisfies   $(\ref{ahyp_theta1})$ and  $\Re(t) \geq \delta >0$ then $\mathbf{s}^0 + t\bs{\theta}\in \mathbf{W}(\delta).$ Indeed, we observe that for all $j \in \{1,\dots,r\},$ \ $\left(\bs{\sigma}^0 + \Re (t) \bs{\theta}\right)\cdot\boldsymbol{\alpha}_{\cdot j}  \geq \Re (t)  \geq \delta$.
   As a result, we  can apply Theorem \ref{aana} to get  a product expansion for $Z (\mathbf{s}^0 + t \bs{\theta})$ whenever $\Re(t) \geq \delta$, for any $\delta > 0.$

 We now state a lemma that will be used several times in the following discussion. 
The proof of this classic result follows from the Weierstrass Preparation Theorem and  can be found in \cite{abh}.


\begin{lemme}\label{athmWeierstrass}
 Let
$f : U \to \mathbb{C}$ be a nonzero holomorphic function  defined on an open set $U\subseteq \mathbb{C}^{n}$. Then the zero locus $f^{-1}(0)$ has   empty interior inside  $\mathbb{C}^n$.
\end{lemme}

\paragraph{Definition of generalized polynomials.} 

 We will need in the following to enlarge in a certain sense the class of classical polynomials of two variables.

\begin{defi}
 We will say that $W(X,Y)$ is a generalized polynomial if for all  $X\in \mathbb{C}\setminus \mathbb{R}_{-}$ and $Y\in \mathbb{C}$ deprived of an half-line we have
$
\ds W(X,Y)=1+\sum_{j=1}^{r}c_j X^{\nu_j}Y^{\mu_j} \  (c_j \in \mathbb{C});
$
where for all $j\in \{1,\dots,r\}$, $\nu_j\in \mathbb{R}_{\geq 0}$ and $\mu_j \in \mathbb{Q}_{>0}$.
\end{defi}
 
Notice that contrary to the classical polynomials, the generalized polynomials  $W(X,Y)$ are a priori defined only for  $X,Y\in \mathbb{C}$ deprived of an half-line since it is necessary to be able to define a logarithm to define them.

Thus when we will speak about a generalized polynomial $W(X,Y)$, it will be understood throughout the remainder that we consider it for $X\in \mathbb{C}\setminus \mathbb{R}_{-}$ and for $Y\in \mathbb{C}$ deprived of an half-line of the form $e^{i\mathfrak{b}}\mathbb{R}_{+}$ with $\mathfrak{b}\in \mathbb{R}$.

%
%
%

 Now we first fix a triplet $\mathbf{\bs{\mu}} = (p, \boldsymbol{\tau}^0, \bs{\theta})$ of parameters. It will always be understood to be the case that $p$ is a prime number, $\bs{\theta}$ satisfies (\ref {ahyp_theta1}) and $\boldsymbol{\tau}^0 \in \mathbb{R}^n$.

Then we consider in particular the  generalized polynomial for $\bs{\sigma}^0$ satisfying (\ref{agenericite}):
$$
W_{\mathbf{\bs{\mu}},\bs{\sigma}^0} \left(X,Y\right) = 1+ a_1 p^{-i\boldsymbol{\tau}^0\cdot\bs{\alpha}_{\cdot 1}}X^{\boldsymbol{\sigma}^0\cdot\bs{\alpha}_{\cdot 1}}Y^{\bs{\theta}\cdot\bs{\alpha}_{\cdot 1}}+\cdots + a_r p^{-i\boldsymbol{\tau}^0\cdot\bs{\alpha}_{\cdot r}}X^{\boldsymbol{\sigma}^0\cdot\bs{\alpha}_{\cdot r}}Y^{\bs{\theta}\cdot\bs{\alpha}_{\cdot r}}.
$$
 Setting $\mathbf{s}^0 = \boldsymbol{\sigma}^0 + i \boldsymbol{\tau}^0$, it follows that $
\ds W_{\mathbf{\bs{\mu}},\bs{\sigma}^0} \left(p^{-1},p^{-t}\right) = h\left(p^{-s^0_1-t\theta_1},\cdots,p^{-s^0_n-t\theta_n}\right).
$

Notice that for a fixed $\bs{\mu}$, $W_{\mathbf{\bs{\mu}},\bs{\sigma}^0}\left(p^{-1},p^{-t}\right)$ is well defined for almost all  $t\in \Xi_{u,\eta}$ (i.e. for all $t\in \Xi_{u,\eta}$ except at most a finite number).
Indeed, $p^{-1}\in \mathbb{C}\setminus \mathbb{R}_{-}$ and if we choose the branch $e^{i\mathfrak{b}}\mathbb{R}_{+}$ and the corresponding determination of the logarithm to define $W_{\mathbf{\bs{\mu}},\bs{\sigma}^0}(X,Y)$ on $\mathbb{C}\setminus \mathbb{R}_{-}\times \mathbb{C}\setminus e^{i\mathfrak{b}}\mathbb{R}_{+}$, then we have $p^{-t}\in e^{i\mathfrak{b}}\mathbb{R}_{+}$ if and only if there exists  $k\in \mathbb{Z}$ such that
\begin{equation}\label{awelldefine}
\Im(t) = \frac{2 k\pi-\mathfrak{b}}{\log(p)};
\end{equation}
hence there is at most a finite number of such $t\in \Xi_{u,\eta}$ verifying (\ref{awelldefine}); which justifies the fact that  $W_{\mathbf{\bs{\mu}},\bs{\sigma}^0}\left(p^{-1},p^{-t}\right)$ is well defined for almost all $t\in \Xi_{u,\eta}$.

%

Now the idea is to use one of the tools which have been developped by M. du Sautoy when he studied Euler products associated to a polynomial of two variables (\cite{sautoy}).
He uses in particular the Puiseux series theory, a generalization of the implicit functions theorem.
However, the classical Puiseux theory fits badly for the class of generalized polynomials since the classical Puiseux algorithm applied to a generalized polynomial don't provide a priori a convergent solution or even a formal solution.
In addition, the using of the multivariable Puiseux theory don't provide satisfactory results because of a very bad control of the domain of convergence of the solutions.

This is the reason why we use the hypothesis of Definition \ref{anouvcond2} which permits to reduce the problem and to apply the implicit functions theorem by considering these generalized polynomials as multivariable polynomials whose  variables are specialized; in this way we justify the existence and the convergence of the solutions for $X\in \mathbb{C}\setminus \mathbb{R}_{-}, |X|$ in the neighborhood of $0$.

%
%
%
%
%

Precisely, we want to express, for $X$ in a neighborhood of $0$, $Y$ as a function of $X$ such that $W_{\mathbf{\bs{\mu}},\bs{\sigma}^0}\left(X,Y\right)=0$.
However, $W_{\mathbf{\bs{\mu}},\bs{\sigma}^0}\left(X,Y\right)$ is not a real polynomial; and we must justify why it is possible to generalize the Puiseux theory (that we could find in \cite{c_alg_pl}) for the polynomials of two variables to this class of generalized polynomials. The main difficulty here comes from the fact that the exponents of the monomials of  $W_{\mathbf{\bs{\mu}},\bs{\sigma}^0}(X,Y)$ are not integers; hence $W_{\mathbf{\bs{\mu}},\bs{\sigma}^0}(X,Y)$ is defined from a determination of a logarithm and does not defined a regular function in  $X=0$.
\paragraph{Puiseux theorem for $W_{\mathbf{\bs{\mu}},\bs{\sigma}^0}\left(X,Y\right)$} 

\begin{defi}
 We put $\boldsymbol{\widehat{\alpha}}_{\cdot e}\in \mathbb{N}^{*}$ the vector collinear with  $\boldsymbol{\alpha}_{\cdot e}$ whose nonzero components are relatively prime.
In this way, if $j\in \Lambda_e$, there exists $q_j\in \mathbb{N}^{*}$ such that $\boldsymbol{\alpha}_{\cdot j} = q_j \boldsymbol{\widehat{\alpha}}_{\cdot e}$.
Then we put $
\ds \widetilde{[h]_e}(T):=1+\sum_{j\in \Lambda_e}a_j T^{q_j}.
$
\end{defi}

Since $\mathcal{F}(\boldsymbol{\alpha}_{\cdot e})$ is supposed to be a non-degenerate face, the polynomial $\widetilde{[h]_e}(T)$ has no multiple root.
For the following we will suppose that $\bs{\theta}\cdot \bs{\widehat{\alpha}}_{\cdot e}$ is a positive integer (we will precise later (see  (\ref{achoixjudicieuxtheta})) that we choose $\bs{\theta}$ so that $\bs{\theta}\cdot \bs{\widehat{\alpha}}_{\cdot e}$ is a positive even integer).

\begin{comment}
 The analytic behaviour of these generalized polynomials $W_{\mathbf{\bs{\mu}},\bs{\sigma}^0}\left(X,Y\right)$ being much more complex than for classical polynomials, this is where we will recourse to the assumption of non-degeneracy of the face $\mathcal{F}(\bs{\alpha}_{\cdot e})$. In addition, this hypothesis permits to have satisfactory results concerning the convergence of the Puiseux branches of $W_{\mathbf{\bs{\mu}},\bs{\sigma}^0}\left(X,Y\right)$.
Note in particular that the domain of convergence $\mathcal{H}$ of these branches obtained in Proposition \ref{athpuiseuxpolgen} does not depend on the parameter $p$, which is essential since we want to express the zeroes of $t\longmapsto h\left(p^{-s^0_1-t\theta_1},\dots, p^{-s^0_n-t\theta_n}\right)$ for infinitely many prime numbers $p$.
\end{comment}

\begin{prop}\label{athpuiseuxpolgen} {\bf Puiseux theorem for $W_{\mathbf{\bs{\mu}},\bs{\sigma}^0}\left(X,Y\right)$}.\\

We fix the parameter vector $\bs{\mu} = (p, \boldsymbol{\tau}^0, \bs{\theta})$.
Let $q\in \mathbb{N}^*$ be the smallest positive integer verifying $q \bs{\theta} \cdot \bs{\alpha}_{\cdot j} \in \mathbb{N}^*$ for all $j=1,\dots,r$.
 We consider the finite set: $$
\ds \mathfrak{r}_{\bs{\mu}}:=\left\{ c_{\bs{\mu}} \in \mathbb{C};~ \exists c \ \textrm{root of} \ \widetilde{[h]_e}(T) \ \textrm{such that} 
\ c_{\bs{\mu}}^{q \bs{\theta}\cdot \bs{\widehat{\alpha}}_{\cdot e}} = p^{i\bs{\tau}^0\cdot \bs{\widehat{\alpha}}_{\cdot e}} c  \right\}.
$$
There exists $\epsilon_1>0$ (not depending on $p$ nor on $\bs{\tau}^0$) such that  for all $X\in \mathcal{H}:=\{X\in \mathbb{C}\setminus \mathbb{R}_{-}, |X|<\epsilon_1\}$ the equation $W_{\bs{\mu}, \bs{\sigma}^0}(X,Y)=0$ admits the set of solutions $
\ds Y = \Omega_{\bs{\mu}, c_{\bs{\mu}},\bs{\sigma}^0}(X) \ (c_{\bs{\mu}} \in \mathfrak{r}_{\bs{\mu}});
$
where for all $ c_{\bs{\mu}} \in \mathfrak{r}_{\bs{\mu}}$ $X\longmapsto \Omega_{\bs{\mu}, c_{\bs{\mu}},\bs{\sigma}^0}(X)$ is an holomorphic function on  $\mathcal{H}$ and satisfies 
$\ds \Omega_{\bs{\mu}, c_{\bs{\mu}},\bs{\sigma}^0}(X) = \sum_{k=0}^{\kappa_{c_{\bs{\mu}}}}\mathfrak{c}_{k}(c_{\bs{\mu}},\bs{\mu})X^{\vartheta(\bs{\sigma}^0)_k},
$
with:
\begin{enumerate}
 \item $\kappa_{c_{\bs{\mu}}}\in \mathbb{N}\cup\{+\infty\}$;
\item $\vartheta(\bs{\sigma}^0)_0=0<\vartheta(\bs{\sigma}^0)_1<\cdots$ is a stricly increasing sequence not depending neither on p prime nor on $\bs{\tau}^0\in \mathbb{R}^n$;
\item $\lim_{k\to +\infty} \vartheta(\bs{\sigma}^0)_k = +\infty$ if $\kappa_{c_{\bs{\mu}}}=+\infty$;
\item there exists two constants $D_{\epsilon_0} >1$ and $A(\boldsymbol{\sigma}^0) >0$ (independing of $p$, $\bs{\tau}^0$ and $k$) such that 
$
\ds \left|\mathfrak{c}_k (c_{\bs{\mu}},\bs{\mu})\right|\ll  D_{\epsilon_0}^{A(\boldsymbol{\sigma}^0)\vartheta(\bs{\sigma}^0)_k}
$
uniformly in $p$ prime and in $k$;
\item $\mathfrak{c}_{0}(c_{\bs{\mu}},\bs{\mu}) = c_{\bs{\mu}}^q$, in particular
$\ds\left|\mathfrak{c}_{0}(c_{\bs{\mu}},\bs{\mu})\right| = |c|^{\frac{1}{\bs{\theta}\cdot \bs{\widehat{\alpha}}_{\cdot e}}}.
$
\end{enumerate}
Moreover $\ds\{ \mathfrak{c}_{0}(c_{\bs{\mu}},\bs{\mu}); c_{\bs{\mu}} \in \mathfrak{r}_{\bs{\mu}} \}
=\{ u \in \mathbb{C}; ~ \exists c \ \textrm{root of} \ \widetilde{[h]_e}(T) \ \textrm{such that}
~ u^{q \bs{\theta}\cdot \bs{\widehat{\alpha}}_{\cdot e}} = p^{i\bs{\tau}^0\cdot \bs{\widehat{\alpha}}_{\cdot e}} c \}.$
\end{prop}

\begin{dem}
We make the change of variable $Y=Y_1^q$.
The problem is reduced to resolve the equation
\begin{equation}\label{aexpreW}
\begin{array}{llll}
\ds {\tilde W_{\mathbf{\bs{\mu}},\bs{\sigma}^0}}(X,Y_1)&:=&W_{\mathbf{\bs{\mu}},\bs{\sigma}^0}(X,Y_1^q)\\
& = & 1+\sum_{j=1}^{r}a_j p^{-i\bs{\tau}^0 \cdot \bs{\alpha}_{\cdot j}}X^{\bs{\sigma}^0 \cdot  \bs{\alpha}_{\cdot j}}
Y_1^{q \bs{\theta} \cdot \bs{\alpha}_{\cdot j}} \\ 
 & = &\ds 1+\sum_{j\in \Lambda_e}a_j p^{-i\bs{\tau}^0 \cdot \bs{\alpha}_{\cdot j}}Y_1^{q \bs{\theta} \cdot \bs{\alpha}_{\cdot j}} + 
 \sum_{j\notin \Lambda_e}a_j p^{-i\bs{\tau}^0 \cdot \bs{\alpha}_{\cdot j}}X^{\bs{\sigma}^0 \cdot  \bs{\alpha}_{\cdot j}}
 Y_1^{q \bs{\theta} \cdot \bs{\alpha}_{\cdot j}} \\ 
 & = &\ds 1+\sum_{j\in \Lambda_e}a_j p^{-i q_j\bs{\tau}^0 \cdot \bs{\alpha}_{\cdot e}}
 Y_1^{q q_j \bs{\theta}\cdot \bs{\widehat{\alpha}}_{\cdot e}} + \sum_{j\notin \Lambda_e}a_j p^{-i\bs{\tau}^0 \cdot \bs{\alpha}_{\cdot j}}X^{\bs{\sigma}^0 \cdot  \bs{\alpha}_{\cdot j}}Y_1^{q \bs{\theta} \cdot \bs{\alpha}_{\cdot j}}.
\end{array}
\end{equation}

Now if $Y_1=Y_1(X)=c_{\bs{\mu}}+o(1)$ for $X\rightarrow 0^{+}$ is a solution of  $W_{\mathbf{\bs{\mu}},\bs{\sigma}^0}(X,Y_1)=0$ when $X\rightarrow 0^{+}$, then we have necessarily $
\ds 1+\sum_{j\in \Lambda_e}a_j p^{-i q_j \bs{\tau}^0 \cdot \bs{\alpha}_{\cdot e}}
c_{\bs{\mu}}^{q q_j \bs{\theta}\cdot \bs{\widehat{\alpha}}_{\cdot e}}=0;
$
and hence 
$
\ds c = p^{-i\bs{\tau}^0 \cdot \bs{\alpha}_{\cdot e}}c_{\bs{\mu}}^{q \bs{\theta}\cdot \bs{\widehat{\alpha}}_{\cdot e}}
$
is a root of the polynomial $\ds\widetilde{[h]_e}(T) = 1+\sum_{j\in \Lambda_e}a_j T^{q_j}$. We deduce that 
$c_{\bs{\mu}}$ is a $q \bs{\theta}\cdot \bs{\widehat{\alpha}}_{\cdot e}$-th root of  
$p^{i\bs{\tau}^0\cdot \bs{\widehat{\alpha}}_{\cdot e}} c$  where 
$c$ is a root of the polynomial $\ds\widetilde{[h]_e}(T)$. 
Thus $c_{\bs{\mu}} \in \mathfrak{r}_{\bs{\mu}}$.

Reciprocally, let $c_{\bs{\mu}} \in \mathfrak{r}_{\bs{\mu}}$. Then there exists a root $c$ of the polynomial
$\ds\widetilde{[h]_e}(T)$ such that $c_{\bs{\mu}}^{q\bs{\theta}\cdot \bs{\widehat{\alpha}}_{\cdot e}}=p^{i\bs{\tau}^0\cdot \bs{\widehat{\alpha}}_{\cdot e}} c$.
The root $c$ of $\widetilde{[h]_e}(T)$ is necessarily nonzero since  $\widetilde{[h]_e}(0)=1\neq 0$.
We make the change of variable
$
\ds Y_1 = c_{\bs{\mu}}\left(1+Y_2\right).
$
We search $Y_2=Y_2(X)$ such that $\ds {\tilde W_{\mathbf{\bs{\mu}},\bs{\sigma}^0}}(X,c_{\bs{\mu}}(1+Y_2(X)))=0$ and $Y_2(X)\rightarrow 0$ when $X\rightarrow 0^{+}$.
We put $
\ds G(X,Y_2):={\tilde W_{\mathbf{\bs{\mu}},\bs{\sigma}^0}}(X,c_{\bs{\mu}}(1+Y_2)).
$
We have:

\begin{displaymath}
\begin{array}{lll}
G(X,Y_2) & = &\ds 1+\sum_{j\in \Lambda_e}a_j p^{-i q_j \bs{\tau}^0 \cdot \bs{\alpha}_{\cdot e}}c_{\bs{\mu}}^{q q_j \bs{\theta}\cdot \bs{\widehat{\alpha}}_{\cdot e}}\left(1+Y_2\right)^{q q_j \bs{\theta}\cdot \bs{\widehat{\alpha}}_{\cdot e}}  \\ 
 &  &\ds + \sum_{j\notin\Lambda_e}a_j p^{-i\bs{\tau}^0 \cdot \bs{\alpha}_{\cdot j}}X^{\bs{\sigma}^0 \cdot  \bs{\alpha}_{\cdot j}}c_{\bs{\mu}}^{q \bs{\theta} \cdot \bs{\alpha}_{\cdot j}}(1+Y_2)^{q \bs{\theta} \cdot \bs{\alpha}_{\cdot j}} \\ 
 & = &\ds 1+\sum_{j\in \Lambda_e}a_j c^{q_j} (1+Y_2)^{q q_j \bs{\theta}\cdot \bs{\widehat{\alpha}}_{\cdot e}}  \\
 &  & \ds +\sum_{j\notin\Lambda_e}a_j p^{-i\bs{\tau}^0 \cdot \bs{\alpha}_{\cdot j}} c_{\bs{\mu}}^{q \bs{\theta}\cdot \bs{\alpha}_{\cdot j}}X^{\bs{\sigma}^0 \cdot  \bs{\alpha}_{\cdot j}}\left(1+Y_2\right)^{q \bs{\theta} \cdot \bs{\alpha}_{\cdot j}}   
\end{array}
\end{displaymath}

Since the $\bs{\theta} \cdot \bs{\alpha}_{\cdot j}\in \mathbb{Q}_{>0}$  and the $\bs{\sigma}^0 \cdot  \bs{\alpha}_{\cdot j}\in \mathbb{R}_{>0}$  ($j\notin \Lambda_e$) and by the choice of $q$, $G(X,Y_2)$ is defined and holomorphic on  $\mathcal{D}_1:=\mathbb{C}\setminus \mathbb{R}_{-}\times \mathbb{C}$.
Put for all $\textbf{X}=(X_j)_{j\notin \Lambda_e}\in \mathbb{C}^{r-\#\Lambda_e}$ and for all $Y_2\in  \mathbb{C}$:
$$
F(\textbf{X},Y_2):= 1+\sum_{j\in \Lambda_e}a_j c^{q_j}(1+Y_2)^{q q_j \bs{\theta}\cdot \bs{\widehat{\alpha}}_{\cdot e}} + 
\sum_{j\notin \Lambda_e}a_j X_j (1+Y_2)^{q \bs{\theta} \cdot \bs{\alpha}_{\cdot j}}.
$$
This function $(\textbf{X},Y_2)\longmapsto F(\textbf{X},Y_2)$ is clearly holomorphic on  $\mathcal{U} = \mathbb{C}^{r-\#\Lambda_e}\times \mathbb{C}$ and we also notice that this polynomial does not depend neither on $p$ nor on $\bs{\tau}^0$.
Moreover, for all $(X,Y_2)\in \mathcal{D}_1$ we have:
\begin{equation}\label{apuispecialis}
G(X,Y_2) = F\left(\left(p^{-i\bs{\tau}^0 \cdot \bs{\alpha}_{\cdot j}} c_{\bs{\mu}}^{q \bs{\theta}\cdot \bs{\alpha}_{\cdot j}}X^{\bs{\sigma}^0 \cdot  \bs{\alpha}_{\cdot j}}\right)_{j\notin \Lambda_e},Y_2\right).
\end{equation}
Furthermore, we easily check that in the neighborhood of  $(\textbf{X},Y_2)=(\textbf{0},0)$ we have:
\begin{displaymath}
\begin{array}{lll}
F(\textbf{X},Y_2) & = &\ds 1+\sum_{j\in \Lambda_e}a_j c^{q_j}\left(1+q q_j \bs{\theta}\cdot \bs{\widehat{\alpha}}_{\cdot e} Y_2 + O(Y_2^2)\right) + O\left(\lVert X\lVert\right) \\ 
 & = &\ds 1+\sum_{j\in \Lambda_e}a_j c^{q_j} + q \left(\bs{\theta}\cdot \bs{\widehat{\alpha}}_{\cdot e}\right)\left(\sum_{j\in \Lambda_e}a_j q_j c^{q_j}\right)Y_2 + O\left(\lVert X\lVert\right) + O\left(Y_2^2\right) \\ 
 & = &\ds \widetilde{[h]_e}(c) + q c  \left(\bs{\theta}\cdot\bs{\widehat{\alpha}}_{\cdot e}\right) \widetilde{[h]'_e}(c) Y_2 + O\left(\lVert X\lVert\right) + O\left(Y_2^2\right).
\end{array}
\end{displaymath}

And since the face $\mathcal{F}(\boldsymbol{\alpha}_{\cdot e})$ is a non-degenerate face by hypothesis, $c$ is a simple root of  $\widetilde{[h]_e}(T)$, and hence $\widetilde{[h]_e}(c)=0$ and $\widetilde{[h]'_e}(c)\neq 0$.
We deduce that $
\ds \frac{\partial F}{\partial Y_2}(\textbf{0},0) = q c \left( \bs{\theta}\cdot \bs{\widehat{\alpha}}_{\cdot e} \right) \widetilde{[h]'_e}(c)\neq 0.
$
Hence, according to the implicit functions theorem, there exists  $\epsilon_0 = \epsilon_0(h)>0$ (independent of $p$) such that
\begin{displaymath}
 \begin{cases}
  \textbf{X}\in D(\textbf{0},\epsilon_0) = \{\textbf{X}\in \mathbb{C}^{r-\#\Lambda_e} \mid \lVert X\lVert<\epsilon_0\} \\
  F(\textbf{X},Y_2)=0
 \end{cases}
\end{displaymath}
is equivalent to
\begin{equation}\label{athmftimpl}
 \begin{cases}
  \textbf{X}\in D(\textbf{0},\epsilon_0) = \{\textbf{X}\in \mathbb{C}^{r-\#\Lambda_e} \mid \lVert X\lVert<\epsilon_0\} \\
  Y_2 = V(\textbf{X});
 \end{cases}
\end{equation}
where
\begin{displaymath}
\begin{array}{lll}
V(\textbf{X}) & = &\ds \sum_{(\nu_j)_{j\notin \Lambda_e}=\bs{\nu}\in \mathbb{N}^{r-\#\Lambda_e}}A(\bs{\nu})\textbf{X}^{\bs{\nu}} \\ 
 & = & \ds \sum_{(\nu_j)_{j\notin \Lambda_e}=\bs{\nu}\in \mathbb{N}^{r-\#\Lambda_e}}A(\bs{\nu})\prod_{j\notin \Lambda_e}X_j^{\nu_j}
\end{array} 
\end{displaymath}
converges and is holomorphic in $D(\textbf{0},\epsilon_0)$.
In particular, we have uniformly in  $p$ prime and $\bs{\tau}^0\in \mathbb{R}^n$:
\begin{equation}\label{aestimA}
\left|A(\bs{\nu})\right| \ll  \left(\frac{2}{\epsilon_0}\right)^{\sum_{j\notin \Lambda_e}\nu_j}.
\end{equation}
Now since $\bs{\sigma}^0 \cdot  \bs{\alpha}_{\cdot j}>0$ for all $j\notin \Lambda_e$ and since 
$|p^{-i\bs{\tau}^0 \cdot \bs{\alpha}_{\cdot j}} c_{\bs{\mu}}^{q \bs{\theta}\cdot \bs{\alpha}_{\cdot j}}|=|c|^{\frac{\bs{\theta}\cdot \bs{\alpha}_{\cdot j}}{\bs{\theta}\cdot \bs{\widehat{\alpha}}_{\cdot e}}}$ for all $j\notin \Lambda_e$, the identity (\ref{apuispecialis}) implies the existence of  $\epsilon_1 = \epsilon_1(h,\boldsymbol{\sigma}^0)>0$ (independent of $p$ and $\bs{\tau}^0$) such that for all $X\in \mathcal{H}:=\{X\in \mathbb{C}\setminus \mathbb{R}_{-} : |X|<\epsilon_1\}$, the point $
\ds {\mathbf X}(\bs{\mu}) =\left(p^{-i\bs{\tau}^0 \cdot \bs{\alpha}_{\cdot j}} c_{\bs{\mu}}^{q \bs{\theta}\cdot \bs{\alpha}_{\cdot j}}X^{\bs{\sigma}^0 \cdot  \bs{\alpha}_{\cdot j}}\right)_{j\notin \Lambda_e}\in D(\textbf{0},\epsilon_0);
$
and
$
\ds V({\mathbf X}(\bs{\mu})) = \sum_{\bs{\nu}=(\nu_j)_{j\notin \Lambda_e}}A(\bs{\nu})
p^{-i\sum_{j\notin \Lambda_e}\nu_j \bs{\tau}^0 \cdot \bs{\alpha}_{\cdot j}} c_{\bs{\mu}}^{q \sum_{j\notin \Lambda_e}\nu_j\bs{\theta}\cdot \bs{\alpha}_{\cdot j}} X^{\sum_{j\notin \Lambda_e}\nu_j \bs{\sigma}^0 \cdot  \bs{\alpha}_{\cdot j}}.
$

Put $
\ds\mathfrak{K}:=\left\{ \sum_{j\notin \Lambda_e}\nu_j \bs{\sigma}^0 \cdot  \bs{\alpha}_{\cdot j} : \bs{\nu}=(\nu_j)_{j\notin \Lambda_e}\in \mathbb{N}^{r-\#\Lambda_e} \right\}
$
which is a discrete part of $\mathbb{R}^{*}_{+}$. So there exists a strictly increasing sequence  (finite or infinite) of $\mathbb{R}^{*}_{+}$ such that $
\ds\mathfrak{K} = \{\vartheta(\bs{\sigma}^0)_k: k=1,\dots,\kappa_{c_{\bs{\mu}}}\}
$
with $\kappa_{c_{\bs{\mu}}}\in \mathbb{N}\cup \{+\infty\}$.
Moreover, it is clear that if $\kappa_{c_{\bs{\mu}}}=+\infty$ then $\lim_{k\to+\infty}\vartheta(\bs{\sigma}^0)_k = +\infty$.
And if we put for all $k\leq \kappa_{c_{\bs{\mu}}}$: 
$$
\ds\mathfrak{c}_{k}(c_{\bs{\mu}},\bs{\mu}) = \sum_{\stackrel{\bs{\nu}=(\nu_j)_{j\notin \Lambda_e}\in \mathbb{N}^{r-\#\Lambda_e}}{\sum_{j\notin \Lambda_e}\nu_j(\bs{\sigma}^0 \cdot  \bs{\alpha}_{\cdot j}) = \vartheta(\bs{\sigma}^0)_k}}A(\bs{\nu})
p^{-i\sum_{j\notin \Lambda_e}\nu_j \bs{\tau}^0 \cdot \bs{\alpha}_{\cdot j}} c_{\bs{\mu}}^{q \sum_{j\notin \Lambda_e}\nu_j\bs{\theta}\cdot \bs{\alpha}_{\cdot j}};
$$
the fact that $\nu_j\leq \frac{\vartheta(\bs{\sigma}^0)_k}{\bs{\sigma}^0 \cdot  \bs{\alpha}_{\cdot j}}$ for all $j\notin \Lambda_e$ since $\sum_{j\notin \Lambda_e}\nu_j\bs{\sigma}^0 \cdot  \bs{\alpha}_{\cdot j}=\vartheta(\bs{\sigma}^0)_k$ and  (\ref{aestimA}) provide the estimation:
\begin{equation}\label{aestimB}
\left|\mathfrak{c}_{k}(c_{\bs{\mu}},\bs{\mu})\right| \ll \left(\frac{2}{\epsilon_0}\right)^{\left(\sum_{j\notin \Lambda_e}\frac{1}{\bs{\sigma}^0 \cdot  \bs{\alpha}_{\cdot j}}\right)\vartheta(\bs{\sigma}^0)_k}\prod_{j\notin \Lambda_e}\frac{\vartheta(\bs{\sigma}^0)_k}{\bs{\sigma}^0 \cdot  \bs{\alpha}_{\cdot j}}
\end{equation}
uniformly in $p$ prime and in $\bs{\tau}^0\in \mathbb{R}^n$.
Consequently we have for all $X\in \mathcal{H} = \{X\in \mathbb{C}\setminus \mathbb{R}_{-}: |X|<\epsilon_1\}$ $
\ds V({\mathbf X}(\bs{\mu})) = \sum_{k=1}^{\kappa_{c_{\bs{\mu}}}}\mathfrak{c}_k(c_{\bs{\mu}},\bs{\mu})X^{\vartheta(\bs{\sigma}^0)_k}.
$
We conclude using the fact that $Y=\left(c_{\bs{\mu}}(1+Y_2)\right)^q = c_{\bs{\mu}}^q (1+Y_2)^q$.
\CQFD
\end{dem}
\begin{rqs}
Notice that the solutions obtained in  $Y$ of $W_{\mathbf{\bs{\mu}},\bs{\sigma}^0}(X,Y)=0$ appear as $q$-th powers.
Thus there will be no problem in the following related to the manipulation of rational powers $\bs{\theta} \cdot \bs{\alpha}_{\cdot j}$ in $Y$ of $W_{\mathbf{\bs{\mu}},\bs{\sigma}^0}(X,Y)$ when we will replace $Y$ by these solutions. 
\end{rqs}
To simplify the writing introduced in the previous proposition, in the whole following we will write any solution of  $W_{\mathbf{\bs{\mu}},\bs{\sigma}^0}(X,Y)=0$ (in finite number) as follows:
$$
\Omega_{\bs{\mu},\bs{\sigma}^0}(X) = c_{\bs{\mu},0} + c_{\bs{\mu},1} X^{\vartheta_{1}} + \cdots + c_{\bs{\mu},N} X^{\vartheta_{N}} + o\left(X^{\vartheta_{N}}\right), \ (N\geq 1); 
$$
where $\vartheta_{N}=\vartheta(\bs{\sigma}^0)_{N}>\cdots,>\vartheta_{1}=\vartheta(\bs{\sigma}^0)_{1}>0$ and $c_{\bs{\mu},m}\in \mathbb{C}$ for all $m\geq 0$.


In particular  $c_{\bs{\mu},0}$ is a root of the one variable polynomial $\ds [W_{\mathbf{\bs{\mu}},\bs{\sigma}^0}]_e(y):=1+\sum_{j\in \Lambda_e}a_j p^{-i\bs{\tau}^0 \cdot \bs{\alpha}_{\cdot j}}y^{\bs{\theta} \cdot \bs{\alpha}_{\cdot j}}.$
We can also notice that if we  set
\begin{equation}\label{adepenpcko}
\tilde{c}_{\mathbf{\bs{\theta}},0} = c_{\mathbf{\bs{\mu}},0} p^{-i\frac{\bs{\tau}^0 \cdot \bs{\alpha}_{\cdot e}}{\bs{\theta}\cdot \bs{\alpha}_{\cdot e}}},
\end{equation}
 then $|\tilde{c}_{\mathbf{\bs{\theta}},0}| = |c_{\mathbf{\bs{\mu}},0}|$ and $\tilde{c}_{\mathbf{\bs{\theta}},0}$ is a root of the polynomial  
$1+\sum_{j \in \Lambda_e}a_j y^{\bs{\theta} \cdot \bs{\alpha}_{\cdot j}}$.


Thus we can describe the zeroes of $t \longmapsto W_{\mathbf{\bs{\mu}},\bs{\sigma}^0}(p^{-1},p^{-t})$ via these Puiseux branches; they can be expressed as follows:
$$
t_{m,\bs{\mu},\bs{\sigma}^0} = -\frac{\log\left(\Omega_{\bs{\mu},\bs{\sigma}^0}\left(p^{-1}\right)\right)}{\log(p)} + \frac{2 \pi m i}{\log(p)};
$$
where $m\in \mathbb{Z}$ and $p$ is a prime number large enough.

%
%

\begin{comment}
 In the following, we want to prove that the branches $\Omega_{\bs{\mu},\bs{\sigma}^0}\left(X\right)$ satisfy some properties, in particular the fact there exists at least one branch such that $\left|\Omega_{\bs{\mu},\bs{\sigma}^0}\left(X\right)\right|<1$ for $X>0$ small enough in a way to have $\Re\left(t_{m,\bs{\mu},\bs{\sigma}^0}\right)>0$ for $p$ large enough.

But to do this, it is necessary to select certain values of the parameters $\bs{\sigma}^0, \bs{\tau}^0$ and $\bs{\theta}$ on which depend the branches $\Omega_{\bs{\mu},\bs{\sigma}^0}\left(X\right)$. 

For the parameters $\bs{\sigma}^0$ and $\bs{\tau}^0$, the major difficulty will be to prove that their possible values that provide the desired properties for $\Omega_{\bs{\mu},\bs{\sigma}^0}\left(X\right)$ are generic, that is are so many that it is always possible, moving the point $\mathbf{s}^0$ inside the ball $\mathcal{B}$ if necessary, to assume that $\mathbf{s}^0=\bs{\sigma}^0 +i\bs{\tau}^0$ with ''good`` values for $\bs{\sigma}^0$ and $\bs{\tau}^0$.

Concerning the direction $\bs{\theta}$, it will be necessary to guarantee the existence of at least one value of this parameter for which we are able to prove the desired properties for $\Omega_{\bs{\mu},\bs{\sigma}^0}\left(X\right)$ provided that the other parameters $\bs{\sigma}^0$ and $\bs{\tau}^0$ belong to the generic set of ''good`` values.

The following lemma \ref{aW-noncyclo} is a good illustration of this principle  which consists in integrating the dynamic aspect of parameters.
\end{comment}

The following technical result will be useful later to show the existence of a particular index $e'\in \{1,\dots,r\}\setminus \{e\}$ (see (\ref{particularindex})) that will be important when we will need to calculate the second term of a Puiseux branch.

\begin{lemme}\label{aW-noncyclo}
Fix a prime number $p$ and a direction $\bs{\theta}$ and consider the parameter vector $\mathbf{\bs{\mu}} = (p, \boldsymbol{\tau}^0, \bs{\theta})$ for some $\bs{\tau}^0\in \mathbb{R}^n$.
Suppose that there exists $\mathfrak{c}\in \mathbb{C}\setminus e^{i\mathfrak{b}}\mathbb{R}_{+}$ such that for all $X\in \mathbb{C}\setminus \mathbb{R}_{-}$ and for all $\bs{\sigma}^0$ in a generic set so that $\mathbf{s}^0\in \mathcal{B}\cap \partial \mathbf{W}(0)$  we have $
\ds W_{\mathbf{\bs{\mu}},\bs{\sigma}^0}\mid_{\boldsymbol{\tau}^0=\textbf{0}}\left(X,\mathfrak{c}\right) = 0
$.
Then we have necessarily $\left|\mathfrak{c}\right|\neq 1$.

\end{lemme}

\begin{dem}

Let us assume that there exists $\mathfrak{c}\in \mathbb{C}^{*}$ such that for all $X\in \mathbb{C}\setminus \mathbb{R}_{-}$ $
\ds W_{\mathbf{\bs{\mu}},\bs{\sigma}^0}\mid_{\boldsymbol{\tau}^0=\textbf{0}}(X,\mathfrak{c})=0.
$ 

By hypothesis, there exists at least an open ball $U\subseteq\mathbb{R}^n$ such that for all $\bs{\sigma}^0\in U\cap\{\textbf{x}\in \mathbb{R}^n: \textbf{x}\cdot\bs{\alpha}_{\cdot e}=0\}$ and for all $X\in \mathbb{C}\setminus \mathbb{R}_{-}$ we have:
\begin{equation}\label{ahtnul}
\begin{array}{lll}
\displaystyle W_{\mathbf{\bs{\mu}},\bs{\sigma}^0}\mid_{\boldsymbol{\tau}^0=\textbf{0}}\left(X,\mathfrak{c}\right) & = & \displaystyle 1+\sum_{j=1}^{r}a_j X^{\bs{\sigma}^0 \cdot  \bs{\alpha}_{\cdot j}}\mathfrak{c}^{\bs{\theta} \cdot \bs{\alpha}_{\cdot j}}  
  =  \displaystyle h\left(X^{\boldsymbol{\sigma}^0_1}\mathfrak{c}^{\theta_1},\dots, X^{\boldsymbol{\sigma}^0_n}\mathfrak{c}^{\theta_n}\right) = 0. 
\end{array}
\end{equation}
Moreover, the only constraint of $\bs{\sigma}^0\in U\cap\{\textbf{x}\in \mathbb{R}^n: \textbf{x}\cdot\bs{\alpha}_{\cdot e}=0\}$ is that its components must verify  $\bs{\sigma}^0 \cdot \bs{\alpha}_{\cdot e} = 0$.
Consequently, since we have assumed (see (\ref{alphaenonnul})) without loss of generality that $\alpha_{n e} \neq 0$, we can consider $\boldsymbol{\sigma}^0\in U\cap\{\textbf{x}\in \mathbb{R}^n: \textbf{x}\cdot\bs{\alpha}_{\cdot e}=0\}$  as a $(n-1)$-uple $\widetilde{\bs{\sigma}}^0=(\widetilde{\bs{\sigma}}^0_1,\dots,\widetilde{\bs{\sigma}}^0_{n-1})$ in an open ball $\widetilde{U}\subseteq \mathbb{R}^{n-1}$  by putting:
\begin{equation*}
 \begin{cases}
 \displaystyle \boldsymbol{\sigma}^0_{\ell} = \widetilde{\bs{\sigma}}^0_{\ell} & (\ell \in \{1,\dots,n-1\}), \\
 \displaystyle \boldsymbol{\sigma}^0_n = -\frac{1}{\alpha_{n e}}\sum_{i=1}^{n-1}\alpha_{i e}\widetilde{\bs{\sigma}}^0_i & .
 \end{cases}
 \end{equation*}
Then define for all $x\in \mathbb{R}_{>0}$:
\begin{displaymath}
\begin{array}{cccc}
 & \widetilde{U} & \longrightarrow & \mathbb{R}^{n-1} \\ 
\Phi_{x}: & & & \\
 & \widetilde{\bs{\sigma}}^0=(\widetilde{\bs{\sigma}}^0_1,\dots,\widetilde{\bs{\sigma}}^0_{n-1}) & \longmapsto & \left(x^{\widetilde{\bs{\sigma}}^0_1},\dots,x^{\widetilde{\bs{\sigma}}^0_{n-1}}\right).
\end{array}
\end{displaymath}
It is clear that $\bigcup_{x>0}\Phi_{x}(\widetilde{U})$ describes a nonempty open set $\widetilde{U}'$ of $(0,\infty)^{n-1}$.
Consequently, for all $(t_1,\dots,t_{n-1})\in \widetilde{U}'$ there exists $x>0$ and $\widetilde{\bs{\sigma}}^0=(\widetilde{\bs{\sigma}}^0_1,\dots,\widetilde{\bs{\sigma}}^0_{n-1})\in \widetilde{U}$ such that $(t_1,\dots,t_{n-1}) = \left(x^{\widetilde{\bs{\sigma}}^0_1},\dots,x^{\widetilde{\bs{\sigma}}^0_{n-1}}\right)$ and we have:
\begin{displaymath}
\begin{array}{lll}
\displaystyle h\left(t_1 \mathfrak{c}^{\theta_1},\dots,t_{n-1} \mathfrak{c}^{\theta_{n-1}}, \mathfrak{c}^{\theta_{n}}\prod_{\ell=1}^{n-1}t_{\ell}^{-\frac{\alpha_{\ell e}}{\alpha_{n e}}}\right) & = & \displaystyle h\left(x^{\boldsymbol{\sigma}^0_1}\mathfrak{c}^{\theta_1},\dots, x^{\boldsymbol{\sigma}^0_n}\mathfrak{c}^{\theta_n}\right) \\ 
 & = & 0 \ \textrm{according to (\ref{ahtnul})}.
\end{array}
\end{displaymath}
Moreover, there exists a nonempty open set $\widetilde{U}''$ of $(0,\infty)^{n-1}$ such that for all $(y_1,\dots,y_{n-1})\in \widetilde{U}''$ there exists $(t_1,\dots,t_{n-1})\in \widetilde{U}'$ verifying for all $i\in \{1,\dots,n-1\}$ $
\ds y_i^{\alpha_{n e}}=t_i.
$
But the function $
\ds \left(y_1,\dots,y_{n-1}\right) \longmapsto \displaystyle h\left(y_1^{\alpha_{n e}} \mathfrak{c}^{\theta_1},\dots,y_{n-1}^{\alpha_{n e}} \mathfrak{c}^{\theta_{n-1}}, \mathfrak{c}^{\theta_{n}}\prod_{\ell=1}^{n-1}y_{\ell}^{-\alpha_{\ell e}}\right)
$ 
is holomorphic on $(\mathbb{C}^{*})^{n-1}$.
And since it vanishes on an open set $\widetilde{U}''$ de $(0,\infty)^{n-1}$, we have in fact
$
\ds\forall (y_1,\dots,y_{n-1})\in (\mathbb{C}^{*})^{n-1}, \ h\left(y_1^{\alpha_{n e}} \mathfrak{c}^{\theta_1},\dots,y_{n-1}^{\alpha_{n e}} \mathfrak{c}^{\theta_{n-1}}, \mathfrak{c}^{\theta_{n}}\prod_{\ell=1}^{n-1}y_{\ell}^{-\alpha_{\ell e}}\right)=0.
$ 
Hence the polynomial $X_1 \cdots X_{n-1} \ h(X_1,\dots,X_n)$ vanishes on $\mathcal{H}\cap (\mathbb{C}^{*})^{n-1}$ where $\mathcal{H}$ is the complex hypersurface defined by the equation $
\ds \mathfrak{c}^{-\bs{\theta}\cdot \bs{\alpha}_{\cdot e}}\textbf{X}^{\boldsymbol{\alpha}_{\cdot e}}-1 = \prod_{\ell=1}^{n}X_{\ell}^{\boldsymbol{\alpha}_{\cdot e}^{\ell}} \mathfrak{c}^{-\theta_{\ell}\boldsymbol{\alpha}_{\cdot e}^{\ell}} -1 =0.
$
We deduce that $X_1 \cdots X_{n-1}h(X_1,\dots,X_n)$ vanishes on the whole hypersurface $\mathcal{H}$ and hence the polynomial $\mathfrak{c}^{-\bs{\theta}\cdot \bs{\alpha}_{\cdot e}} \textbf{X}^{\boldsymbol{\alpha}_{\cdot e}} -1$
divides a power of the polynomial $X_1\dots X_{n-1} ~h(X_1,\dots,X_n)$. 
Since the polynomials $\mathfrak{c}^{-\bs{\theta}\cdot \bs{\alpha}_{\cdot e}} \textbf{X}^{\boldsymbol{\alpha}_{\cdot e}} -1$ and $X_1\dots X_{n-1}$ are relatively prime, we deduce that the polynomial  
\begin{equation}\label{adiviseurh}
P_{\mathfrak{c}}(\textbf{X}):= \mathfrak{c}^{-\bs{\theta}\cdot \bs{\alpha}_{\cdot e}}\textbf{X}^{\boldsymbol{\alpha}_{\cdot e}}-1
\end{equation}
 necessarily divides a power of $h$; and hence $P_{\mathfrak{c}}(\textbf{X})$ divides also $h$ because all irreducible factors of $P_{\mathfrak{c}}(\textbf{X})$ are of multiplicity $1$.
And since $h$ is with rational coefficients and  $\mathfrak{c}$ is an algebraic number, the polynomial $
\ds Q(\textbf{X}):=\prod_{\mathfrak{c}'}P_{\mathfrak{c}'}(\textbf{X}) \in \mathbb{Q}[\textbf{X}]
$
(where the product is done over all the conjugates $\mathfrak{c}'$ of $\mathfrak{c}$)
also divides $h$.
Remark that $Q(\textbf{X})$ can be reduce in fact to a one variable polynomial (by the change of variable $T:=\textbf{X}^{\boldsymbol{\alpha}_{\cdot e}}$). 
Moreover, if we assume by absurd that $|\mathfrak{c}|=1$, we could apply the criterion ii) of cyclotomy of  Estermann's result to the one variable polynomial  $Q(\textbf{X})$ to deduce that this polynomial is cyclotomic; which is not possible since  $h$ does not contain any cyclotomic factor by hypothesis.
\CQFD
\end{dem}

 Now the problem consists in finding an infinite number of  $t_{m,\bs{\mu},\bs{\sigma}^0}$ of positive real part. This reduces  to finding Puiseux branches  $\Omega_{\bs{\mu},\bs{\sigma}^0} (X)$ such that $|\Omega_{\bs{\mu},\bs{\sigma}^0} (X)|<1$ for $|X|$  small enough. For in such an event, it would follow  that for $p$ sufficiently large,  $t_{m,\bs{\mu},\bs{\sigma}^0} \in \Xi_{u,\eta}$. 
We will show in the following (see Lemma \ref{accumulation-zeros-h}) that we can always find  such a branch  $\Omega_{\bs{\mu},\bs{\sigma}^0}  (X)$ such that $\left|\Omega_{\bs{\mu},\bs{\sigma}^0}(X)\right|<1$ for $|X|$ small enough.
Note in passing that for $p$ large enough $W_{\mathbf{\bs{\mu}},\bs{\sigma}^0}(X,Y)$ is well defined by putting $X=p^{-1}$ and $Y=p^{-t_{m,\bs{\mu},\bs{\sigma}^0}}$ since $p^{-1}\in \mathbb{C}\setminus \mathbb{R}_{-}$ and $p^{-t_{m,\bs{\mu},\bs{\sigma}^0}}\in \mathbb{C}\setminus e^{i\mathfrak{b}}\mathbb{R}_{+}$ since:
\begin{displaymath}
\begin{array}{lll}
\Im(t_{m,\bs{\mu},\bs{\sigma}^0}) & = & \ds\Im\left(-\frac{\log\left(\Omega_{\bs{\mu},\bs{\sigma}^0}(p^{-1})\right)}{\log(p)}+\frac{2 i \pi m}{\log(p)}\right)  
  =  \ds\frac{-\arg(c_{\bs{\mu},0})+O\left(p^{-\vartheta_{k,1}}\right)+2 \pi m}{\log(p)} \\ 
 & = & \ds\frac{-\arg(\tilde{c}_{\mathbf{\bs{\theta}},0})+O\left(p^{-\vartheta_{k,1}}\right)+2 \pi m}{\log(p)} - \frac{\bs{\tau}^0 \cdot \bs{\alpha}_{\cdot e}}{\bs{\theta}\cdot \bs{\alpha}_{\cdot e}} \notin \frac{-\mathfrak{b} + 2 \mathbb{Z}\pi}{\log(p)} 
\end{array}
\end{displaymath}
if we choose $\boldsymbol{\tau}^0\in \mathbb{R}^n$ generically.

\vspace{0.2cm}

Now if the one variable polynomial $1+\sum_{j \in \Lambda_e}a_j y^{\bs{\theta} \cdot \bs{\alpha}_{\cdot j}}$ is {\it not } cyclotomic (which is equivalent to have  $[h]_e$ and $[W_{\mathbf{\bs{\mu}},\bs{\sigma}^0}]_e$ not cyclotomic), the fact that its coefficients are integers implies  there exists at least  one root $c^*$ whose norm is strictly less than $1$. As a result, it follows that there  also exists a Puiseux branch   $\Omega_{\bs{\mu},\bs{\sigma}^0}^*(X)$ with constant term $c_{\mathbf{\bs{\mu}},0}^*$ such that $|c_{\mathbf{\bs{\mu}},0}^*|<1$. Thus,   $|\Omega_{\bs{\mu},\bs{\sigma}^0}^* (X)|<1$ for $|X|$ sufficiently small, 
which is what we need to complete Part 1.

 The situation is therefore more complicated when: 
$$[W_{\mathbf{\bs{\mu}},\bs{\sigma}^0}]_e \quad \text{{\it is  cyclotomic.}}$$ 
{\it We will therefore assume this property   for the rest of this first part}. 
 In this event, it  will be necessary to study the second term which appears in any  Puiseux series  $\Omega_{\bs{\mu},\bs{\sigma}^0}(X)$ of $W_{\mathbf{\bs{\mu}},\bs{\sigma}^0}(X,Y)$.



Consider a set  $J$ made of representatives of each class of the following equivalence relation  $\sim$:
$$
\boldsymbol{\alpha}_{\cdot j} \sim \bs{\alpha}_{\cdot j'} \Longleftrightarrow \boldsymbol{\alpha}_{\cdot j}-\bs{\alpha}_{\cdot j'}\in \mathbb{Q}\boldsymbol{\alpha}_{\cdot e}.
$$

Then we write $
\ds W_{\mathbf{\bs{\mu}},\bs{\sigma}^0}(X,Y) = [W_{\mathbf{\bs{\mu}},\bs{\sigma}^0}]_e(Y) + \sum_{j_0\in J; j_0\not\sim e}X^{\bs{\sigma}^0\cdot\bs{\alpha}_{\cdot j_0}}R_{\bs{\mu},j_0}(Y)
$
with
$
\ds R_{\bs{\mu},j_0}(Y)=\sum_{j\sim j_0}a_j p^{-i\bs{\tau}^0 \cdot \bs{\alpha}_{\cdot j}} Y^{\bs{\theta} \cdot \bs{\alpha}_{\cdot j}}.
$
Recall that since we suppose here that  $[W_{\mathbf{\bs{\mu}},\bs{\sigma}^0}]_e(Y)$ is a cyclotomic polynomial, all its roots are of modulus  $1$.

Let $c_{\bs{\mu},0}$ be a root of $[W_{\mathbf{\bs{\mu}},\bs{\sigma}^0}]_e$ of multiplicity $m_0= 1$ since $\mathcal{F}(\boldsymbol{\alpha}_{\cdot e})$ is non-degenerate in the sense of Definition \ref{anouvcond2}.
Consider an index $e'\in \{1,\dots,r\}\setminus\{e\}$ satisfying the following property:

\begin{equation}\label{particularindex} 
\ds \bs{\sigma}^0\cdot \bs{\alpha}_{\cdot e'}>0 \
\textrm{is minimal among the} \ \bs{\sigma}^0\cdot\bs{\alpha}_{\cdot j_0}>0 \  (j_0\in J, j_0 \not\sim e) \  \textrm{such that} \  R_{\bs{\mu},j_0}(c_{\bs{\mu},0})\neq 0.
\end{equation}
\vspace{0.2cm}

This index $e'$ play a special role in the computation of the second term of a Puiseux branch of main term the root $c_{\bs{\mu},0}$.
Furthermore, remark that such an index $e'$ exists according to lemma \ref{aW-noncyclo}.
Indeed $
\ds R_{\bs{\mu},j_0}(c_{\bs{\mu},0}) = \sum_{j \sim j_0} a_j p^{-i\boldsymbol{\tau}^0\cdot\left(\boldsymbol{\alpha}_{\cdot j}-\boldsymbol{\alpha}_{\cdot e}\frac{\bs{\theta} \cdot \bs{\alpha}_{\cdot j}}{\bs{\theta}\cdot \bs{\alpha}_{\cdot e}}\right)} \left(\tilde{c}_{\mathbf{\bs{\theta}},0}\right)^{\bs{\theta} \cdot \bs{\alpha}_{\cdot j}} \ \textrm{according to}\ (\ref{adepenpcko}).
$
But if $\boldsymbol{\alpha}_{\cdot j} = \boldsymbol{\alpha}_{\cdot j_0}+q \boldsymbol{\alpha}_{\cdot e}$, we obtain
$\ds \boldsymbol{\alpha}_{\cdot j}-\boldsymbol{\alpha}_{\cdot e}\frac{\bs{\theta} \cdot \bs{\alpha}_{\cdot j}}{\bs{\theta}\cdot \bs{\alpha}_{\cdot e}}  =  \boldsymbol{\alpha}_{\cdot j_0} + q \boldsymbol{\alpha}_{\cdot e} -\boldsymbol{\alpha}_{\cdot e}\frac{\bs{\theta}\cdot \bs{\alpha}_{\cdot j_0}+q\bs{\theta}\cdot \bs{\alpha}_{\cdot e}}{\bs{\theta}\cdot \bs{\alpha}_{\cdot e}}  
  =  \boldsymbol{\alpha}_{\cdot j_0}-\boldsymbol{\alpha}_{\cdot e}\frac{\bs{\theta}\cdot \bs{\alpha}_{\cdot j_0}}{\bs{\theta}\cdot \bs{\alpha}_{\cdot e}}.$
Consequently since the $\boldsymbol{\alpha}_{\cdot j}-\boldsymbol{\alpha}_{\cdot e}\frac{\bs{\theta} \cdot \bs{\alpha}_{\cdot j}}{\bs{\theta}\cdot \bs{\alpha}_{\cdot e}}$ are all equal for  $j\sim j_0$  we have $
\ds R_{\bs{\mu},j_0}(c_{\bs{\mu},0}) = 0
$
is equivalent to 
\begin{equation}\label{aequationnulle?}
R_{j_0}(\tilde{c}_{\mathbf{\bs{\theta}},0}):=\sum_{j\sim j_0} a_j  \left(\tilde{c}_{\mathbf{\bs{\theta}},0}\right)^{\bs{\theta} \cdot \bs{\alpha}_{\cdot j}} = 0.
\end{equation}
Thus if $e'$ does not exist, we would have for all $X\in \mathbb{C}\setminus \mathbb{R}_{-}$ and for all $\bs{\sigma}^0\in \mathbb{R}^n$ satisfying (\ref{agenericite}), $
\ds W_{\mathbf{\bs{\mu}},\bs{\sigma}^0}\mid_{\boldsymbol{\tau}^0=\textbf{0}}(X,\tilde{c}_{\mathbf{\bs{\theta}},0}) = 1+\sum_{j\in \Lambda_e} a_j \left(\tilde{c}_{\bs{\theta},0}\right)^{\bs{\theta}\cdot \bs{\alpha}_{\cdot j}} + \sum_{j_0\in J; j_0\not\sim e} X^{\bs{\sigma}^0\cdot\bs{\alpha}_{\cdot j_0}}R_{j_0}(\tilde{c}_{\mathbf{\bs{\theta}},0})= [W_{\bs{\mu},\bs{\sigma}^0}]_e(c_{\bs{\mu},0})+ \sum_{j_0\in J; j_0\not\sim e} X^{\bs{\sigma}^0\cdot\bs{\alpha}_{\cdot j_0}}R_{j_0}(\tilde{c}_{\mathbf{\bs{\theta}},0})=0,
$
which is impossible according to Lemma \ref{aW-noncyclo} since here $|\tilde{c}_{\mathbf{\bs{\theta}},0}|=|c_{\bs{\mu},0}|=1$.
Obviously, it is possible to have some $j_0$ such that
\begin{equation}\label{aminimalindice}
\bs{\sigma}^0\cdot \bs{\alpha}_{\cdot e'} = \bs{\sigma}^0\cdot\bs{\alpha}_{\cdot j_0}.
\end{equation}
However, if ${\boldsymbol{\sigma}^0}\in \mathbb{R}^n$ is chosen generically so that ${\mathbf{s}^0}\in \mathcal{B}\cap \partial \mathbf{W}(0)$, the equality (\ref{aminimalindice}) implies necessarily that $j_0\sim e'$. 

\vspace{0.2cm}

 For the following it will be necessary to impose two supplementary conditions on $\bs{\theta}$ (in addition to (\ref{ahyp_theta1})) that depend upon our choice for the index $e'$: 
\begin{equation}\label{achoixjudicieuxtheta}
\begin{array}{l}
(a) \ \bs{\theta}\cdot \bs{\widehat{\alpha}}_{\cdot e} \in \mathbb{Z}_+ \ \text{{\it is even}};  \\ 
(b) \ \bs{\theta}\cdot \bs{\alpha}_{\cdot e'}  \in \mathbb{Z}_+ \ \text{\it {is odd}}.
\end{array}
\end{equation}
Remark that these two conditions are not in conflict because $\boldsymbol{\widehat{\alpha}}_{\cdot e}$ and $\boldsymbol{\alpha}_{\cdot e'}$ are not collinear (since $\bs{\sigma}^0\cdot \bs{\alpha}_{\cdot e'}>0$).

 Two other elementary properties will also be used in the following. 
\begin{enumerate}
\item  Since $\boldsymbol{\sigma}^0$ satisfies  the generic  condition  (\ref {agenericite}) it follows that 
$
\ds \bs{\sigma}^0\cdot\bs{\alpha}_{\cdot j_0} = \bs{\sigma}^0\cdot \bs{\alpha}_{\cdot e'} \quad \text{{\it implies}} \quad  \boldsymbol{\alpha}_{\cdot j_0} \in \boldsymbol{\alpha}_{\cdot e'} + \langle \boldsymbol{\alpha}_{\cdot e} \rangle\,.$ 
Thus if $j_0 \sim e',$ then $\bs{\theta}\cdot \bs{\alpha}_{\cdot j_0}\in \mathbb{Z}_{+}$ is  odd.
\item Since $j\in \Lambda_e$ implies  $\bs{\theta} \cdot \bs{\alpha}_{\cdot j}$ is an   even integer, it follows that  $-c_{\mathbf{\bs{\mu}}}$ is also a root of $[W_{\mathbf{\bs{\mu}},\bs{\sigma}^0}]_e$.  
\end{enumerate}

\begin{comment}
 The following lemma \ref{arg-generique} should be understood  as  a  genericity condition  that we need to impose upon $\boldsymbol{\tau}^0$ in order to prove our main result of this part.
\end{comment}

 The existence of an accumulation of zeroes  $t_{m,\bs{\mu},\bs{\sigma}^0}\in \Xi_{u,\eta}$ will be shown in Lemma \ref{accumulation-zeros-h}. A  property
that is essential to  prove   Lemma \ref{accumulation-zeros-h} is the following.
\begin{lemme}\label{arg-generique}  Let $\mathbf{\bs{\mu}} = (p, \boldsymbol{\tau}^0, \bs{\theta})$ be a parameter vector and consider $c_{\mathbf{\bs{\mu}},0}$ a root of the cyclotomic polynomial $[W_{\mathbf{\bs{\mu}},\bs{\sigma}^0}]_e$. 
 Let $\Omega_{\bs{\mu},\bs{\sigma}^0} (X) = c_{\mathbf{\bs{\mu}},0} + c_{\mathbf{\bs{\mu}},1} X^{\vartheta_{1}} + o\left(X^{\vartheta_{1}}\right)$ be a Puiseux branch of $W_{\mathbf{\bs{\mu}},\bs{\sigma}^0}(X,Y) = 0$ of main term the root  $c_{\mathbf{\bs{\mu}},0}$.  
Then there exists a generic subset $\mathcal{G}\subseteq \mathbb{R}^n$ such that $\boldsymbol{\tau}^0\in \mathcal{G}$ implies $ 
\ds \arg\left(\frac{c_{\mathbf{\bs{\mu}},1}}{c_{\mathbf{\bs{\mu}},0}}\right) \neq \frac{\pi}{2} \mod(\pi).
$
\end{lemme}


\begin{lemme}\label{accumulation-zeros-h}
 Assume that $[W_{\mathbf{\bs{\mu}},\bs{\sigma}^0}]_e$ is a cyclotomic polynomial.
There exists a branch $\Omega_{\bs{\mu},\bs{\sigma}^0}^*(X)$ of  $W_{\mathbf{\bs{\mu}},\bs{\sigma}^0}(X,Y)=0$ such that $
\ds |\Omega_{\bs{\mu},\bs{\sigma}^0}^*(X)|<1$ for $X$ sufficiently small and positive.
Thus, for some $u, \eta > 0,$ there exist infinitely many zeroes $t_{m,\bs{\mu},\bs{\sigma}^0}$ of  $\prod_{p\leq M_{\delta}}h(p^{-s^0_1-t\theta_1},\dots, p^{-s^0_n-t\theta_n})$  inside  $\Xi_{u,\eta}$. 
\end{lemme}
 We first show how Lemma \ref{accumulation-zeros-h} follows from Lemma \ref{arg-generique}, and then present the  proof of Lemma \ref{arg-generique}.

\begin{dem} ({\it of Lemma \ref{accumulation-zeros-h} assuming Lemma \ref{arg-generique}})\ 
 Consider a   branch of $W_{\mathbf{\bs{\mu}},\bs{\sigma}^0}$ with main term $c_{\mathbf{\bs{\mu}},0}$, a root of the cyclotomic polynomial  $[W_{\mathbf{\bs{\mu}},\bs{\sigma}^0}]_e$, which we   write as $\ds \Omega_{\bs{\mu},\bs{\sigma}^0}\left(X\right) = c_{\mathbf{\bs{\mu}},0} + c_{\mathbf{\bs{\mu}},1}X^{\vartheta_{1}} + o\left(X^{\vartheta_{1}}\right).$
 For  $\boldsymbol{\tau}^0 \in \mathcal{G}$  we know  that $
\ds\arg\left(\frac{c_{\mathbf{\bs{\mu}},1}}{c_{\mathbf{\bs{\mu}},0}}\right) \neq \frac{\pi}{2} \mod(\pi).
$

 Thus $
\ds \frac{\pi}{2} < \arg\left(\frac{c_{\mathbf{\bs{\mu}},1}}{c_{\mathbf{\bs{\mu}},0}}\right) < \frac{3 \pi}{2}$ or  $\ds \frac{\pi}{2} < \arg\left(-\frac{c_{\mathbf{\bs{\mu}},1}}{c_{\mathbf{\bs{\mu}},0}}\right) < \frac{3 \pi}{2}.
$
 Since $[W_{\mathbf{\bs{\mu}},\bs{\sigma}^0}]_e$ is   cyclotomic, the first  term of  $\Omega_{\bs{\mu},\bs{\sigma}^0}$ satisfies $|c_{\mathbf{\bs{\mu}},0}|=1$.
 If we first assume  that  $\frac{\pi}{2} < \arg\left(\frac{c_{\mathbf{\bs{\mu}},1}}{c_{\mathbf{\bs{\mu}},0}}\right) < \frac{3 \pi}{2}$, we have $\ds\left|\Omega_{\bs{\mu},\bs{\sigma}^0}(X)\right| = \left|1 + \frac {c_{\mathbf{\bs{\mu}},1}}{c_{\mathbf{\bs{\mu}},0}} X^{\vartheta_{1}} + o(X^{\vartheta_{1}})\right| < 1$ for $X$ small. 
 And in this event, it suffices to set  $\ds \Omega_{\bs{\mu},\bs{\sigma}^0}^* (X) = \Omega_{\bs{\mu},\bs{\sigma}^0} (X)$ to finish the proof of the Lemma.

Now assume that $
\ds\frac{\pi}{2} < \arg\left(-\frac{c_{\mathbf{\bs{\mu}},1}}{c_{\mathbf{\bs{\mu}},0}}\right) < \frac{3 \pi}{2}.
$

 Choosing  $\bs{\theta}$ to satisfy  (\ref{achoixjudicieuxtheta}), we will show  there exists a branch $\Omega_{\bs{\mu},\bs{\sigma}^0}^*$ of $W_{\bs{\mu},\bs{\sigma}^0} = 0$ whose first  term equals  $-c_{\mathbf{\bs{\mu}},0}$ and  whose   second term is the same as that for  $\Omega_{\bs{\mu},\bs{\sigma}^0}.$  Applying the above reasoning, it will follow that this suffices to finish the proof of Lemma \ref{accumulation-zeros-h}.

 The first observation uses the fact that  $-c_{\mathbf{\bs{\mu}},0}$ is also a root of  $[W_{\mathbf{\bs{\mu}},\bs{\sigma}^0}]_e$. This follows from hypothesis (a) of (\ref{achoixjudicieuxtheta}). 
And since each root  $c_{\bs{\mu},0}$ and $-c_{\bs{\mu},0}$ provides a corresponding Puiseux series solution of  $W_{\mathbf{\bs{\mu}},\bs{\sigma}^0}(X,Y)=0$ according to Proposition \ref{athpuiseuxpolgen}, there exists necessarily a branch $\Omega_{\bs{\mu},\bs{\sigma}^0}'$ of $W_{\mathbf{\bs{\mu}},\bs{\sigma}^0}(X,Y) = 0$ with Puiseux series expansion given by $\ds\Omega_{\bs{\mu},\bs{\sigma}^0}' (X) = -c_{\mathbf{\bs{\mu}},0} + c_{\mathbf{\bs{\mu}},1}' X^{\vartheta_{1}'} + o\left(X^{\vartheta_{1}'}\right)\,.$

 We now show the following equalities:
\begin{equation}\label{secondorder}
 \begin{cases}
 \displaystyle  \vartheta_{1}  =  \vartheta_{1}'  \\
 \displaystyle c_{\mathbf{\bs{\mu}},1} = c_{\mathbf{\bs{\mu}},1}'. 
 \end{cases}
 \end{equation}
 This will suffice to finish the proof in this second case since  (\ref{secondorder})  tells us that there is a branch $\Omega_{\bs{\mu},\bs{\sigma}^0}^*$ of $W_{\mathbf{\bs{\mu}},\bs{\sigma}^0} = 0$ such that $\ds\Omega_{\bs{\mu},\bs{\sigma}^0}^* = -c_{\mathbf{\bs{\mu}},0} + c_{\mathbf{\bs{\mu}},1} X^{\vartheta_{1}} + o\left(X^{\vartheta_{1}}\right)\,;$
which then permits to  argue exactly as in the first possibility.

 To prove (\ref{secondorder}), we observe  that the minimality property satisfied by $\bs{\sigma}^0\cdot \bs{\alpha}_{\cdot e'}$  implies that the terms of lowest degree in $X$ of $W_{\mathbf{\bs{\mu}},\bs{\sigma}^0}\left(X,\Omega_{\bs{\mu},\bs{\sigma}^0}\right)$ resp. $W_{\mathbf{\bs{\mu}},\bs{\sigma}^0}\left(X,\Omega_{\bs{\mu},\bs{\sigma}^0}'\right)$
coincide  with those of $W_{\mathbf{\bs{\mu}},\bs{\sigma}^0}\left(X, c_{\mathbf{\bs{\mu}},0} + c_{\mathbf{\bs{\mu}},1} X^{\vartheta_{1}}\right)$ resp. $W_{\mathbf{\bs{\mu}},\bs{\sigma}^0}\left(X, -c_{\mathbf{\bs{\mu}},0} + c_{\mathbf{\bs{\mu}},1}' X^{\vartheta_{1}'}\right)$. This leads to  the following two identities for all $\epsilon>0$:
\begin{displaymath}
 [W_{\mathbf{\bs{\mu}},\bs{\sigma}^0}]_e\left(c_{\mathbf{\bs{\mu}},0}  + c_{\mathbf{\bs{\mu}},1} X^{\vartheta_{1}}\right)  + X^{\bs{\sigma}^0\cdot \bs{\alpha}_{\cdot e'}} R_{\bs{\mu},e'}\left(c_{\mathbf{\bs{\mu}},0}  + c_{\mathbf{\bs{\mu}},1} X^{\vartheta_{1}}\right) = 0 \mod\left(X^{\bs{\sigma}^0\cdot \bs{\alpha}_{\cdot e'}+\epsilon}\right).  
\end{displaymath} 
\begin{displaymath}
 [W_{\mathbf{\bs{\mu}},\bs{\sigma}^0}]_e\left(-c_{\mathbf{\bs{\mu}},0} + c_{\mathbf{\bs{\mu}},1}' X^{\vartheta_{1}'}\right)  + X^{\bs{\sigma}^0\cdot \bs{\alpha}_{\cdot e'}} R_{\bs{\mu},e'}\left(-c_{\mathbf{\bs{\mu}},0} + c_{\mathbf{\bs{\mu}},1}' X^{\vartheta_{1}'}\right) = 0 \mod\left(X^{\bs{\sigma}^0\cdot \bs{\alpha}_{\cdot e'}+\epsilon}\right).  
\end{displaymath} 
 Consequently if $[W_{\mathbf{\bs{\mu}},\bs{\sigma}^0}]'_e$ designates the derivative of $[W_{\mathbf{\bs{\mu}},\bs{\sigma}^0}]_e$, we obtain:
\begin{equation}\label{thetwoequations}
 \begin{cases}
 \displaystyle c_{\mathbf{\bs{\mu}},1}[W_{\mathbf{\bs{\mu}},\bs{\sigma}^0}]'_e\left(c_{\mathbf{\bs{\mu}},0}\right) X^{\vartheta_{1}} +  X^{\bs{\sigma}^0\cdot \bs{\alpha}_{\cdot e'}}R_{\bs{\mu},e'}\left(c_{\mathbf{\bs{\mu}},0}\right) =0 , \\
  \displaystyle c_{\mathbf{\bs{\mu}},1}'[W_{\mathbf{\bs{\mu}},\bs{\sigma}^0}]'_e\left(-c_{\mathbf{\bs{\mu}},0}\right) X^{\vartheta_{1}'} +  X^{\bs{\sigma}^0\cdot \bs{\alpha}_{\cdot e'}}R_{\bs{\mu},e'}\left(-c_{\mathbf{\bs{\mu}},0}\right) =0.
 \end{cases}
 \end{equation}

But according to (\ref{achoixjudicieuxtheta}), we have:
$$
c_{\mathbf{\bs{\mu}},1}'[W_{\mathbf{\bs{\mu}},\bs{\sigma}^0}]'_e\left(-c_{\mathbf{\bs{\mu}},0}\right) X^{\vartheta_{1}'} +  X^{\bs{\sigma}^0\cdot \bs{\alpha}_{\cdot e'}}R_{\bs{\mu},e'}\left(-c_{\mathbf{\bs{\mu}},0}\right) = -c_{\mathbf{\bs{\mu}},1}'[W_{\mathbf{\bs{\mu}},\bs{\sigma}^0}]'_e\left(c_{\mathbf{\bs{\mu}},0}\right) X^{\vartheta_{1}'} - X^{\bs{\sigma}^0\cdot \bs{\alpha}_{\cdot e'}}R_{\bs{\mu},e'}\left(c_{\mathbf{\bs{\mu}},0}\right).
$$

 Hence from the equations (\ref{thetwoequations}), we deduce that $\vartheta_{1}'=\vartheta_{1}$ and
\begin{equation}\label{aalgopuiseux1}
 c_{\bs{\mu},1}'=c_{\bs{\mu},1} = -\frac{R_{\bs{\mu},e'}(c_{\bs{\mu},0})}{[W_{\mathbf{\bs{\mu}},\bs{\sigma}^0}]'_e(c_{\bs{\mu},0})};
\end{equation}
which proves (\ref{secondorder}).

 As a result,  when $\frac{\pi}{2} < \arg\left(-\frac{c_{\mathbf{\bs{\mu}},1}}{c_{\mathbf{\bs{\mu}},0}}\right) < \frac{3 \pi}{2}$,   there exists  a Puiseux series $\ds\Omega_{\bs{\mu},\bs{\sigma}^0}^* (X) = -c_{\mathbf{\bs{\mu}},0} + c_{\mathbf{\bs{\mu}},1}X^{\vartheta_{1}} + o\left(X^{\vartheta_{1}}\right)$
such that  $|\Omega_{\bs{\mu},\bs{\sigma}^0}^*(X)|<1$ for $X$ positive in a neighborhood of $0$. 
Thus, in either of the two cases for $\arg\left(-\frac{c_{\mathbf{\bs{\mu}},1}}{c_{\mathbf{\bs{\mu}},0}}\right)$ we obtain zeroes of $W_{\mathbf{\bs{\mu}},\bs{\sigma}^0}(p^{-1}, p^{-t})$ whenever $
\ds t_{m,\bs{\mu},\bs{\sigma}^0} = -\frac{\log\left(\Omega_{\bs{\mu},\bs{\sigma}^0}^* \left(p^{-1}\right)\right)}{\log(p)} + \frac{2 \pi m i}{\log(p)},
$
where $m\in \mathbb{Z}$ and $p$ is a prime number large enough.
 It follows that  $t_{m,\bs{\mu},\bs{\sigma}^0} \in \Xi_{u,\eta}$ if $u<\Im(t_{m,\bs{\mu},\bs{\sigma}^0})<u+\eta$, that is, if $
\ds u<\frac{2\pi m}{\log\left(p\right)} - \frac{\arg\left(\Omega_{\bs{\mu},\bs{\sigma}^0}^* \left(p^{-1}\right)\right)}{\log\left(p\right)} < u+\eta,
$
which is equivalent to:
\begin{equation}\label{boite}
\frac{u\log\left(p\right)}{2\pi} + \frac{\arg\left(\Omega_{\bs{\mu},\bs{\sigma}^0}^* \left(p^{-1}\right)\right)}{2\pi}<m<\frac{\left(u+\eta\right)\log\left(p\right)}{2\pi} + \frac{\arg\left(\Omega_{\bs{\mu},\bs{\sigma}^0}^* \left(p^{-1}\right)\right)}{2\pi}.
\end{equation}
For $p$ large enough, we will therefore have   zeroes of $t \to W_{\mathbf{\bs{\mu}},\bs{\sigma}^0} (p^{-1}, p^{-t})$ inside $\Xi_{u,\eta}$. Allowing $p \to \infty$ will then produce an infinite set of zeroes within   $\Xi_{u,\eta}$ of positive real part tending to $0$ as $p$ tends to infinity;
which completes  the proof that Lemma \ref{arg-generique}
implies Lemma \ref{accumulation-zeros-h}. \qed
\end{dem}
Let us prove now Lemma \ref{arg-generique}:
\begin{dem}[\textbf{Lemma \ref{arg-generique}}]
 The idea of the proof is to exhibit clearly how $\arg \left(c_{\mathbf{\bs{\mu}},1}/c_{\mathbf{\bs{\mu}},0}\right)$ depends upon $\boldsymbol{\tau}^0.$ To do this, we need to revisit the construction of $c_{\mathbf{\bs{\mu}},1}$. 
According to (\ref{aalgopuiseux1}) we have $
c_{\bs{\mu},1} = -\frac{R_{\bs{\mu},e'}(c_{\bs{\mu},0})}{[W_{\mathbf{\bs{\mu}},\bs{\sigma}^0}]'_e(c_{\bs{\mu},0})}.
$
 Recall also that   (\ref{adepenpcko}) implies that $
\tilde{c}_{\mathbf{\bs{\theta}},0} := c_{\mathbf{\bs{\mu}},0} p^{-i\frac{\bs{\tau}^0 \cdot \bs{\alpha}_{\cdot e}}{\bs{\theta}\cdot \bs{\alpha}_{\cdot e}}} 
$
is {\it independent of} $\boldsymbol{\tau}^0$ and $p$.
As a result, we first observe:

 {\it The denominator $c_{\mathbf{\bs{\mu}},0} [W_{\mathbf{\bs{\mu}},\bs{\sigma}^0}]'_e\left(c_{\mathbf{\bs{\mu}},0}\right)$ of $\frac{c_{\mathbf{\bs{\mu}},1}}{c_{\mathbf{\bs{\mu}},0}}$ is   independent of  $\boldsymbol{\tau}^0$ and $p$.}

 To justify this assertion we use the fact that if $j \in \Lambda_e$ then $\frac{\bs{\tau}^0 \cdot \bs{\alpha}_{\cdot e}}{\bs{\theta}\cdot \bs{\alpha}_{\cdot e}}\bs{\theta} \cdot \bs{\alpha}_{\cdot j} = \bs{\tau}^0 \cdot \bs{\alpha}_{\cdot j}.$
Consequently we have:
\begin{displaymath}
\begin{array}{lll}
\displaystyle c_{\bs{\mu},0} [W_{\mathbf{\bs{\mu}},\bs{\sigma}^0}]'_e\left(c_{\bs{\mu},0}\right) & = & \displaystyle \tilde{c}_{\mathbf{\bs{\theta}},0} p^{i\frac{\bs{\tau}^0 \cdot \bs{\alpha}_{\cdot e}}{\bs{\theta}\cdot \bs{\alpha}_{\cdot e}}}\sum_{j\in \Lambda_e} a_j (\bs{\theta} \cdot \bs{\alpha}_{\cdot j}) \tilde{c}_{\mathbf{\bs{\theta}},0}^{\bs{\theta} \cdot \bs{\alpha}_{\cdot j}-1}p^{i\left(\frac{\bs{\tau}^0 \cdot \bs{\alpha}_{\cdot e}}{\bs{\theta}\cdot \bs{\alpha}_{\cdot e}}\left(\bs{\theta} \cdot \bs{\alpha}_{\cdot j} -1\right)-\bs{\tau}^0 \cdot \bs{\alpha}_{\cdot j}\right)} \\ 
 & = & \displaystyle \sum_{j\in \Lambda_e}a_j (\bs{\theta} \cdot \bs{\alpha}_{\cdot j}) \tilde{c}_{\mathbf{\bs{\theta}},0}^{\bs{\theta} \cdot \bs{\alpha}_{\cdot j}},
\end{array}
\end{displaymath} 
 We next deal with the numerator for $c_{\mathbf{\bs{\mu}},1}.$ For each  $k$ such that $k\sim e'$, we define $
\ds \lambda_k:= \frac{- \ a_k \ \tilde c_{\mathbf{\bs{\mu},0}}^{\bs{\theta}\cdot \bs{\alpha}_{\cdot k}}}{\sum_{j\in \Lambda_e}a_j (\bs{\theta} \cdot \bs{\alpha}_{\cdot j}) \tilde{c}_{\mathbf{\bs{\theta}},0}^{\bs{\theta} \cdot \bs{\alpha}_{\cdot j}}}\,. 
$
Notice that each $\lambda_k$ is independent of  $\boldsymbol{\tau}^0$ and  $p$. Moreover, a straightforward calculation shows: 
\begin{displaymath}
\begin{array}{lll}
\displaystyle \frac{c_{\mathbf{\bs{\mu}},1}}{c_{\mathbf{\bs{\mu}},0}} & = & \displaystyle \sum_{\{k :  \bs{\alpha}_{\cdot k}-\boldsymbol{\alpha}_{\cdot e'}\in \langle \boldsymbol{\alpha}_{\cdot e} \rangle\}} \lambda_k \, p^{i\left(\bs{\tau}^0 \cdot \bs{\alpha}_{\cdot e}\frac{\bs{\theta}\cdot \bs{\alpha}_{\cdot k}}{\bs{\theta}\cdot \bs{\alpha}_{\cdot e}}-\bs{\tau}^0\cdot \bs{\alpha}_{\cdot k}\right)}  
  = \displaystyle \sum_{\{k :  \bs{\alpha}_{\cdot k}-\boldsymbol{\alpha}_{\cdot e'}\in \langle \boldsymbol{\alpha}_{\cdot e} \rangle\}} \lambda_k \, p^{i\bs{\tau}^0\cdot \mathbf{w}_k},
\end{array}
\end{displaymath} 
where $ 
\mathbf{w}_k:=\frac{\bs{\theta}\cdot \bs{\alpha}_{\cdot k}}{\bs{\theta}\cdot \bs{\alpha}_{\cdot e}}\boldsymbol{\alpha}_{\cdot e} - \bs{\alpha}_{\cdot k}.
$
 We next observe that $\mathbf{w}_k$ is independent of $k.$ Indeed,  if $k$ and $k'$ are such that $\bs{\alpha}_{\cdot k}-\boldsymbol{\alpha}_{\cdot e'}\in \langle \boldsymbol{\alpha}_{\cdot e} \rangle$ and $\bs{\alpha}_{\cdot k'}-\boldsymbol{\alpha}_{\cdot e'}\in \langle \boldsymbol{\alpha}_{\cdot e} \rangle$, then $\bs{\alpha}_{\cdot k}-\bs{\alpha}_{\cdot k'}\in \langle \boldsymbol{\alpha}_{\cdot e} \rangle$. Thus,  there exists $q\in \mathbb{Q}$ such that $ 
\bs{\alpha}_{\cdot k}-\bs{\alpha}_{\cdot k'} = q  \boldsymbol{\alpha}_{\cdot e}.
$
But then $
\ds\frac{\bs{\theta}\cdot (\bs{\alpha}_{\cdot k}-\bs{\alpha}_{\cdot k'})}{\bs{\theta}\cdot \bs{\alpha}_{\cdot e}}\boldsymbol{\alpha}_{\cdot e} = q  \boldsymbol{\alpha}_{\cdot e} = \bs{\alpha}_{\cdot k}-\bs{\alpha}_{\cdot k'};
$
implies $\mathbf{w}_k=\mathbf{w}_{k'}$. Thus, in particular $\mathbf{w}_k = \mathbf{w}_{e'}$. 
We conclude that $
\ds\frac{c_{\mathbf{\bs{\mu}},1}}{c_{\mathbf{\bs{\mu}},0}} = p^{i\bs{\tau}^0\cdot \mathbf{w}_{e'}}\sum_{\{k : \bs{\alpha}_{\cdot k}-\boldsymbol{\alpha}_{\cdot e'}\in \langle \boldsymbol{\alpha}_{\cdot e} \rangle \}}\lambda_k.
$

 If we now define $
\ds\varphi_{e'} := \arg\left(\sum_{\{k : \bs{\alpha}_{\cdot k}-\boldsymbol{\alpha}_{\cdot e'}\in \langle \boldsymbol{\alpha}_{\cdot e} \rangle \}} \lambda_k\right),
$ 
then $\varphi_{e'}$ is independent of $\boldsymbol{\tau}^0$ and $p$. The condition that $\arg (c_{\mathbf{\bs{\mu}},1}/c_{\mathbf{\bs{\mu}},0}) = \pi/2$ then reduces to exactly one of the following two conditions:   
\begin{equation*}
   \left(\bs{\tau}^0\cdot \mathbf{w}_{e'}\right)\log(p) + \varphi_{e'} \mod(2 \pi) \in \left\{ \frac{\pi}{2};\frac{3\pi}{2} \right\}.
 \end{equation*}
But the set $
\ds M:=\bigcup_{p}\bigcup_{u\in \mathbb{Z}}\left\{\boldsymbol{\tau}^0\in \mathbb{R}^n : \left(\bs{\tau}^0\cdot \mathbf{w}_{e'}\right)\log(p) + \varphi_{e'} =\frac{\pi}{2} + u\pi \right\} 
$
is of empty interior according to Baire's Theorem since it is a countable union of sets of empty interior.
 Now it suffices to define $\mathcal{G}$ as the {\it complement} of $M$ in order to obtain a generic set of $\boldsymbol{\tau}^0$ for which $\arg (c_{\mathbf{\bs{\mu}},1}/c_{\mathbf{\bs{\mu}},0}) \neq \pi/2$. This completes the proof of Lemma \ref {arg-generique}. 
\CQFD
\end{dem}

\noindent  {\bf Proof of Part 2.}

In this part we still assume that $h$ is not cyclotomic and has no cyclotomic factor, and here $[h]_e$ is not necessarily a cyclotomic polynomial.

   We have, so far,   found an infinite number of zeroes $t_{m,\bs{\mu},\bs{\sigma}^0}$ of $\ds t \longmapsto h\left(p^{-s^0_1-t \bs{\theta}_1}, \cdots,p^{-s^0_n-t \bs{\theta}_n}\right) = W_{\mathbf{\bs{\mu}},\bs{\sigma}^0} (p^{-1}, p^{-t})$ which accumulate inside $\Xi_{u,\eta}$. This is possible provided  $\boldsymbol{\sigma}^0$ satisfies (\ref{agenericite}), $\boldsymbol{\tau}^0$ belongs to the generic set $\mathcal{G}$ (see Lemma \ref {arg-generique}), and $\bs{\theta}$ satisfies  (\ref{ahyp_theta1}) and (\ref{achoixjudicieuxtheta}).

To prove Theorem \ref{aresultatprincipal}, it is necessary to ensure that these zeroes are {\it not} cancelled by possible zeroes/poles  coming from the factor   $t \to A_{M_{\delta}}(\mathbf{s}^0 + t \bs{\theta})$ of $Z(\mathbf{s}^0 + t \bs{\theta})$,  determined by    Theorem \ref{aana}. Such   zeroes/poles  must be  of the form $
\ds t(\bs{\beta},\rho) = \frac{\rho - \mathbf{s}^0\cdot\bs{\alpha}\cdot {}^t \!\bs{\beta}}{\bs{\theta}\cdot\bs{\alpha}\cdot {}^t \!\bs{\beta}} \ \big(\Re (t(\bs{\beta}, \rho)) > 0 \big),
$
where $\boldsymbol{\beta}\in \Nr$ and $\rho$ is the pole $1$ or a non-trivial zero of the Riemann zeta function. 

To do this, it will be convenient to work locally in a ``sufficiently small" neighborhood $\mathcal {B}$ of  {\it any}  generic  point $\hat {\mathbf s}$ whose real part  satisfies (\ref{agenericite}). By ``sufficiently small"  we mean the following: 
\begin{equation} \label{genericite+}
\text{{\it there exists   $\kappa = \kappa(\mathcal{B}) > 0$ such that  
$\mathbf s = \boldsymbol{\sigma}^0 + i \boldsymbol{\tau}^0 \in \mathcal {B}$ \ \text{{\it implies}} $\bs{\sigma}^0\cdot \bs{\alpha}_{\cdot j} \ge \kappa \ \ \forall j \notin \Lambda_e.$}}
\end{equation}
In addition, we must also  distinguish between the two possibilities:   
$$\textrm {Case A:} \quad  \bs{\beta}\notin B_e;   \qquad \textrm{Case B:} \quad \boldsymbol{\beta}\in B_e\,.$$
The reason for this is that we need to prove that cancellation with the zeroes $t(\bs{\beta}, \rho)$ does not occur uniformly in $\bs{\beta}$, and the arguments that do this depend upon  membership in $B_e.$ 

 Case A is dealt with in   Lemma \ref{perturbation-sigma0} where we show that for  $p$ large enough (i.e. $p>p_0$ where $p_0$ is a  constant depending {\it only} upon the chosen neighborhood of $\hat {\mathbf s}$ and is, in particular, independent of $\bs{\beta}\notin B_e$) the quantity $\ds h\left(p^{-s^0_1-t(\bs{\beta},\rho)\theta_1},\cdots,p^{-s^0_n-t(\bs{\beta},\rho)\theta_n}\right) \neq 0$
for all $\bs{\beta}\notin  B_e$ and $\rho.$ This, however, will   be possible only after finding an additional  genericity condition that the components $\boldsymbol{\sigma}^0$ and $\boldsymbol{\tau}^0$  should satisfy. Once this conclusion is established, we know, at least for the $\bs{\beta}\notin B_e,$ that 
there must be infinitely many zeroes/poles of $W_{\mathbf{\bs{\mu}},\bs{\sigma}^0}(p^{-1}, p^{-t})$ in $\Xi_{u,\eta}$ that are {\it not} cancelled by the poles/zeroes  $t(\bs{\beta}, \rho)$.

 Case B will be dealt with in  Lemma  \ref{lemmeBe}. The principal difficulty we encounter here will be that  $h\left(p^{-s^0_1-t\theta_1},\cdots,p^{-s^0_n-t\theta_n}\right)$ {\it can vanish} for some $t = t(\bs{\beta}, \rho).$ This is possible, in particular,  when $W_{\mathbf{\bs{\mu}},\bs{\sigma}^0}$ admits a root in $Y$ that is independent of $X$ and of modulus strictly less than $1$. We overcome this problem by    establishing   a weaker property that still suffices to prove Theorem \ref{aresultatprincipal}. Indeed, what saves us is the fact that the set of such $t(\bs{\beta}, \rho)$ is distinctly smaller than the set of all   zeroes $t_{m,\bs{\mu},\bs{\sigma}^0}$ constructed in Part 1.

 The reader should therefore also  be alert to   the change in perspective between Parts 1 and  2. In Part 2, the genericity conditions we   impose will be {\it locally} applied in a suitable   neighborhood (i.e. satisfying (\ref{genericite+})) of a given point $\hat {\mathbf s}$. In Part 1, the genericity conditions were essentially global in nature.

\begin{rqs} \label{parameternu}  It will be useful in this part to consider the components of the parameter vector $\bs{\mu}=(p, \bs{\tau}^0, \bs{\theta})$ introduced previously. This is because we will need to construct an explicit function of $\boldsymbol{\sigma}^0$ in order to discover the additional genericity  condition. Thus, in Part 2, it will be understood that we consider the parameter vector:
$$\bs{\mu} = (p, \boldsymbol{\tau}^0, \bs{\theta}),$$ 
where $p$ is a prime number, $\boldsymbol{\tau}^0 \in \mathcal{G}$ and $ \theta$ satisfies (\ref{ahyp_theta1}) and (\ref{achoixjudicieuxtheta}).
\end{rqs}
 Before stating the two main lemmas of Part 2, we   first derive a concise (and useful) expression for $ h\left(p^{-s^0_1-t\theta_1},\cdots,p^{-s^0_n-t\theta_n}\right)$ that distinguishes   clearly  between   dependance upon   $\bs{\mu}, \bs{\beta}, \rho$ (as parameters) and $\boldsymbol{\sigma}^0$ as a domain variable for a function. 

Let $\boldsymbol{\beta}\in \Nr$,     $\rho$ a zero of  $\zeta(s)=0$ with $\Re (\rho) \in ]0, 1[$ (or  the pole at $s = 1$), and $p$ a prime number.  Given $\bs{\mu} = (p,  \boldsymbol{\tau}^0, \bs{\theta})$ a parameter vector as above, we begin by writing:
\begin{displaymath}
\begin{array}{lll}
f_{\bs{\mu}, \bs{\beta}, \rho}(\boldsymbol{\sigma}^0) := h\left(p^{-s^0_1-t(\bs{\beta},\rho)\theta_1},\cdots,p^{-s^0_n-t(\bs{\beta},\rho)\theta_n}\right) & = & 1+\summ_{k=1}^{r}a_k p^{-\mathbf{s}^0\cdot \bs{\alpha}_{\cdot k}-\bs{\theta}\cdot \bs{\alpha}_{\cdot k}\left( \frac{\rho - \mathbf{s}^0\cdot\bs{\alpha}\cdot {}^t \!\bs{\beta}}{\bs{\theta}\cdot\bs{\alpha}\cdot {}^t \!\bs{\beta}}\right)}  \\ 
 & = & 1+\summ_{k=1}^{r}a_k p^{-u_{k,\bs{\theta}} (\boldsymbol{\sigma}^0, \bs{\beta}) - v_{k, \bs{\theta}, \bs{\tau}^0,\rho} (\bs{\beta})}
\end{array}
\end{displaymath}
where
$$
u_{k,\bs{\theta}} (\boldsymbol{\sigma}^0, \bs{\beta}) =  \bs{\sigma}^0\cdot \bs{\alpha}_{\cdot k} - \bs{\sigma}^0\cdot\bs{\alpha}\cdot {}^t \!\bs{\beta} \frac{\bs{\theta}\cdot \bs{\alpha}_{\cdot k} }{\bs{\theta}\cdot\bs{\alpha}\cdot {}^t \!\bs{\beta}}; \ 
v_{k, \bs{\theta}, \bs{\tau}^0,\rho} (\bs{\beta}) = i \left\{\bs{\tau}^0\cdot \bs{\alpha}_{\cdot k} - \bs{\tau}^0\cdot\bs{\alpha}\cdot {}^t \!\bs{\beta} \frac{\bs{\theta}\cdot \bs{\alpha}_{\cdot k}}{\bs{\theta}\cdot\bs{\alpha}\cdot {}^t \!\bs{\beta}}\right\} + \rho \frac{\bs{\theta}\cdot \bs{\alpha}_{\cdot k}}{\bs{\theta}\cdot\bs{\alpha}\cdot {}^t \!\bs{\beta}}.
$$

 We  note that $v_{k, \bs{\theta}, \bs{\tau}^0,\rho}(\bs{\beta})$ is independent of $\boldsymbol{\sigma}^0$ and $u_{k,\bs{\theta}}(\bs{\sigma}^0,\bs{\beta})$ is independent of $\boldsymbol{\tau}^0$ and linear in $\boldsymbol{\sigma}^0.$  
 We first simplify the right side by grouping together terms indexed by $k$ that give the same $u_{k,\bs{\theta}}$ function. Thus, we define the equivalence relation $\mathcal{R}_\beta$ on $\{1,\dots, r\}$:
$$k \ \mathcal{R}_\beta \ k'   \Longleftrightarrow  \displaystyle \textrm{for all} \ \boldsymbol{\sigma}^0 \quad u_{k,\bs{\theta}}(\boldsymbol{\sigma}^0,\bs{\beta}) = u_{k',\bs{\theta}}(\boldsymbol{\sigma}^0,\bs{\beta}).$$

 We denote the equivalence class of $k$ by $[k]$ and let   $\mathcal{V}$ denote a set of representatives of these classes.

\begin{rqs}\label{rqBe}
  When $\boldsymbol{\beta}\in B_e$ and $\boldsymbol{\sigma}^0$ satisfies (\ref{agenericite}), it follows that \ $u_{k,\bs{\theta}}(\boldsymbol{\sigma}^0, \bs{\beta})  = \bs{\sigma}^0\cdot \bs{\alpha}_{\cdot k}.$ We conclude that  $\boldsymbol{\beta}\in B_e$ implies $\ds\bs{\alpha}_{\cdot k} \ \mathcal{R}_\beta \ \bs{\alpha}_{\cdot k'}$ if and only if $\bs{\alpha}_{\cdot k}-\bs{\alpha}_{\cdot k'}\in \langle \boldsymbol{\alpha}_{\cdot e} \rangle.$  
\end{rqs}
 Thus we have the following expression for $f_{\bs{\mu}, \bs{\beta},\rho}\left( \boldsymbol{\sigma}^0 \right)$:  
 \begin{equation}\label{fonctionidentnulle}
 f_{\bs{\mu}, \bs{\beta},\rho}\left( \boldsymbol{\sigma}^0 \right) = 1+\summ_{\nu \in \mathcal{V}}\left(\summ_{k \in [\nu]}a_k p^{-v_{k, \bs{\theta}, \bs{\tau}^0,\rho} (\bs{\beta})} \right) p^{-u_{\nu,\bs{\theta}}(\boldsymbol{\sigma}^0, \bs{\beta})},
 \end{equation}
 where now, the  $-u_{\nu,\bs{\theta}}(\boldsymbol{\sigma}^0, \bs{\beta})$ are 
 pairwise distinct.

\begin{lemme}\label{perturbation-sigma0}
Let $\hat {\mathbf s}$ be any point such that $\Re (\hat {\mathbf s})$ is a generic point of $\partial \mathbf{W}(0)$ (i.e. satisfying (\ref{agenericite})). Let $\mathcal B$ denote a   neighborhood of $\hat {\mathbf s}$ satisfying (\ref{genericite+}). 

There is a generic set $\mathcal {G}^* \subset \{\mathbf s : \bs{\sigma}\cdot\bs{\alpha}_{\cdot e} = 0\}=\mathcal{F}(\bs{\alpha}_{\cdot e})$
 satisfying the following property:
{\it There exists $p_0 = p_0 (\mathcal {B})$ such that $p > p_0$ implies $\ds h\left(p^{-s^0_1-t(\bs{\beta},\rho)\theta_1},\dots,p^{-s^0_n-t(\bs{\beta},\rho)\theta_n}\right) \neq 0 $
if $\mathbf s^0 \in \mathcal B \cap \mathcal {G}^*$, $\bs{\beta}\notin B_e,$ and  $\rho$ satisfies   $\Re(t(\bs{\beta},\rho))>0$.}
\end{lemme}

\begin{dem}  Since $h\left(p^{-s^0_1-t(\bs{\beta},\rho)\theta_1},\cdots,p^{-s^0_n-t(\bs{\beta},\rho)\theta_n}\right) = f_{\bs{\mu},\bs{\beta},\rho}(\boldsymbol{\sigma}^0)$, we will   use (\ref{fonctionidentnulle}) to  verify the assertion.  The main result that we need to do this is as follows:

{\bf Claim:}  {\it There exist $\nu \in \mathcal{V}$, a generic set $\mathcal{G}' \subset \mathcal{G}$ (see Lemma \ref{arg-generique}),  and $p_0 = p_0(\mathcal{B})$  such that for all $\bs{\beta}\notin B_e$, and $\rho$ satisfying  $\Re \left(t(\bs{\beta}, \rho)\right) > 0$ we have $ \sum_{k \in [\nu]} a_k p^{-v_{k, \bs{\theta}, \bs{\tau}^0,\rho}(\bs{\beta})} \neq 0$
whenever  $\bs{\mu} = (p, \boldsymbol{\tau}^0, \bs{\theta})$ is a parameter vector  for which $p>p_0$ and  $\boldsymbol{\tau}^0 \in \mathcal{G}'$.}

{\bf Proof of Claim.}   For each $k\in\{1,\dots,r\}$ we define the following expressions $$\mathcal{A}_{k,\bs{\theta}}(\bs{\beta}) := \bs{\alpha}_{\cdot k}-\frac{\bs{\theta}\cdot \bs{\alpha}_{\cdot k}}{\bs{\theta}\cdot\bs{\alpha}\cdot {}^t \!\bs{\beta}}\boldsymbol{\alpha} \cdot {}^t \! \boldsymbol{\beta}; \qquad   \mathcal{A}_{k,\bs{\theta}}'(\bs{\beta}, \rho) := \rho \cdot \frac{\bs{\theta}\cdot \bs{\alpha}_{\cdot k}}{\bs{\theta}\cdot\bs{\alpha}\cdot {}^t \!\bs{\beta}}. 
$$
It follows that $\ds\sum_{k\in [\nu]} a_k p^{-v_{k, \bs{\theta}, \bs{\tau}^0,\rho}(\bs{\beta})} =\sum_{k\in [\nu]} a_k p^{-\mathcal{A}_{k,\bs{\theta}}' (\bs{\beta},\rho) -i \ \bs{\tau}^0\cdot\mathcal{A}_{k,\bs{\theta}}(\bs{\beta})}.$

For each $\nu$ and each $\bs{\beta}\notin B_e,$ we next define the relation $\mathfrak{R}_\beta$ on $[\nu]$ by setting
$$k_1 \ \mathfrak{R}_\beta \ k_2 \quad \text{{\it iff }} \quad \mathcal{A}_{k_1, \theta}(\bs{\beta}) = \mathcal{A}_{k_2, \theta}(\bs{\beta})\,,$$
and  denote  a fixed set of representatives by $\mathcal{M}_{\nu,\bs{\beta}}$. Thus,   two distinct elements of $\mathcal{M}_{\nu,\bs{\beta}}$ determine distinct 
values of $\mathcal{A}_{\cdot, \bs{\theta}}(\bs{\beta}).$   
  
Then we easily check that if we put $
\ds w_{j,k,\bs{\theta},\rho}(\bs{\beta}) = \rho \cdot \frac{\bs{\theta}\cdot \left(\bs{\alpha}_{\cdot j}-\bs{\alpha}_{\cdot k}\right)}{\bs{\theta}\cdot\bs{\alpha}\cdot {}^t \!\bs{\beta}}\,,
$ we have:
\begin{equation}\label{exprgene}
 \sum_{k\in [\nu]} a_k p^{-v_{k,\bs{\mu}, \rho}(\bs{\beta})} = \sum_{j\in \mathcal{M}_{\nu,\bs{\beta}}} \left(\sum_{k \in [\nu]; \  \mathcal{A}_{k,\bs{\theta}}(\bs{\beta}) = \mathcal{A}_{j,\bs{\theta}}(\bs{\beta})} a_{k} p^{w_{j,k,\bs{\theta},\rho}(\bs{\beta})}\right) p^{-\mathcal{A}_{j,\bs{\theta}}'(\bs{\beta}, \rho)} p^{-i \left(\boldsymbol{\tau}^0 \cdot \mathcal{A}_{j,\bs{\theta}}(\bs{\beta})\right)}.
\end{equation}
 We now proceed to a proof of the Claim by contradiction. Thus, we assume the assertion is false. This means that  for any generic subset $\widetilde {\mathcal {G}}$, each $\nu \in \mathcal {V}$ and each prime $p_0,$ there exist $p > p_0$ and some $\bs{\beta}\notin B_e$ so that if $\boldsymbol{\tau}^0 \in \widetilde {\mathcal{G}}$ and $\bs{\mu} = (p, \boldsymbol{\tau}^0, \bs{\theta}),$ then $\sum_{k \in [\nu]} a_k p^{-v_{k, \bs{\theta}, \bs{\tau}^0,\rho}(\bs{\beta})} = 0$
for some $\rho $ such that  $\Re \left(t(\bs{\beta},\rho)\right) > 0.$ 

Since the sum over $k$ on the left side of (\ref{exprgene}) is assumed to equal $0$, exactly one of two possibilities can occur for the inner sum on the right side of (\ref{exprgene}). Either: 
\begin{enumerate}
\item $\sum_{k \in [\nu]; \  \mathcal{A}_{k,\bs{\theta}}(\bs{\beta}) = \mathcal{A}_{j,\bs{\theta}}(\bs{\beta})} a_{k} p^{w_{j,k,\bs{\theta},\rho}(\bs{\beta})} = 0$ \  for each $j \in \mathcal{M}_{\nu,\bs{\beta}}\,;$

or 

\item  there exists $j \in \mathcal{M}_{\nu,\bs{\beta}}$ such that $\sum_{k \in [\nu]; \  \mathcal{A}_{k,\bs{\theta}}(\bs{\beta}) = \mathcal{A}_{j,\bs{\theta}}(\bs{\beta})} a_{k} p^{w_{j,k,\bs{\theta},\rho}(\bs{\beta})} \neq 0\,.$  

\end{enumerate}

 {\it We first show that possibility (1) cannot occur. }

 If possibility (1) did occur, then we first notice that the condition  $\mathcal{A}_{j, \bs{\theta}}(\bs{\beta}) = \mathcal{A}_{k,\bs{\theta}}(\bs{\beta})$ implies $
\boldsymbol{\alpha}_{\cdot j} - \bs{\alpha}_{\cdot k} = \frac{\bs{\theta}\cdot \left(\bs{\alpha}_{\cdot j}-\bs{\alpha}_{\cdot k}\right)}{\bs{\theta}\cdot\bs{\alpha}\cdot {}^t \!\bs{\beta}} \cdot \boldsymbol{\alpha} \cdot {}^t \! \boldsymbol{\beta}.
$
We now take the inner product with any $\boldsymbol{\sigma}^0 = \Re (\mathbf s^0)$ such that $\mathbf s^0 \in \mathcal {B}.$ We obtain the equality $
\bs{\sigma}^0\cdot \left(\bs{\alpha}_{\cdot j}-\bs{\alpha}_{\cdot k}\right) = \frac{\bs{\theta}\cdot \left(\bs{\alpha}_{\cdot j}-\bs{\alpha}_{\cdot k}\right)}{\bs{\theta}\cdot\bs{\alpha}\cdot {}^t \!\bs{\beta}} \cdot  \sum_{i\notin \Lambda_e}\beta_i(\bs{\sigma}^0\cdot \bs{\alpha}_{\cdot i}).
$
Since  $\boldsymbol{\beta}\in \mathbb{N}^r - B_e$    implies 
 $\sum_{i\notin \Lambda_e}\beta_i(\bs{\sigma}^0\cdot \bs{\alpha}_{\cdot i})\neq 0$, we conclude:
\begin{equation}\label{refw}
w_{j,k,\bs{\theta},\rho}(\bs{\beta}) = \rho \cdot \frac{\bs{\sigma}^0\cdot \left(\bs{\alpha}_{\cdot j}-\bs{\alpha}_{\cdot k}\right)} {\sum_{i\notin \Lambda_e}\beta_i(\bs{\sigma}^0\cdot \bs{\alpha}_{\cdot i})}.
\end{equation}
 It is now important to observe that  the set 
\begin{equation}\label{equationfini}
\mathcal{E}:=\{\beta_i :  i\notin \Lambda_e, \gamma(\boldsymbol{\beta})\neq 0, \Re(t(\bs{\beta},\rho))\geq 0 \} \ \textrm{{\it is finite.}}
\end{equation}
Indeed, since the $t(\bs{\beta},\rho)$, which could cancel the zeroes $t_{m,\bs{\mu},\bs{\sigma}^0}$ from Part 1, are necessarily of positive real part, we have
\begin{equation}\label{preellepos}
\Re(t(\bs{\beta},\rho)) = \frac{\Re(\rho)-\sum_{i\notin \Lambda_e}\beta_i(\bs{\sigma}^0\cdot \bs{\alpha}_{\cdot i})}{\bs{\theta}\cdot\bs{\alpha}\cdot {}^t \!\bs{\beta}} \geq 0.
\end{equation}
Consequently (using the bound from (\ref{genericite+})), $
\kappa \sum_{i \notin \Lambda_e} \beta_i \le \sum_{i\notin \Lambda_e}\beta_i(\bs{\sigma}^0\cdot \bs{\alpha}_{\cdot i}) \leq \Re(\rho) <1.
$
This shows  that $\mathcal{E}$ is a finite set (and  depends upon $\mathcal {B}$).

Now if $p$ is a prime number satisfying the equation in possibility (1), then necessarily $p^{\frac{\rho}{\sum_{i\notin \Lambda_e}\beta_i(\bs{\sigma}^0\cdot \bs{\alpha}_{\cdot i})}}$ is a solution of the equations (indexed by $j \in \mathcal{M}_{\nu,\bs{\beta}}$) 
\begin{equation}\label{eqenX}
\sum_{k \in [\nu], \ \mathcal{A}_{k,\bs{\theta}}(\bs{\beta}) = \mathcal{A}_{j,\bs{\theta}}(\bs{\beta})} a_{k} X^{\bs{\sigma}^0\cdot \left(\bs{\alpha}_{\cdot j}-\bs{\alpha}_{\cdot k}\right)} = 0.
\end{equation}

 These equations do, in fact, depend upon $X$, that is, they cannot
reduce to a linear relation among the coefficients $\{a_k\}_{k\in [\nu]},$ which could only occur if the $\bs{\sigma}^0\cdot \bs{\alpha}_{\cdot k},$  $(k \in [\nu])$ were not pairwise distinct. However, this cannot happen.  

Indeed, if $\bs{\sigma}^0\cdot  \bs{\alpha}_{\cdot k_1} = \bs{\sigma}^0\cdot \bs{\alpha}_{\cdot k_2}$ then $\mathcal{A}_{k_1, \bs{\theta}}(\bs{\beta}) = \mathcal{A}_{k_2, \bs{\theta}} (\bs{\beta})$ implies $\mathcal{A}_{k_1,\bs{\theta}}'(\bs{\beta}, \rho) = \mathcal{A}_{k_2,\bs{\theta}}'(\bs{\beta}, \rho)$. Thus, 
$$\rho \ \bs{\alpha}_{\cdot k_1} = \rho \mathcal{A}_{k_1, \bs{\theta}}(\bs{\beta}) + \mathcal{A}_{k_1,\bs{\theta}}'(\bs{\beta}, \rho) \bs{\alpha}\cdot {}^t \!\bs{\beta} = \rho \mathcal{A}_{k_2, \bs{\theta}}(\bs{\beta}) + \mathcal{A}_{k_2,\bs{\theta}}'(\bs{\beta}, \rho) \bs{\alpha}\cdot {}^t \!\bs{\beta} = \rho \ \bs{\alpha}_{\cdot k_2},$$
and  $\bs{\alpha}_{\cdot k_1}=\bs{\alpha}_{\cdot k_2}$ follows.

It is then clear that there exists $\mathfrak{M} = \mathfrak{M}(\mathcal {B}) > 0$ such that all common solutions to these equations belong to the interval  $|X| \le  \mathfrak{M}.$  

Furthermore, we have a positive lower bound  for $\Re(\rho)$  since $\Re(t(\bs{\beta},\rho))\geq 0$ implies, by (\ref{preellepos}), that $
\Re(\rho) > \min_{j\notin\Lambda_e}\bs{\sigma}^0 \cdot  \bs{\alpha}_{\cdot j}> \kappa.
$

It is then clear that  there exists a  constant $p_0 = p_0(\mathcal{B})$ such that  $p>p_0$ implies 
$$
\left|p^{\frac{\rho}{\sum_{j\notin \Lambda_e} \beta_j(\bs{\theta}\cdot \bs{\alpha}_{\cdot j})}}\right|  \geq \exp\left(\log(p)\frac{\min_{j\notin\Lambda_e}\bs{\sigma}^0 \cdot  \bs{\alpha}_{\cdot j}}{\max_{\beta_j \in \mathcal {E}}\beta_j\sum_{j\notin \Lambda_e}\bs{\sigma}^0 \cdot  \bs{\alpha}_{\cdot j}}\right) >\mathfrak{M}.
$$
Thus,  possibility  (1) {\it cannot} occur.

{\it We now show possibility (2) cannot occur.}

 Given that $\sum_{k \in [\nu]; \  \mathcal{A}_{k,\bs{\theta}}(\bs{\beta}) = \mathcal{A}_{j,\bs{\theta}}(\bs{\beta})} a_{k} p^{w_{j,k,\bs{\theta},\rho}(\bs{\beta})} \neq 0$ for some $j \in \mathcal{M}_{\nu,\bs{\beta}}$,  we first choose 
  $\bs{\phi} \in \mathbb{R}^n$ so that the scalar products $\bs{\phi}\cdot\mathcal{A}_{j, \bs{\theta}}(\bs{\beta})$ \ ($j \in \mathcal{M}_{\nu,\bs{\beta}}$)  are {\it pairwise distinct}.  
   Since the $\mathcal{A}_{j, \bs{\theta}}(\bs{\beta})$ \ (for $j \in \mathcal{M}_{\nu,\bs{\beta}}$) \ are distinct by definition,  elementary linear algebra shows that for each $\bs{\beta},$ there exists a finite  union $\mathcal {L}_\beta$ of {\it proper} linear subspaces such that if $\bs{\phi} \notin \mathcal{L}_\beta$    then these scalar products are  pairwise distinct. Thus, $\ds\mathcal{G}_{\bs{\tau}} := \mathbb{R}^n - \cup_{\bs{\beta}\notin B_e} \mathcal {L}_\beta$
is a generic subset of $\mathbb{R}^n$.  

Given $\bs{\phi} \in \mathcal{G}_{\bs{\tau}},$  we   set  $\boldsymbol{\tau}^0:=x \,\bs{\phi}$. It follows that if $x$ remains  outside a subset of $\mathbb{R}$ of empty interior, then the scalar products $x\,\bs{\phi}\cdot\mathcal{A}_{j,\bs{\theta}}(\bs{\beta})$ \ ($j \in \mathcal{M}_{\nu,\bs{\beta}}$) are pairwise distinct {\it for all $\bs{\beta}\notin B_e$}. We can then think of (\ref{exprgene}) as an identity between functions of $x$ that says the following for all $x$ outside this set of empty interior and all $\bs{\beta}\notin B_e$:
\begin{equation} \label{rewrite}
0 = \sum_{j \in \mathcal{M}_{\nu,\bs{\beta}}} m_j(\bs{\beta},\bs{\theta}, \rho)  p^{-i x\,\bs{\phi}\cdot\mathcal{A}_{j,\bs{\theta}}(\bs{\beta})}
\end{equation}
where 
$\ds m_j (\bs{\beta}, \bs{\theta}, \rho) = \bigg(\sum_{k \in [\nu]; \  \mathcal{A}_{k,\bs{\theta}}(\bs{\beta}) = \mathcal{A}_{j,\bs{\theta}}(\bs{\beta})} a_{k} p^{w_{j,k,\bs{\theta},\rho}(\bs{\beta})}\bigg) p^{-\mathcal{A}_{j,\bs{\theta}}'(\bs{\beta}, \rho)}.$

Denoting the elements of $\mathcal{M}_{\nu,\bs{\beta}}$ as $\{h_1 < h_2 < \dots < h_R\},$ \ where $R = \# \mathcal{M}_{\nu,\bs{\beta}},$ \ and differentiating (both sides of  (\ref{rewrite})) $R-1$ times with respect to $x$  we obtain a set of $R$ linear equations, which in matrix form is:  
\begin{equation} \label{matrixform}
\mathbf V_{\nu, \beta} \cdot \mathbf M = \mathbf 0 
\end{equation}
where 
$\ds {}^t \mathbf M = \left(m_{h_1}(\bs{\beta},\bs{\theta},\rho) p^{-i x\,\bs{\phi}\cdot\mathcal{A}_{h_1,\bs{\theta}}(\bs{\beta})},\dots, m_{h_R}(\bs{\beta},\bs{\theta},\rho) p^{-i x\,\bs{\phi}\cdot\mathcal{A}_{h_R,\bs{\theta}}(\bs{\beta})} \right)$ and
$\ds  \mathbf V_{\nu, \beta} = \left(v_{h_a,h_b}\right)_{a, b \in \{1,\dots, R\}}$ with $\ds v_{h_a, h_b} = \big(-ix\,\bs{\phi}\cdot\mathcal{A}_{h_b,\bs{\theta}}(\bs{\beta})\big)^{a-1}\,.$

By hypothesis, the second possibility implies that    $\mathbf M \neq \mathbf 0.$ In addition, we recognize $\mathbf V_{\nu, \beta}$ as a Vandermonde matrix
whose Vandermonde determinant does not equal $0$ {\it for all $\bs{\beta}\notin B_e$} precisely because
$\bs{\phi} \in \mathcal{G}_{\bs{\tau}},$ i.e. the scalar products $x\,\bs{\phi}\cdot\mathcal{A}_{h_b,\bs{\theta}}(\bs{\beta})$ are pairwise distinct. 
Thus, (\ref{matrixform}) {\it cannot} occur. So,   possibility (2)  is also impossible for all $\bs{\beta}\notin B_e$, provided $\boldsymbol{\tau}^0 \in \mathcal{G}_{\bs{\tau}}.$ With this contradiction, it suffices to set 
$\mathcal{G}' = \mathcal{G} \cap \mathcal{G}_{\bs{\tau}}$
to complete the proof of the above Claim.\qed 

{\it Finishing the proof of Lemma \ref{perturbation-sigma0}.}

 Having found both $p_0 (\mathcal {B})$ and the generic set $\mathcal{G}'$ during the proof of the Claim, we   now finish the proof of the lemma by  finding the generic set  $\mathcal{G}^*$ with the asserted property. 

We  assume
\begin{equation} \label{nucondition}
 \text{{\it $ \bs{\mu} = (p, \boldsymbol{\tau}^0, \bs{\theta})$ is a parameter vector
such that $p > p_0(\mathcal {B})$ and $\boldsymbol{\tau}^0 \in \mathcal{G}'.$ }}
\end{equation}

It is now convenient to  use the coordinates $\widetilde{\bs{\sigma}}^0$ on the hyperplane  $\{\bs{\sigma}^0\cdot \bs{\alpha}_{\cdot e} = 0\}$ that have also been used in the proof of Lemma \ref{aW-noncyclo} (Assertion 2).
We can also put $\widetilde{\bs{\alpha}}=\left(\widetilde{\alpha}_{\ell j}\right)_{(\ell,j)\in \{1,\dots,n-1\}\times\{1,\dots,r\}} \in \mathbb{M}_{n-1,r}(\mathbb{Z})$ such that for for $j\in \{1,\dots,r\}$ and $\ell \in \{1,\dots,n-1\}$  $\widetilde{\alpha}_{\ell j}=\alpha_{\ell j}-\frac{\alpha_{n j}}{\alpha_{n e}}\alpha_{\ell e}$;  so that we have for all $j\in \{1,\dots,r\}$ $\bs{\sigma}^0 \cdot  \bs{\alpha}_{\cdot j}=\widetilde{\bs{\sigma}}^0 \cdot \widetilde{\bs{\alpha}}_{\cdot j}$.
 In this way, the function $u_{\nu, \bs{\theta}}(\boldsymbol{\sigma}^0, \bs{\beta})$ appearing on the right side of (\ref{fonctionidentnulle}) then becomes a function $\tilde u_{\nu,\bs{\theta}}(\widetilde{\bs{\sigma}}^0, \bs{\beta})=  \widetilde{\bs{\sigma}}^0\cdot\widetilde{\bs{\alpha}}_{\cdot \nu} - \widetilde{\bs{\sigma}}^0\cdot\widetilde{\bs{\alpha}}\cdot {}^t \!\bs{\beta} \frac{\bs{\theta}\cdot\bs{\alpha}_{\cdot \nu} }{\bs{\theta}\cdot\bs{\alpha}\cdot {}^t \!\bs{\beta}}$ which is linear in $\widetilde{\sigma}^0$.  

For each $\bs{\beta}\notin B_e,$ we now choose  and fix $\tilde{\bs{\omega}} = \tilde{\bs{\omega}} (\bs{\beta}) \in \mathbb{R}^{n-1}$ to be a vector whose components are $\mathbb{Q}$-linearly independent and so that    $\tilde u_{\nu, \bs{\theta}}(\tilde{\bs{\omega}}, \bs{\beta})$ are {\it pairwise distinct} numbers  indexed by $\nu \in \mathcal{V}$.  We then set $
\widetilde{\bs{\sigma}}^0 = t \tilde{\bs{\omega}}.
$

Since $\widetilde{u}_{\nu,\bs{\theta}} (t \tilde{\bs{\omega}}, \bs{\beta}) = t \, \widetilde{u}_{\nu, \bs{\theta}}(\tilde{\bs{\omega}}, \bs{\beta})$, it follows that (\ref{fonctionidentnulle}) can then be rewritten as follows: $
\ds f_{\bs{\mu}, \bs{\beta}, \rho} \left(t \tilde{\bs{\omega}} \right) = 1+\summ_{\nu \in \mathcal{V}} A_{\nu, \bs{\mu},\rho}(\bs{\beta}) \exp\left(-t \log(p) \, \tilde{u}_{\nu, \bs{\theta}}(\tilde{\bs{\omega}}, \bs{\beta})\right),
$
where some $A_{\nu, \bs{\mu},\rho}(\bs{\beta}) \neq 0$ ({\it for each} $\bs{\beta}, \rho$) by the Claim. 
Since the $\tilde{u}_{\nu, \bs{\theta}}(\tilde{\bs{\omega}}, \bs{\beta})$  are  pairwise distinct,  the  $\# \mathcal{V}$  functions (defined for each $\bs{\mu}, \bs{\beta}, \rho$)
$\ds t \longmapsto \exp\left(-t\log(p) \,\tilde{u}_{\nu, \bs{\theta}}(\tilde{\bs{\omega}}, \bs{\beta})\right); \ (\nu \in \mathcal{V})$ 
are linearly independent.  It follows that  $t \to f_{\bs{\mu}, \bs{\beta}, \rho} \left(t \tilde{\bs{\omega}}\right)$ is not identically zero for each $(\bs{\beta}, \rho)$  and $\bs{\mu}$ satisfying (\ref{nucondition}). 

Thus, for each such $\bs{\mu},$  we conclude  that $f_{\bs{\mu}, \bs{\beta}, \rho} \left(\widetilde{ \boldsymbol{\sigma}}^0 \right)$ is not identically zero for each $(\bs{\beta}, \rho).$ An application of Lemma \ref{athmWeierstrass} then tells us that each hypersurface  $\{f_{\bs{\mu}, \bs{\beta}, \rho}(\widetilde{\bs{\sigma}}) = 0\}$ has empty interior inside  $\mathbb{R}^{n-1}$.

Define $
 \ds M_{\bs{\tau}^0,\bs{\theta}} = \bigcup_{\bs{\beta}, \rho, p>p_0(\mathcal{B})} f_{\bs{\mu}, \bs{\beta}, \rho}^{-1}(0).
 $
It follows from Baire's theorem that  each $M_{\bs{\tau}^0,\bs{\theta}}$ defines a set of  empty interior inside  $\mathbb{R}^{n-1}.$ We set
$\widetilde{\mathcal {G}}_{\boldsymbol{\tau}^0, \bs{\theta}} = \mathbb{R}^{n-1} - M_{\bs{\tau}^0,\bs{\theta}}.$
Via the  map $\widetilde{\boldsymbol{\sigma}}^0 \to \boldsymbol{\sigma}^0$ each $\widetilde {\mathcal {G}}_{\boldsymbol{\tau}^0, \bs{\theta}}$ defines a generic subset $\mathcal {G}_{\boldsymbol{\tau}^0, \bs{\theta}}$ of 
the hyperplane $\{\bs{\sigma}^0\cdot \bs{\alpha}_{\cdot e} = 0\}.$ 

We then finish the proof of Lemma \ref{perturbation-sigma0}   by setting (for any $\bs{\theta}$ satisfying (\ref{ahyp_theta1}) and (\ref{achoixjudicieuxtheta}))
$
\ds\mathcal{G}^{*} = \left\{\mathbf{s}=\boldsymbol{\sigma}^0+i\boldsymbol{\tau}^0 \mid  \boldsymbol{\tau}^0\in \mathcal{G}' \ \textrm{and} \ \boldsymbol{\sigma}^0 \in \widetilde {\mathcal {G}}_{\boldsymbol{\tau}^0, \bs{\theta}}\right\}.
$
\qed

\end{dem}
We now must address Case B to complete the proof of Theorem \ref{aresultatprincipal}.
Using the expression (\ref{fonctionidentnulle}), we first observe that there exist two a priori possibilities if Case B (i.e. $\boldsymbol{\beta}\in B_e$) occurs:
\begin {enumerate}
\item The property asserted by the Claim in the proof of Lemma \ref{perturbation-sigma0} is true if $\boldsymbol{\beta}\in B_e$ and $\rho$ satisfies $\Re (t(\bs{\beta},\rho)) > 0$;
\item The property asserted by the Claim in the proof of Lemma \ref{perturbation-sigma0} is false for some $\boldsymbol{\beta}\in B_e$ and $\rho$ satisfying $\Re (t(\bs{\beta}, \rho)) > 0.$
\end{enumerate}
If possibility 1 occurs, then we can repeat the proof of Lemma \ref{perturbation-sigma0} to complete the proof of Theorem \ref{aresultatprincipal}. Indeed, this would mean that the conclusion of Lemma \ref{perturbation-sigma0} applied to all $\bs{\beta}$. Thus, there would be infinitely many zeroes/poles of $W_{\mathbf{\bs{\mu}},\bs{\sigma}^0}(p^{-1}, p^{-t})$ in $\Xi_{u,\eta}$ that were not cancelled by the poles/zeroes $t(\bs{\beta}, \rho)$ {\it for any} $\bs{\beta}$,  not just $\bs{\beta}\notin B_e.$    
Since this would be the case for $\mathbf{\bs{\mu}}$ such that $(\boldsymbol{\sigma}^0, \boldsymbol{\tau}^0)$ is a generic vector, Theorem \ref{aresultatprincipal} would follow.

So, the difficulty occurs only when possibility 2 happens. Our argument in this case is quite different from that used to prove Lemma 
\ref{perturbation-sigma0}. Essentially, it  reduces to a counting argument.

For each  $\varepsilon \in ]0, 1[ $  we first  define  the region $\Xi_{u,\eta}^{\varepsilon}$  as follows:
\begin{center}
\begin{tabular}{ll}
$\Xi_{u,\eta}^{\varepsilon}:$ & $\Re\left(t\right) > \varepsilon $\\ 
 & $0<u<\Im\left(t\right)<u+\eta.$
\end{tabular}
\end{center}
In light of the preceding discussion, the proof of the next lemma will therefore finish the proof of Theorem  \ref{aresultatprincipal}. 

\begin{lemme}\label{lemmeBe}

Moving $\boldsymbol{\sigma}^0$ so that $\mathbf{s}^0\in \partial \mathbf{W}(0)\cap \mathcal{B}$ if necessary, there are inside $\Xi_{u,\eta}^{\delta}$ (as $\delta \longrightarrow 0$)  some zeroes  $t_{m,\bs{\mu},\bs{\sigma}^0}$ coming from $W_{\mathbf{\bs{\mu}},\bs{\sigma}^0}\left(p^{-1},p^{-t}\right)$ which are not poles of the $\zeta$-factors of $A_{M_{\delta}}$ corresponding to the $\boldsymbol{\beta}\in B_e$.
In particular there exists a infinite number of zeroes of $Z(\mathbf{s}^0+t \bs{\theta})$ inside $\Xi_{u,\eta}$ which are not cancelled.
\end{lemme}

\begin{dem} Let  $\bs{\mu}$ be a parameter vector as defined  in 
Remark \ref{parameternu}. Let $\boldsymbol{\beta}\in B_e$ and assume $\rho$ satisfies $\Re (t(\bs{\beta},\rho)) > 0.$ We assume $\bs{\mu}, \bs{\beta}, \rho$ are such that for each $\nu \in \mathcal {V}$ $ \sum_{k \in [\nu]} a_k p^{-v_{k, \bs{\theta}, \bs{\tau}^0,\rho}(\bs{\beta})} = 0\,.$ 
Applying the observation made in   Remark \ref{rqBe} to the expression for $v_{k,\bs{\mu}, \rho}(\bs{\beta})$, we see that $\boldsymbol{\beta}\in B_e$ implies  that the set $\mathcal {V}$ can be identified with a set of representatives for the equivalence relation $\mathcal{R}$:  
$$\bs{\alpha}_{\cdot k_1} \  \mathcal{R} \ \bs{\alpha}_{\cdot k_2} \quad \text{ if and only if } \quad \bs{\alpha}_{\cdot k_1} \in \bs{\alpha}_{\cdot k_2} + \langle \boldsymbol{\alpha}_{\cdot e} \rangle\,.$$ 
We  can therefore express the sum as follows:
$$
\summ_{k \in [\nu]} a_k p^{-v_{k, \bs{\theta}, \bs{\tau}^0,\rho}(\bs{\beta})} = \sum_{\{k : \bs{\alpha}_{\cdot k} -\bs{\alpha}_{\cdot \nu}\in \langle \boldsymbol{\alpha}_{\cdot e} \rangle \}} a_k p^{-\rho \, \frac{\bs{\theta}\cdot \bs{\alpha}_{\cdot k}}{\bs{\theta}\cdot\bs{\alpha}\cdot {}^t \!\bs{\beta}} - i\big(\bs{\tau}^0\cdot \bs{\alpha}_{\cdot k}-\bs{\theta}\cdot \bs{\alpha}_{\cdot k}\frac{\bs{\tau}^0\cdot\bs{\alpha}\cdot {}^t \!\bs{\beta}}{\bs{\theta}\cdot\bs{\alpha}\cdot {}^t \!\bs{\beta}}\big)}.
$$

Moreover, since $\boldsymbol{\beta}\in B_e,$ the expressions involving a product with
$\bs{\alpha}\cdot {}^t \!\bs{\beta}$ equal sums over $i \in \Lambda_e$. The  fact that $i\in\Lambda_e$ implies $\bs{\alpha}_{\cdot i} = q_i\boldsymbol{\alpha}_{\cdot e}.$ A simple calculation now shows  that

 $
\bs{\tau}^0\cdot \bs{\alpha}_{\cdot k}-\bs{\theta}\cdot \bs{\alpha}_{\cdot k}\frac{\bs{\tau}^0\cdot\bs{\alpha}\cdot {}^t \!\bs{\beta}}{\bs{\theta}\cdot\bs{\alpha}\cdot {}^t \!\bs{\beta}} = \bs{\tau}^0\cdot \bs{\alpha}_{\cdot k}-\bs{\theta}\cdot \bs{\alpha}_{\cdot k}\frac{\bs{\tau}^0 \cdot \bs{\alpha}_{\cdot e}}{\bs{\theta}\cdot \bs{\alpha}_{\cdot e}}.
$

As a result, $
0 = \summ_{k \in [\nu]} a_k p^{-v_{k, \bs{\theta}, \bs{\tau}^0,\rho}(\bs{\beta})} = \sum_{\{ k : \bs{\alpha}_{\cdot k} -\bs{\alpha}_{\cdot \nu}\in \langle \boldsymbol{\alpha}_{\cdot e} \rangle\}} a_k p^{-i \bs{\tau}^0\cdot \bs{\alpha}_{\cdot k}} p^{-\bs{\theta}\cdot \bs{\alpha}_{\cdot k} \big(\frac{\rho}{\bs{\theta}\cdot\bs{\alpha}\cdot {}^t \!\bs{\beta}} - i \frac{\bs{\tau}^0 \cdot \bs{\alpha}_{\cdot e}}{\bs{\theta}\cdot\bs{\alpha}_{\cdot e}}\big)}.
$ 
This now implies that {\it for each}  $\nu \in\mathcal{V}$, 
$p^{-\frac{\rho}{\bs{\theta}\cdot\bs{\alpha}\cdot {}^t \!\bs{\beta}} + i \frac{\bs{\tau}^0 \cdot \bs{\alpha}_{\cdot e}}{\bs{\theta}\cdot\bs{\alpha}_{\cdot e}}}$ is a root of the generalized polynomial $
\sum_{\{k: \bs{\alpha}_{\cdot k} -\bs{\alpha}_{\cdot \nu}\in \langle \boldsymbol{\alpha}_{\cdot e} \rangle \}} a_k p^{-i \bs{\tau}^0\cdot \bs{\alpha}_{\cdot k}} Y^{\bs{\theta}\cdot \bs{\alpha}_{\cdot k}}\,.
$

Applying the above interpretation of $\mathcal{V}$,   we further observe  that for any $\mathbf{\bs{\mu}} = (p,  \boldsymbol{\tau}^0,\bs{\theta})$ such that $\boldsymbol{\sigma}^0$ satisfies (\ref{agenericite}):  
\begin{displaymath}
\begin{array}{lll}
W_{\mathbf{\bs{\mu}},\bs{\sigma}^0}\left(X,Y\right) & = & \displaystyle \sum_{\nu \in\mathcal{V}} \ \sum_{\{k:   \bs{\alpha}_{\cdot k} -\bs{\alpha}_{\cdot \nu}\in \langle \boldsymbol{\alpha}_{\cdot e} \rangle\}} a_k p^{-i \bs{\tau}^0\cdot \bs{\alpha}_{\cdot k}} X^{\bs{\sigma}^0\cdot \bs{\alpha}_{\cdot k}} Y^{\bs{\theta}\cdot \bs{\alpha}_{\cdot k}} \\ 
 & = & \displaystyle \sum_{\nu \in\mathcal{V}} X^{\boldsymbol{\sigma}^0 \cdot \bs{\alpha}_{\cdot \nu}} \sum_{\{k : \bs{\alpha}_{\cdot k} -\bs{\alpha}_{\cdot \nu}\in \langle \boldsymbol{\alpha}_{\cdot e} \rangle \}} a_k p^{-i \bs{\tau}^0\cdot \bs{\alpha}_{\cdot k}} Y^{\bs{\theta}\cdot \bs{\alpha}_{\cdot k}}. 
\end{array} 
\end{displaymath}
Thus,   $W_{\mathbf{\bs{\mu}},\bs{\sigma}^0} \left(X, p^{-\frac{\rho}{\bs{\theta}\cdot\bs{\alpha}\cdot {}^t \!\bs{\beta}} + i \frac{\bs{\tau}^0 \cdot \bs{\alpha}_{\cdot e}}{\bs{\theta}\cdot \bs{\alpha}_{\cdot e}}}\right) = 0$  {\it  for any $X$}.

In addition, we observe that {\it this root is of modulus strictly less than } $1$. 

Indeed $
\left|p^{-\frac{\rho}{\bs{\theta}\cdot\bs{\alpha}\cdot {}^t \!\bs{\beta}}+i \frac{\bs{\tau}^0 \cdot \bs{\alpha}_{\cdot e}}{\bs{\theta}\cdot \bs{\alpha}_{\cdot e}}} \right| = \exp\left(-\frac{\Re(\rho)}{\bs{\theta}\cdot\bs{\alpha}\cdot {}^t \!\bs{\beta}}\log(p)\right)<1
$
because $\Re(\rho)>0$.  

Adapting  the discussion (and using the notations) from Part 1,    we conclude that there exists a Puiseux branch $\Omega_{\bs{\mu},\bs{\sigma}^0}$ of $W_{\mathbf{\bs{\mu}},\bs{\sigma}^0}(X,Y) = 0$ that reduces to the constant $ c_{\mathbf{\bs{\mu}},0}= c_{\mathbf{\bs{\mu}},0}(\bs{\beta},\rho) = p^{-\frac{\rho}{\bs{\theta}\cdot\bs{\alpha}\cdot {}^t \!\bs{\beta}}+i \frac{\bs{\tau}^0 \cdot \bs{\alpha}_{\cdot e}}{\bs{\theta}\cdot \bs{\alpha}_{\cdot e}}}$ whose norm is both independent of $p$ (by (\ref{adepenpcko})) {\it and}  strictly less than $1$.

Now, what interests us are the  zeroes $t_{m,\bs{\mu},\bs{\sigma}^0}$ of  the function $t \to W_{\mathbf{\bs{\mu}},\bs{\sigma}^0} (p^{-1},p^{-t})$  that are determined  by this branch $\Omega_{\bs{\mu},\bs{\sigma}^0}$, and which belong to  $\Xi^{\varepsilon}_{u,\eta}$. Such $t_{m,\bs{\mu},\bs{\sigma}^0}$ are necessarily of the form
$
t_{m,\bs{\mu},\bs{\sigma}^0} =-\frac{\log(c_{\mathbf{\bs{\mu}},0})}{\log(p)}+\frac{2 i \pi m}{\log(p)}\,, 
$
where $m\in\mathbb{Z}$ and $p$ is a prime number. 


\vspace{0.2cm}

To say that  $t_{m,\bs{\mu},\bs{\sigma}^0}\in \Xi^{\varepsilon}_{u,\eta}$ is to say that  $\Re(t_{m,\bs{\mu},\bs{\sigma}^0})>\varepsilon$ and $u<\Im(t_{m,\bs{\mu},\bs{\sigma}^0})<u+\eta$. Thus,
$$
p<\exp\left(-\frac{\log|c_{\mathbf{\bs{\mu}},0}|}{\varepsilon}\right) \quad  \textrm{and} \quad \frac{u \log(p)}{2 \pi} + \arg(c_{\mathbf{\bs{\mu}},0})<m<\frac{(u+\eta) \log(p)}{2 \pi} + \arg(c_{\mathbf{\bs{\mu}},0}).
$$

For fixed $\boldsymbol{\sigma}^0$ (satisfying (\ref{agenericite}),  $\boldsymbol{\tau}^0$ (belonging to $\mathcal{G}'$) and  $\bs{\theta}$ (satisfying (\ref{ahyp_theta1}) and (\ref{achoixjudicieuxtheta})), we can count the number of such $t_{m,\bs{\mu},\bs{\sigma}^0}$ (ignoring multiplicities) if they are  distinct.
Indeed, {\it the $t_{m,\bs{\mu},\bs{\sigma}^0}$ are pairwise distinct}. Suppose $t_{m,\bs{\mu},\bs{\sigma}^0} = t_{m',\bs{\mu}', \bs{\sigma}^0}$ with $\mathbf{\bs{\mu}}' = (p', \boldsymbol{\tau}^0, \bs{\theta})$. By taking the real parts, we get $
-\frac{\log|c_{\mathbf{\bs{\mu}},0}|}{\log(p)} = -\frac{\log|c_{\mathbf{\bs{\mu}}',0}|}{\log(p')}.$
 And since $|c_{\mathbf{\bs{\mu}},0}| = |c_{\mathbf{\bs{\mu}}',0}|$ (by (\ref{adepenpcko})), we obtain   $\log(p)=\log(p')$, i.e. $p = p'.$
Comparing the imaginary parts then gives  $m=m'$.

Consequently the number of such  $t_{m,\bs{\mu},\bs{\sigma}^0}$ (counted without  multiplicities) inside $\Xi^{\varepsilon}_{u,\eta}$ is given by 
$
 N_{(\boldsymbol{\sigma}^0,\boldsymbol{\tau}^0,\bs{\theta})}(\varepsilon) := \sum_{p<\exp\left(-\frac{\log|c_{\mathbf{\bs{\mu}},0}|}{\varepsilon}\right)} \left(\frac{\eta\log(p)}{2 \pi}+ \xi \right)
$
where $\xi = \xi(\boldsymbol{\sigma}^0, \boldsymbol{\tau}^0, \bs{\theta})$ and satisfies $|\xi|\leq 1$.

The prime number theorem tells us
$
\sum_{p\leq x}\log(p) \sim x \  (x \longrightarrow \infty)
$
and
$
\sum_{p\leq x}\xi = O(\pi(x)) = o(x) \  (x \longrightarrow \infty).
$
We conclude with the basic asymptotic:
\begin{equation} \label{countdelta}
N_{(\boldsymbol{\sigma}^0,\boldsymbol{\tau}^0, \bs{\theta})} (\varepsilon) \sim \frac{\eta}{2 \pi}\exp\left(-\frac{\log|c_{\mathbf{\bs{\mu}},0}|}{\varepsilon}\right) \ \ \ (\varepsilon \longrightarrow 0).
\end{equation}

The last step is to estimate the number of   possible poles/zeroes (counted without the multiplicities) coming from $t \to A_{M_{\delta}}\left(\boldsymbol{\sigma}^0 + i \boldsymbol{\tau}^0 + t \bs{\theta}\right)$  (see Theorem \ref{aana}).
If $t_0$ is a zero/pole of $t \to A_{M_{\delta}}\left(\boldsymbol{\sigma}^0 + i \boldsymbol{\tau}^0 + t \bs{\theta}\right)$ inside $\Xi^{\varepsilon}_{u,\eta}$, then there exists $\boldsymbol{\beta}\in \Nr$ such that $ \left(\boldsymbol{\sigma}^0 + i \boldsymbol{\tau}^0 +   t_0 \bs{\theta}\right)\cdot\bs{\alpha}\cdot {}^t \! \bs{\beta}$ is a zero or a pole of $\zeta\left(\cdot\right)$. 
This implies $
\Re\left(t_0\right)\bs{\theta}\cdot\bs{\alpha}\cdot {}^t \!\bs{\beta}  \leq \Re\left( \left(\boldsymbol{\sigma}^0 + i \boldsymbol{\tau}^0 +   t_0 \bs{\theta}\right)\cdot\bs{\alpha}\cdot {}^t \! \bs{\beta} \right) \leq 1.$
Thus, the bounds
 $\varepsilon < \Re\left(t_0\right) \leq \frac{1}{\bs{\theta}\cdot\bs{\alpha}\cdot {}^t \!\bs{\beta}}$  
and  
$
\bs{\theta}\cdot\bs{\alpha}\cdot {}^t \!\bs{\beta} < \frac{1}{\varepsilon}
$  follow. 
And since $\bs{\theta} \cdot \bs{\alpha}_{\cdot j}  > 0$ for each $j$, there exists a constant $c = c(\bs{\theta}) > 0$   such that
$
 \bs{\theta}\cdot\bs{\alpha}\cdot {}^t \!\bs{\beta} \geq c \, \Vert \bs{\beta}\Vert.
$
Hence 
\begin{equation}\label{***}
 \Vert \bs{\beta}\Vert = O(\varepsilon^{-1}),
\end{equation}
from which, the condition 
$ \Im\left(t_0\right) < u+\eta$ 
then implies
$\Im\left( \left(\boldsymbol{\sigma}^0 + i \boldsymbol{\tau}^0 +   t_0 \bs{\theta}\right)\cdot\bs{\alpha}\cdot {}^t \! \bs{\beta} \right)  
  = O_{u,\eta, \boldsymbol{\tau}^0}\left(\varepsilon^{-1}\right).$  

Having fixed $\eta>0$, the number of zeros or singularities of a $\zeta$-factor of $t\longmapsto A_{M_{\delta}}\left(\boldsymbol{\sigma}^0 + i \boldsymbol{\tau}^0 + t \bs{\theta}\right)$ is given by $
O\left(\frac{1}{\varepsilon}\log\left(\frac{1}{\varepsilon}\right)\right),
$
with regard to a classical result concerning the estimation of the number of nontrivial zeros of the Riemann zeta function having the imaginary part less than $\frac{1}{\delta}$.
Moreover, the same zero or singularity can, according to $\left(\ref{***}\right)$, appear in at most $\left(\frac{1}{\varepsilon}\right)^r$ terms; which gives at most $
O\left(\left(\frac{1}{\varepsilon}\right)^{r+1}\log\left(\frac{1}\varepsilon\right)\right)
$
zeros or singularities coming from $t\longmapsto A_{M_{\delta}}\left(\boldsymbol{\sigma}^0 + i \boldsymbol{\tau}^0 + t \bs{\theta}\right)$ inside $\Xi^{\varepsilon}_{u,\eta}$ (counted without their multiplicities).

Hence this estimation is negligible regardless to $N_{(\boldsymbol{\sigma}^0, \boldsymbol{\tau}^0, \bs{\theta})} (\varepsilon)$, which achieves the proof of this lemma and completes the proof of the main theorem  \ref{aresultatprincipal}.
\CQFD
\end{dem}

%
%

\section{Study on the possibility of a continuation of dimension strictly inferior beyond $\partial \mathbf{W}(0)$.}

Let us start with this example coming from N. Kurokawa in \cite{kurokawa}.
Put $
h(X_1,X_2,X_3) = 1- X_1 X_2 - X_2 X_3 - X_3 X_1 + 2 X_1 X_2 X_3.
$
One can easily check that $h$ is not cyclotomic and that the corresponding Euler product $
Z(s_1,s_2,s_3) = \prod_{p}\left(1-p^{-s_1 - s_2}-p^{-s_2-s_3}-p^{-s_3-s_1}+2 p^{-s_1-s_2-s_3} \right)$
continues to $
\mathbf{W}(0) = \{(s_1,s_2,s_3) \in \mathbb{C}^3 \mid \sigma_1+\sigma_2>0, \sigma_2+\sigma_3>0, \sigma_3+\sigma_1>0\}.
$
%
%
 According to the previous results, we know that there does not exist any meromorphic continuation to an open ball of complex dimension $3$ beyond any point of $\partial \mathbf{W}(0)$.

However $
Z(s_1,s_2,0) = \prod_{p}\left(1-p^{-s_1}\right)\left(1-p^{-s_2}\right) = \frac{1}{\zeta(s_1)\zeta(s_2)}
$
 is meromorphic on $\mathbb{C}^2$. So here there is a continuation on a complex hypersurface beyond the point $\textbf{0} \in \partial \mathbf{W}(0)$.
This example shows that Theorem \ref{aresultatprincipal} is optimal from the point of view of the complex dimension of a possible meromorphic extension beyond $\partial \mathbf{W}(0)$.

\subsection{On the existence of a continuation on a real hypersurface beyond $\partial \mathbf{W}(0)$.}

We have just seen that the natural boundary, when existing, makes sense only for meromorphic continuations of maximal complex dimension.
Thus we cannot expect to make  sense to the natural boundary in a general way if we come down from a complex dimension.

However, we can wonder if it is possible to improve the previous results only by coming down from one real dimension; in other words if there can  exist or not a continuation on a real hypersurface beyond $\partial \mathbf{W}(0)$ in a sense needing to be precised since in this case the notion of holomorphy does not a priori make sense.
The answer is given in Theorem 3.
For that we need to appeal to the theory of C-R functions (Cauchy-Riemann) on a real hypersurface which generalises the class of holomorphic functions.
To begin let us recall the following classical result:

\begin{lemme}[\cite{thiebaut}]\label{S-partie}
 Let $D$ be a connected open ball of $\mathbb{C}^n$ and $f$ and $g$ be analytic functions on $D$.
If $f$ coincides with $g$ on a part $S$ of $D$ for which there exists a connected open $V$ of $D$ such that $V\setminus S$ is not connected, then  $f$ coincides with $g$ on $D$.
In particular, this result is true if $S$ is a real hypersurface of $D$.
\end{lemme}

The reader could refer to \cite{thiebaut} page 15 for a proof.
\vspace{0.2cm}

Now the key point is a fondamental result concerning C-R functions of which we will find a proof in \cite{range}:

\begin{lemme}\label{ext-CR}
 Let $\frak{H}$ be a real-analytic   hypersurface in an open of $\mathbb{C}^n$, and $f:\frak{H}\longrightarrow \mathbb{C}$ be a real-analytic   C-R function.
Then, if $v \in \frak{H}$, there exists a neighborhood $U$ of $v$ and $F$ an holomorphic function on $U$ such that $F=f$ on $U\cap \frak{H}$.
\end{lemme}

Finally, the proof of Theorem 3 follows directly from the two previous lemmas:

\vspace{0.2cm}
%

\noindent{\bf Proof of Theorem 3:}
\vspace{0.2cm}

 Consider $\frak{H}$ to be a real-analytic hypersurface which intersects across $\mathcal{F}(\boldsymbol{\alpha}_{\cdot e})$ and assume by absurd that $f$ is a continuation of  $Z(\mathbf{s})$ on $\frak{H}$ in a neighborhood of a point $\mathbf{s}^0 \in \frak{H}\cap\mathcal{F}(\boldsymbol{\alpha}_{\cdot e})$.
Moreover put $S=\frak{H}\cap \mathbf{W}(0)$.
Then, thanks to lemma \ref{ext-CR}, there exists a neighborhood $U\subset\mathbb{C}^n$ of $\mathbf{s}^0$ and $F$ an holomorphic function on $U$ such that $F=f$ on $U\cap \frak{H}$.
But since $f$ is an extension of $Z(\mathbf{s})$, we also have $Z=f=F$ on $S$.
According to lemma \ref{S-partie}, we have  $Z=F$ on $U\cap \mathbf{W}(0)\neq \emptyset$.
But that means that there exists an open ball $\mathcal{B}\subset U\subset\mathbb{C}^n$ centered in $\mathbf{s}^0 \in \mathcal{F}(\boldsymbol{\alpha}_{\cdot e})$ such that  $F$ extends $Z$ to $\mathcal{B}$; which is impossible in accordance with Theorem \ref{aresultatprincipal}. 
\CQFD

\end{document}